\documentclass[12pt]{amsart}
\input{diagrams}

\newcommand{\nc}{\newcommand}

\renewcommand{\top}{\mathrm{top}}

\nc{\one}{\mbox{\bf 1}}

\nc{\invtensor}{\underset{\leftarrow}{\otimes}}

\nc{\ad}{\operatorname{ad}}

\nc{\tr}{\operatorname{tr}}
\nc{\str}{\operatorname{str}}

\nc{\rank}{\operatorname{rank}}
\nc{\rnk}{\operatorname{rank}}

\nc{\corank}{\operatorname{corank}}

\nc{\Sym}{\operatorname{Sym}}

\nc{\sym}{\operatorname{sym}}

\nc{\id}{\operatorname{id}}
\nc{\Id}{\operatorname{Id}}

\nc{\htt}{\operatorname{ht}}

\nc{\Norm}{\operatorname{Norm}}
\nc{\Ker}{\operatorname{Ker}}
\nc{\rker}{\operatorname{rKer}}

\nc{\im}{\operatorname{Im}}

\nc{\osp}{\operatorname{osp}}
\nc{\sgn}{\operatorname{sgn}}
\nc{\F}{\operatorname{F}}
\nc{\Mod}{\operatorname{Mod}}
\nc{\Mat}{\operatorname{Mat}}
\nc{\sMat}{\operatorname{sMat}}
\nc{\Soc}{\operatorname{Soc}}
\nc{\Inj}{\operatorname{Inj}}
\nc{\Hom}{\operatorname{Hom}}
\nc{\End}{\operatorname{End}}
\nc{\supp}{\operatorname{supp}}
\nc{\Card}{\operatorname{Card}}
\nc{\Ann}{\operatorname{Ann}}
\nc{\Ind}{\operatorname{Ind}}
\nc{\Coind}{\operatorname{Coind}}
\nc{\wt}{\operatorname{wt}}
\nc{\spn}{\operatorname{span}}
\nc{\ch}{\operatorname{ch}}
\nc{\codim}{\operatorname{codim}}
\nc{\Stab}{\operatorname{Stab}}

\nc{\Sch}{{\mathcal S}\mbox{\em ch}}

\nc{\Irr}{\operatorname{Irr}}

\nc{\Spec}{\operatorname{Spec}}

\nc{\Prim}{\operatorname{Prim}}

\nc{\Aut}{\operatorname{Aut}}

\nc{\Fract}{\operatorname{Fract}}

\nc{\gr}{\operatorname{gr}}

\nc{\HC}{\operatorname{HC}}

\nc{\wdchi}{\widetilde{\chi}}

\nc{\wdH}{\widetilde{H}}

\nc{\wdN}{\widetilde{N}}

\nc{\wdM}{\widetilde{M}}

\nc{\wdO}{\widetilde{O}}

\nc{\wdR}{\widetilde{R}}

\nc{\wdS}{\widetilde{S}}

\nc{\wdV}{\widetilde{V}}

\nc{\wdC}{\widetilde{C}}

\nc{\cN}{\operatorname{\mathcal M}^{\#}}
\nc{\cM}{\operatorname{\mathcal M}}
\nc{\Ob}{\operatorname{\mathcal Ob}}

\nc{\cO}{\operatorname{\mathcal O}}
\nc{\cC}{\operatorname{\mathcal C}}
\nc{\cD}{\operatorname{\mathcal D}}
\nc{\Dglie}{\operatorname{{\mathcal D}glie}}

\nc{\Fin}{\operatorname{{\mathcal F}in}}

\nc{\Sg}{{\mathcal S}({\mathfrak g})}

\nc{\Ug}{{\mathcal U}({\mathfrak g})}

\nc{\Zg}{{\mathcal Z}({\mathfrak g})}

\nc{\tZg}{{\widetilde{\mathcal Z}({\mathfrak g})}}

\nc{\Zk}{{\mathcal Z}({\mathfrak k})}

\nc{\Sh}{{\mathcal S}({\mathfrak h}_0)}

\nc{\Uh}{{\mathcal U}({\mathfrak h})}

\nc{\Up}{{\mathcal U}({\mathfrak p})}

\nc{\Ub}{{\mathcal U}({\mathfrak b})}

\nc{\Zh}{{\mathcal Z}({\mathfrak h})}
\nc{\Ah}{{\mathcal A}({\mathfrak h})}

\nc{\Ag}{{\mathcal A}({\mathfrak g})}

\nc{\Ap}{{\mathcal A}({\mathfrak p})}

\nc{\Zp}{{\mathcal Z}({\mathfrak p})}

\nc{\cZ}{\mathcal Z}

\nc{\cS}{\mathcal S}

\nc{\cB}{\mathcal B}
\nc{\cP}{\mathcal P}

\nc{\cA}{\mathcal A}

\nc{\cU}{\mathcal U}

\nc{\cH}{\mathcal H}

\nc{\cL}{\mathcal L}

\nc{\cF}{\mathcal F}

\nc{\fg}{\mathfrak g}

\nc{\CO}{\mathcal O}

\nc{\Cl}{\mathcal {C}\ell}

\nc{\fn}{\mathfrak n}

\nc{\fm}{\mathfrak m}

\nc{\fp}{\mathfrak p}

\nc{\fh}{\mathfrak h}

\nc{\ft}{\mathfrak t}

\nc{\fk}{\mathfrak k}

\nc{\fb}{\mathfrak b}

\nc{\fI}{\mathfrak I}

\nc{\veps}{\varepsilon}
\nc{\vareps}{\varepsilon}

\nc{\fs}{\mathfrak s}

\nc{\fsl}{\mathfrak{sl}}
\nc{\fgl}{\mathfrak{gl}}

\nc{\fpq}{\mathfrak{pq}}
\nc{\fq}{\mathfrak q}
\nc{\fsq}{\mathfrak{sq}}
\nc{\fpsq}{\mathfrak{psq}}

\nc{\fpo}{\mathfrak{po}}

\nc{\dirlim}{\underset{\rightarrow}{\lim}\,}

\nc{\nen}{\newenvironment}

\nc{\ol}{\overline}

\nc{\ul}{\underline}

\nc{\ra}{\rightarrow}

\nc{\lra}{\longrightarrow}

\nc{\Lra}{\Longrightarrow}

\nc{\Lla}{\Longleftarrow}

\nc{\Llra}{\Longleftrightarrow}

\nc{\thla}{\twoheadleftarrow}

\nc{\hra}{\hookrightarrow}

\nc{\iso}{\overset{\sim}{\lra}}

\nc{\ssubset}{\underset{\not=}{\subset}}

\nc{\Thm}[1]{Theorem~\ref{#1}}

\nc{\Prop}[1]{Proposition~\ref{#1}}

\nc{\Lem}[1]{Lemma~\ref{#1}}

\nc{\Cor}[1]{Corollary~\ref{#1}}

\nc{\Conj}[1]{Conjecture~\ref{#1}}

\nc{\Claim}[1]{Claim~\ref{#1}}

\nc{\Defn}[1]{Definition~\ref{#1}}

\nc{\Exa}[1]{Example~\ref{#1}}

\nc{\Rem}[1]{Remark~\ref{#1}}

\nc{\Note}[1]{Note~\ref{#1}}

\nc{\Quest}[1]{Question~\ref{#1}}

\nc{\Hyp}[1]{Hypoth\`ese~\ref{#1}}

\nen{thm}[1]{\label{#1}{\bf Theorem.\ } \em}{}

\nen{prop}[1]{\label{#1}{\bf Proposition.\ } \em}{}

\nen{lem}[1]{\label{#1}{\bf Lemma.\ } \em}{}

\nen{cor}[1]{\label{#1}{\bf Corollary.\ } \em}{}

\nen{conj}[1]{\label{#1}{\bf Conjecture.\ } \em}{}

\nen{claim}[1]{\label{#1}{\bf Claim.\ } \em}{}

\nen{defn}[1]{\label{#1}{\bf Definition.\ } }{}

\nen{exa}[1]{\label{#1}{\bf Example.\ } }{}

\nen{rem}[1]{\label{#1}{\em Remark.\ } }{}

\nen{note}[1]{\label{#1}{\em Note.\ } }{}

\nen{exer}[1]{\label{#1}{\em Exercise.\ } }{}

\nen{sket}[1]{\label{#1}{\em Sketch of proof.\ } }{}

\nen{quest}[1]{\label{#1}{\bf Question.\ } \em}{}

\nen{hyp}[1]{\label{#1}{\bf Hypoth\`ese.\ } \em}{}

\begin{document}

\title[Shapovalov determinants of $Q$-type Lie superalgebras]
{Shapovalov determinants of $Q$-type Lie superalgebras}
\author[Maria Gorelik]{Maria Gorelik}
%Incumbent of the Frances and Max Hersh career development chair,
\address{
Dept. of Mathematics, The Weizmann Institute of Science,
Rehovot 76100, Israel
{\tt email: maria.gorelik@weizmann.ac.il} 
}

%\thanks{}

\begin{abstract}
We define an analogue of Shapovalov forms for $Q$-type Lie superalgebras
and factorize the corresponding Shapovalov determinants which are responsible
for simplicity of highest weight modules. We apply the factorization
to obtain a description of the centres of $Q$-type Lie superalgebras.
\end{abstract}

\maketitle


2000 Mathematics Subject Classification 17B20, 17B35.

\newpage

\section{Introduction}

\subsection{}
\label{intro1}
In 1972 N.~Shapovalov~\cite{sh} suggested a powerful method
for studying highest weight modules of a finite
dimensional simple Lie algebra. He elucidated the description
of a bilinear form on the enveloping algebra of
a simple finite dimensional complex Lie algebra $\fg$ introduced
by Gelfand and Kirillov in~\cite{GK}. 
The kernel of this form (Shapovalov form) at
 a given point $\lambda\in\fh^*$ determines the maximal
submodule $\ol{M(\lambda)}$ of a Verma module $M(\lambda)$.
In particular, a Verma module $M(\lambda)$ is simple
if and only if the kernel of Shapovalov form at $\lambda$ is equal to zero.
The Shapovalov form can be realized as a direct sum of forms
$S_{\nu}$; for each $S_{\nu}$ one can define its determinant 
(Shapovalov determinant). The zeroes of Shapovalov determinants determine
when a Verma module is reducible. 
N.~Shapovalov computed
these determinants for the finite
dimensional simple Lie algebras: he presented them as products
of polynomials of degree one. As a consequence, a Verma module
$M(\lambda)$ is simple if and only if $\lambda$ does not belong to a
union of hyperplanes.

Shapovalov's method was generalized by V.~Kac, D.~Kazhdan 
in~\cite{kk} to Kac-Moody Lie algebras with symmetrizable Cartan matrix,
by  V.~Kac (\cite{KFF}, \cite{k13}) to Lie superalgebras 
with symmetrizable Cartan matrix, 
and by C.~De Concini, V.~Kac (\cite{dk}) and A.~Joseph (\cite{jbook})
 to quantum case.
The formula  for Shapovalov determinants for 
 Lie superalgebras with symmetrizable
Cartan matrix is given in~\cite{k13}.

\subsection{}
By the term ``$Q$-type superalgebras'' we mean four series of
Lie superalgebras: $\fq(n)$ ($n\geq 2$) and its subquotients 
$\fsq(n)$, $\fpq(n)$, $\fpsq(n)$ (the last one is a simple 
Lie superalgebra for $n\geq 3$ in the notation of~\cite{kadv} it is $Q(n)$).
The $Q$-type Lie superalgebras are rather special. First, their
Cartan subalgebras are not abelian and have  non-trivial odd components.
Second, they  possess a non-degenerate invariant
bilinear form which is {\em odd} and they do not have
quadratic Casimir elements.

\subsubsection{}\label{NlMl}
The first peculiarity leads to the existence of two different
candidates for a role of Verma module of the highest weight 
$\lambda\in\fh^*_{\ol{0}}$:
a module $M(\lambda)$ which is induced from a simple $\fh_{\ol{0}}$-module 
$\mathbb{C}_{\lambda}$
and a module $N(\lambda)$ which is induced from a simple $\fh$-module.
The character of $M(\lambda)$ nicely
depends on $\lambda$; we call $M(\lambda)$ a {\em Verma module}.
We call  $N(\lambda)$ a {\em Weyl module}. Observe that
each Verma module $M(\lambda)$ has a finite filtration
with the factors isomorphic to $N(\lambda)$ up to a parity change.
Each Weyl module has a unique simple quotient.

\subsubsection{}
In this paper we define a Shapovalov map for $Q$-type superalgebras.
Its kernel at a given point $\lambda\in\fh^*_{\ol{0}}$ determine the maximal
submodule $\ol{M(\lambda)}$ which does not meet the highest weight space.
The above observation implies that the Weyl module $N(\lambda)$
is simple if and only if $\ol{M(\lambda)}=0$.

\subsubsection{}
It turns out that the Shapovalov determinants again admit 
linear factorization (i.e., are the products
of polynomials of degree one) and so a Weyl module
$N(\lambda)$ is simple if and only if $\lambda$ does not belong to a
union of hyperplanes.

\subsubsection{}
In all cases mentioned in~\ref{intro1}
 the calculation of Shapovalov determinants
uses an explicit formula for a quadratic Casimir element
which implies a linear factorization for Shapovalov determinants.

In the present work the calculation
is based on an observation which allows one to deduce
the linear factorizability of the Shapovalov determinants without
using the quadratic Casimir elements--- see~\ref{introcal}.

\subsubsection{Determinants versus reduced norm.}
Let $\fg$ be a Lie superalgebra of $Q$-type and 
$\fh=\fh_{\ol{0}}\oplus\fh_{\ol{1}}$
be its Cartan subalgebra.
Each Shapovalov map $S_{\nu}$ is a map between two 
bimodules over the non-commutative
algebra $R:=\Uh$. As left and as right $R$-modules the source and the target
of $S_{\nu}$ are free of the same finite rank.
Viewing $S_{\nu}$ as an $\Sh$-homomorphism between
free $\Sh$-modules we define $\det S_{\nu}\in \Sh$.
Similarly to the case of endomorphisms
of modules over an Azumaya algebra, $\det S_{\nu}\in \Sh$ 
turns out to be a power
of another polynomial $\Norm S_{\nu}$ (reduced norm) which we propose
as an analogue of Shapovalov determinants for $Q$-type
superalgebras. Notice that the resulting formulas for $\Norm S_{\nu}$
look like the formulas for Shapovalov determinants for contragredient Lie
superalgebras (see~\cite{k13}).
We leave to Appendix a thorough
definition of reduced norm which would cover our setup.

\subsection{Computation of $\det S_{\nu}$.}
\label{introcal}
The computation of Shapovalov determinants in~\cite{kk},\cite{jbook} has the 
following steps. 
The first one is to show that each determinant admits a linear factorization;
this easily follows from the existence of a quadratic Casimir.
The second step is to construct the Jantzen filtration on Verma modules
which provides some information about the multiplicity of each linear factor.
Finally, one computes the leading term of $\det S_{\nu}$ and
then obtains  the multiplicities.

\subsubsection{Linear factorizability}\label{linfactint}
Let $\fg$ be a classical Lie superalgebra which is not of type $P$.

Denote by $W$ the Weyl group of $\fg_{\ol{0}}$ and by $\Zg$ the centre of  the
universal enveloping algebra $\Ug$.
A Harish-Chandra projection identifies $\Zg$
with a subalgebra $Z$ of $W$-invariant polynomials on $\fh_{\ol{0}}^*$;
if $\fg$ is a semisimple Lie algebra one has $Z=\Sh^W$. 
For any $\fg$ there exists a non-zero homogeneous polynomial
$z_a\in Z$ such that
$$z_a\Sh^W\subset Z. \ \eqno{(**)}$$
The property (**) could be easily deduced (see~\cite{pe})
from the explicit description of $Z$ given in~\cite{s1},
\cite{klapl},\cite{sq}. One can also obtain (**)
by describing the anticentre of $\Ug$ which is much easier
to describe than $Z$ itself (see~\cite{ghost}) and then taking
$z_a:=T^2$ where $T$ is an anticentral element 
of the minimal degree.

Now  the linear factorizability of Shapovalov determinants can be 
achieved as follows. Let
$C\in\Sh^W$ be the standard quadratic element
given by $C(\lambda)=(\lambda,\lambda)$. By (**) $Z$ contains $z_aC$.
If $\det S_{\mu}(\lambda)=0$ for some $\mu$, then a Verma module $M(\lambda)$
has a primitive vector of weight $\lambda-\nu$ for some $0<\nu\leq \mu$ 
that is
$(z_aC)(\lambda)=(zC)(\lambda-\nu)$ (here $\leq$ is the standard partial
order on $\fh_{\ol{0}}^*$). If $z_a(\lambda)\not=0$
this implies $(\lambda,\lambda)=(\lambda-\nu,\lambda-\nu)$ that
is $2(\lambda,\nu)=(\nu,\nu)$. 
Let $Q^+$ be the positive part of root lattice $Q(\pi)$.
For each $\nu\in Q^+$ the last equation
defines a hyperplane.
Hence $\det S_{\mu}(\lambda)=0$ 
implies $2(\lambda,\nu)=(\nu,\nu)$ for some $\nu\in Q^+$ or 
$z_a(\lambda)=0$. This means that $\det S_{\mu}$ admits a factorization
where each factor is either linear or one of irreducible factors
of $z_a$. Notice that instead of Casimir $C$ we could use all elements
of $\Sh^W$. This would reduce the above set of hyperplanes
to those corresponding to $\nu=n\alpha$ where $\alpha$ is a positive root.

The computation of the leading term of $\det S_{\nu}$
shows that it admits a linear factorization.
Hence an irreducible homogeneous polynomial of a higher degree can not
be a factor of $\det S_{\nu}$. Since $z_a$ is homogeneous, this implies
that $\det S_{\nu}$ admits a linear factorization.

\subsubsection{}
To find the multiplicities of linear factors
we compute the leading term of $\det S_{\nu}$.
Then we
define a Jantzen-type filtration on a Verma  module
$M(\lambda)$ and prove a sum formula for the multiplicities.
Comparing the leading term and the sum formula we
determine the multiplicities.

The leading term can be computed by various methods.
In this text we use a reduction to the minimal rank case
which is $\fsq(2)$ for $Q$-type algebras.

\subsection{Applications}
\label{introcent}
The computation of Shapovalov determinants gives us immediately
a criterion of irreducibility of a Weyl module.

We also obtain a description of  Jantzen filtration at the points 
$\lambda\in\fh_{\ol{0}}^*$ corresponding to a generic reducible Weyl module.
This is essential for the computation of $\Zg$, see below.

The Harish-Chandra projection $\HC$ provides an embedding
$\Zg\to\Sh$.
Explicit knowledge of the Shapovalov determinants allows us to describe
the image of this embedding following a Kac approach ~\cite{klapl};
we give some details in~\ref{km1},\ref{km2}.
The result is 
similar to the one for contragredient Lie superalgebras (see~\cite{klapl}).
For the case $\fg=\fq(n)$ the centre $\Zg$ was described 
in~\cite{sq}, \cite{ns}.

\subsubsection{}\label{km1}
The main idea is to recover a central element by its action on Verma 
modules. 

In Section~\ref{secthat} we introduce a certain completion $\hat{U}$
of $\Ug$; roughly speaking,
$\hat{U}$ is an algebra which acts on all $\fg$ modules
which are locally nilpotent over $\fn^+$.
 We show that the centre $\cZ(\hat{U})$
coincides with $\Zg$. This follows from the following statement
suggested to the author by J.~Bernstein: 
 for any $a\in \cZ(\hat{U})$ one has $\deg\HC(a)=\deg(a)$.
This statement can be viewed as an analogue
of Chevalley's theorem stating  that
for a semisimple Lie algebra the restriction of a non-zero $\fg$-invariant
regular function on $\fg$ to $\fh$ is non-zero.
The formula $\deg\HC(a)=\deg(a)$ implies that $\cZ(\hat{U})=\Zg$
if $\fg$ is a finite-dimensional contragredient 
or $Q$-type Lie superalgebra. 

In~\ref{IR1} we define a 
$\fg$--$\fh$ bimodule
$\cM$ which plays role of a generic Verma module (so that
a Verma module $M(\lambda)$
 can be viewed as the evaluation of $\cM$ at $\lambda$).
The $\fg$-action on $\cM$ can be extended to an
action of $\hat{U}$. It turns out that $\cM$ is a
faithful $\hat{U}$-module; moreover,
$z\in\hat{U}$ is central iff the action of $z$
 on $\cM$ coincides with the
right action of $\HC(z)$: $zv=v\HC(z)$ for all $v\in\cM$.
Using this property, we compute the centre $\cZ(\hat{U})$, 
see~\ref{km2} for details.

\subsubsection{}\label{km2}
Knowledge of submodules of a generic reducible Weyl module
gives us necessary conditions on $\HC(z)$ for $z\in\cZ(\hat{U})$.
Then for each $\phi\in\Sh$ satisfying these necessary conditions
we construct an element $z=\sum z_{\nu}\in\cZ(\hat{U})$ 
with $\HC(z)=\phi$ by a recursive procedure introduced in~\cite{klapl}.
The key ingredient is that $S_{\nu}$ is invertible
over the field of fractions of $\Sh$ and
that $S_{\nu}^{-1}$ has poles of order at most one
at a subset of codimension two in $\fh^*_{\ol{0}}$.
The fact that $S_{\nu}^{-1}$ has poles of order at most one
at $\lambda$ is equivalent to the statement that
the Jantzen filtration of $M(\lambda)$ has length 
at most two; the latter holds for the {\em regular} and {\em subregular} 
points $\lambda\in\fh^*_{\ol{0}}$.

In~\cite{gcent} we have checked
the finiteness of the recursive procedure by an estimation
of degrees and thus show that the central element $z\in\cZ(\hat{U})$
lies in $\Zg$. In the present paper we use the equality
 $\cZ(\hat{U})=\Zg$ which is proven independently.

\subsection{Construction of Shapovalov maps.}
\label{secintro1}
The Shapovalov forms can be naturally interpreted in terms
of {\em Shapovalov maps} which we define below.
This approach was suggested to us by J.~Bernstein.

In~\cite{kadv} V.~Kac introduced a notion of contragredient
Lie superalgebra. These are Lie superalgebras
which can be constructed by a standard procedure from their Cartan
matrices. Let $\fg=\fn^-+\fh+\fn^+$ be a contragredient or a $Q$-type
Lie superalgebra. Denote by $\sigma$
the antiautomorphism of $\Ug$ equal to $-\id$ on $\fg$.

Let $\cC$ be the category of $\fh$-modules 
and $\cD$ be the category of $\fg$-modules graded by elements of $Q^-:=-Q^+$
(the grading is consistent with the natural $Q(\pi)$-grading on $\Ug$).
We denote by $N_{\nu}$ the $\nu$th homogeneous component of $N$.
Let $\Phi_0: \cD\to\cC$ be the functor given by $\Phi_0(N)=N_0$.
The functor $\Phi_0$ admits a left adjoint functor $\Ind:\cC\to\cD$
and  a right adjoint functor $\Coind$. For any $L\in\cC$ 
the adjunction morphisms $L\to \Phi_0(\Ind(L))$ and
$\Phi_0(\Coind(L))\to L$ are, in fact, isomorphisms. In particular,
$$\Hom_{\cD}(\Ind(L),\Coind(L))=\Hom_{\cC}(L,L).$$
Let $\Xi(L): \Ind(L)\to\Coind(L)$ correspond to
the identity map $L\to L$; in this way we obtain  
a morphism of functors $\Xi:\Ind\to\Coind$. 
The kernel of $\Xi(L)$ is the maximal graded
submodule of $\Ind(L)$ which does not meet its zero component.

\subsubsection{}\label{IR}
Set $R:=\Uh, A:=\Sh$ ($A=R$ if $\fg$  is not of $Q$-type).
View $R$ as an object of $\cC$.
We check that the canonical morphisms
$\Ind(E)\to\Ind(R)\otimes_{R}E,\ \ \ \Coind(R)\otimes_R E\to
\Coind(E)=\Coind(R)\otimes_R E$ are isomorphisms; they allow to identify
$\Xi(E)$ with $\Xi(R)\otimes \id_E.$
We call $S:=\Xi(R)$ a {\em Shapovalov map}. 

Both ${\fg}$--$\fh$ bimodules $\Ind(R),\Coind(R)$ viewed as $\fh_{\ol{0}}$-bimodules 
are isomorphic to $\cU(\fb^-)$ (where $\fb^-:=\fh+\fn^-$). 
Their homogeneous components are $R$-bimodules which are
free $A$-modules of the same finite rank.
Thus we can decompose $S=\sum S_{\nu}$ where 
$$S_{\nu}:\Ind(R)_{-\nu}\to\Coind(R)_{-\nu}$$ 
is an $R$-bimodule homomorphism. Viewing the source and the target
as left $A$-modules we realize
$S_{\nu}$ as an $A$-homomorphism 
between two free $A$-modules of the same rank. 
A matrix of $S_{\nu}$ (with entries in $A$)
is called a {\em Shapovalov matrix} and its determinant
is called a {\em Shapovalov determinant} (this is an element in $A=\Sh$
which is defined up to an invertible scalar).

\subsubsection{}\label{IR1}
Define $M(\lambda)$ and $N(\lambda)$ as in~\ref{NlMl}
($M(\lambda)=N(\lambda)$ if $\fg$ is not of $Q$-type).
The family of $M(\lambda)$ can be obtained from $\cM:=\Ind(R)$
by evaluation.
Denote by $S(\lambda)$ the evaluation
of the Shapovalov map $S$ at $\lambda$. The kernel of
$S(\lambda)$ is the maximal submodule of $M(\lambda)$ which does not 
meet the highest weight space. We denote this
submodule  by $\ol{M(\lambda)}$ and define similarly $\ol{N(\lambda)}$.

\subsubsection{}
Recall that $M(\lambda)$ (resp., $\ol{M(\lambda)}$)
has a finite filtration
whose factors are isomorphic to $N(\lambda)$ (resp., to $\ol{N(\lambda)}$)
up to parity change.
As a consequence, 
$$\Ker S(\lambda)=0\ \Longleftrightarrow\ 
\ol{M(\lambda)}=0\ \Longleftrightarrow\ \ol{N(\lambda)}=0.$$

Each $N(\lambda)$ has a unique simple quotient $V(\lambda)$
and $V(\lambda)=N(\lambda)/\ol{N(\lambda)}$. Hence
$$\begin{array}{c}
V(\lambda)=N(\lambda)
\ \Longleftrightarrow\ \det S_{\nu}(\lambda)\not=0
\text{ for all }\nu.
\end{array}$$

\subsection{Shapovalov forms.}\label{secintro2}
Historically, the Shapovalov map was introduced as a bilinear form.
This can be described as follows. 
The module $\Ind(R)$ identifies with $\cU(\fb^-)$ as 
$\fb^-$-$R$ bimodule.
The module $\Coind(R)$ can be realized
(up to a parity change) as a graded dual of $\cU(\fb^+)$.
More precisely,  $\Coind(R)$ identifies with the maximal graded submodule of
$\Hom_R(\cU(\fb^+), R^{\sigma})$ where $R^{\sigma}$
is an $R$ bimodule obtained from $R$ by the shift by $\sigma$.
Using this identification, we realize the Shapovalov map
as $S: \cU(\fb^-)\to\Hom_{R_r}(\cU(\fb^+),R^{\sigma})$;
the formula for $S$ is
$$S(u_-)(u_+)=(-1)^{p(u_-)p(u_+)}\HC\bigl(\sigma(u_+)u_-\bigr)$$
where $u_{\pm}\in\cU(\fb^{\pm})$.

If $\fh=\fh_{\ol{0}}$ (that is $R=A$) then $R^{\sigma}=R^*:=\Hom_A(R,A)$.

If  $\fh\not=\fh_{\ol{0}}$ (i.e., $\fg$ is of $Q$-type) the right-hand side
of the above formula is an element of a non-commutative algebra $R$
which is not very convenient. Fortunately, there exists a map
$\int:R\to A$ (see~\ref{tintin}) which induces an isomorphism
$R^{\sigma}\to\Pi^{\dim\fh_{\ol{1}}}(R^*)$. In this way,
$\Coind(R)$ identifies (up to a parity change) with the maximal graded
submodule of $\Hom_A(\cU(\fb^+),A)$. The last identification
gives rise to another realization of
the Shapovalov map $B:\cU(\fb^-)\to\Pi^{\dim\fh_{\ol{1}}}\Hom_A(\cU(\fb^+),A)$.
It is given by the formula 
$$B(u_-)(u_+)=(-1)^{p(u_-)p(u_+)}\int
\HC\bigl(\sigma(u_+)u_-\bigr).$$
The map $B$ is instrumental  in the computation of the centre of $\Ug$
(see~\ref{introcent}).

\subsubsection{}\label{tintin}
Let $\fg$ be a $Q$-type Lie superalgebra.
The algebra $R=\Uh$ is a Clifford superalgebra over 
the polynomial algebra $A=\Sh$. For each $\lambda\in\fh_{\ol{0}}^*$
the evaluation of $R$ at $\lambda$ is a complex Clifford superalgebra.
Notice that a non-degenerate complex Clifford superalgebra
is either the matrix algebra  (if $\dim\fh_{\ol{1}}$ is even) 
or the algebra $Q(n)$ (this is an associative algebra
whose Lie algebra is $\fq(n)$). In particular,
it possesses a supertrace
which is even if $\dim\fh_{\ol{1}}$ is even and odd 
if $\dim\fh_{\ol{1}}$ is odd.
In both cases, there exists a map
$$
\int: R\to A. 
$$
satisfying $\int [R,R]=0$; the evaluation of $\int$ at $\lambda$
is proportional to supertrace on the complex 
Clifford superalgebra if the latter is non-degenerate.

\subsection{Content of the paper}
In Section~\ref{prelim}  we recall definitions and some properties
of main objects. 

In Section~\ref{sdefshap} we propose definition of Shapovalov
map for $Q$-type Lie superalgebras 
which was briefly explained in~\ref{secintro1},~\ref{secintro2}.

In Section~\ref{conn} we compare Shapovalov determinants for various algebras.

In Section~\ref{XY} we construct a non-graded isomorphism $\cM\to\cN$.

In Section~\ref{sectq2} we consider an example  $\fg=\fsq(2)$.

In Section~\ref{leadterm} we calculate the leading terms of Shapovalov
determinants.

In Section~\ref{sectjan} we adapt the definition of Jantzen filtration
(see~\cite{ja}) to the $Q$-type Lie superalgebras.  
As in the contragredient case, the Jantzen filtration is instrumental 
for computations of Shapovalov determinants.

In Section~\ref{sectantie} we describe the anticentres of $Q$-type Lie
superalgebras.

In Section~\ref{shapf} we compute Shapovalov determinants 
(see~\Thm{thmdetsh}). We also show that the Jantzen filtration
have length $2$ for {\em subregular} (see~\ref{gamm})
values of $\lambda$.

In Section~\ref{secthat} we describe a certain completion $\hat{U}$
of $\Ug$.  We show that for any $a\in\cZ(\hat{U})$
one has $\deg\HC(a)=\deg(a)$. As a consequence, $\cZ(\hat{U})=\Zg$
if $\fg$ is a finite-dimensional contragredient 
or $Q$-type Lie superalgebra. 

In Section~\ref{sectcentre} we  describe the centre $\cZ(\hat{U})=\Zg$. 
In~\ref{centcon} we prove that
$\cZ(\fsq(n))=\cZ(\fq(n))$ and $\cZ(\fpsq(n))=\cZ(\fpq(n))$.

In the appendix~\ref{appn} we analyze the structure of
 $\Uh$ which is a Clifford algebra
over the polynomial algebra
$\Sh$. We recall some basic facts on Clifford algebras and
introduce the map $\int:\Uh\to\Sh$. We adapt
a notion of reduced norm to $\Uh$.

\subsection{Acknowledgment}
I am grateful to J.~Bernstein for explaining
the functorial approach to construction of Shapovalov maps
and for numerous suggestions. A part of this work was done
during  my stay at Max-Planck Institut f\"ur Mathematik at Bonn.
I am grateful to this institution for stimulating atmosphere and excellent
working conditions.

\section{Index of notations}
\label{index}
Symbols used frequently are given below under the section number
where they are first defined. Notation used in Appendix are defined
there.

$$\begin{array}{llccll}
\ref{prelim1}&\mathbb{Z}_{\geq 0}, \mathbb{Z}_{>0}, p(u), \Pi,
\deg u, V^{\oplus r}, \gr & & &
& \\
\ref{sigma}&\sigma & & &\ref{Phi+0} & \Phi_0, \Ind, \Coind, \Xi\\
\ref{halp} & h_{\alpha}, h_{\ol{\alpha}}, e_{\alpha},
f_{\alpha}, H_{\alpha}, E_{\alpha}, F_{\alpha} & & &
\ref{cN} & \cM,\cM_{\nu}, R^{\sigma}, \cN, N^*, S, B\\
\ref{QQ}   & Q^+, \nu\geq\mu & & &
\ref{Snu} & S_{\nu}, B_{\nu}, \Norm S_{\nu}\\
\ref{lack} & A, R, \Cl(\lambda), c(\lambda), E(\lambda)& & &
\ref{partition} & \tau(\nu), 
\mathbf{k}, \cP(\nu), \tau_{{\alpha}}(\nu), |\mathbf{k}|\\
\ref{Mlam}& M(\lambda), N(\lambda), V(\lambda) & & &\ref{Tp}& T_{\fp}\\
\ref{olN}& \ol{N} & & &\ref{TUg} & t_{\fh},  t_{\fg}\\
\ref{HCh} & \HC  & & & \ref{gamm} &\Gamma,\gamma_{h,c} \\
\ref{cCO}& \cC, \cD, \cD_+, N_{\nu} & & & 
\ref{integral} &\int.
\end{array}$$

\section{Preliminaries}
\label{prelim}

\subsection{}
\label{prelim1}
The symbol $\mathbb{Z}_{\geq 0}$ stands 
for the set of non-negative integers and 
$\mathbb{Z}_{>0}$ for the set of positive integers. We denote
by $|X|$ the number of elements in a finite set $X$.

Let $V=V_{\ol{0}}\oplus V_{\ol{1}}$ be a  $\mathbb{Z}_2$-graded
vector space.
We denote by $\dim V$ the total dimension of $V$.
For a homogeneous element $u\in V$ 
we denote by $p(u)$ its 
$\mathbb{Z}_2$-degree; in all formulae where this notation is used,
$u$ is assumed to be $\mathbb{Z}_2$-homogeneous. 
For a subspace $N\subset V$ we set $N_i:=N\cap V_i$ for $i=0,1$.
Let $\Pi$ be the functor which switches parity, i.e. $(\Pi V)_{\ol{0}}=V_{\ol{1}},
(\Pi V)_{\ol{1}}=V_{\ol{0}}$. We denote by $V^{\oplus r}$ the direct sum of $r$-copies
of $V$.

For a Lie superalgebra $\fg$ we denote by $\Ug$ its universal enveloping
algebra and  by $\Sg$ its symmetric algebra. 
Recall that $\Sg=\gr\Ug$ with respect to the canonical filtration 
$\cF^k(\fg):=\fg^k$. 
For $u\in\Ug$ denote by $\deg u$ the degree 
of $\gr u$ in $\Sg$.

Throughout the paper the base field is $\mathbb{C}$ and
$\fg=\fg_{\ol{0}}\oplus\fg_{\ol{1}}$ denote one 
(unless otherwise specified, an arbitrary one)
of $Q$-type Lie superalgebras $\fq(n),\fsq(n)$  for $n\geq 2$,
$\fpq(n),\fpsq(n)$ for $n\geq 3$.

Throughout the paper we are often dealing with homomorphisms
between two isomorphic free $R$-modules where $R$ is a commutative algebra.
The determinant of such a homomorphism is defined up to an 
invertible element of $R$;
in all our examples the set of invertible elements of $R$ is $\mathbb{C}^*$.

\subsection{$Q$-type Lie superalgebras}
\label{qtype}
Recall that $\fq(n)$ consists of the matrices with the block form
$$X_{A,B}:=\begin{pmatrix}
\ \ A \ &| &B\\
-- &-&--\\
\ \ B \ & | & A 
\end{pmatrix}$$
where $A,B$ are arbitrary $n\times n$ matrices; 
$\fq(n)_{\ol{0}}=\{X_{A,0}\}\cong \fgl(n),\ \fq(n)_{\ol{1}}=\{X_{0,B}\}$ and
$$[X_{A,0},X_{A',0}]=X_{[A,A'],0},\ \ [X_{A,0},X_{0,B}]=X_{0,[A,B]},\ \ 
[X_{0,B},X_{0,B'}]=X_{0,BB'+B'B}.$$

Define $\tr': \fq(n)\to \mathbb{C}$ by $\tr'(X_{A,B})=\tr B$.
In this notation, 
$$\begin{array}{rl}
\fsq(n):&=\{x\in \fq(n)|\ \tr' x=0\},\\
\fpq(n):&=\fq(n)/(\Id),\\
\fpsq(n):&=\fsq(n)/(\Id),
\end{array}$$
where $\Id$ is the identity matrix.

These definitions are illustrated by the following diagram
$$
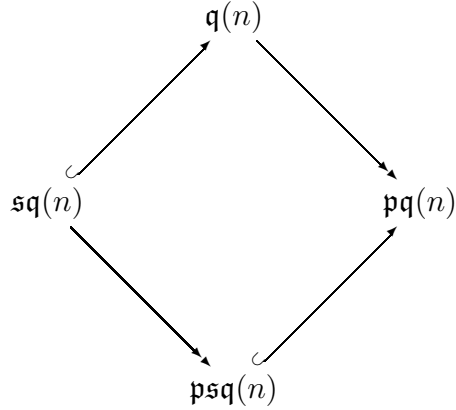
\begin{diagram}
       &         & \fq(n) &        &       \\
       &\ruInto  &        &\rdOnto &       \\
\fsq(n)&         &        &        &\fpq(n) \\
       &\rdOnto  &        &\ruInto &        \\
       &         &\fpsq(n)&        &         \\
\end{diagram}
$$

Clearly, the category of $\fpq(n)$-modules (resp., $\fpsq(n)$-modules) 
is the subcategory of $\fq(n)$-modules (resp., of $\fsq(n)$-modules) 
which are killed by the identity matrix $\Id$. 

The map $(x,y)\mapsto \tr'(xy)$ gives an odd non-degenerate
invariant symmetric bilinear form on $\fq(n)$ and on $\fpsq(n)$.

For the quotient algebras $\fpq(n),\fpsq(n)$ we denote by $X_{A,B}$
the image of the corresponding element in the appropriate
algebra.

\subsubsection{}
\label{sigma}
Recall that a linear map $\sigma$ is called an {\em antiautomorphism}
of a Lie superalgebra (resp., of an associative superalgebra) if it satisfies
 the rule $\sigma([x,y])=(-1)^{p(x)p(y)}[\sigma(y),\sigma(x)]$ (resp.,
$\sigma(xy)=(-1)^{p(x)p(y)}\sigma(y)\sigma(x)$). A Lie superalgebra $\fg$
admits an antiautomorphism $\sigma$ given by $\sigma(x)=-x$;
we denote by $\sigma$ also the induced antiautomorphisms of $\Ug$.

If $\fg$ is a classical Lie superalgebra which is not of type $P$,
it admits a ``naive antiautomorphism'' $x\mapsto x^t$
satisfying the rule $[x,y]^t=[y^t,x^t]$. In an appropriate basis,
this antiautomorphism is given by the matrix transposition.
It preserves the elements of a Cartan subalgebra.

\subsubsection{}
\label{trdec}
For $Q$-type Lie superalgebras the set of even roots 
($\Delta_{\ol{0}}^+$) coincides 
with the set of odd roots ($\Delta_{\ol{1}}^+$). 
This phenomenon has two obvious
consequence. The first one is that all triangular decompositions of a 
$Q$-type Lie superalgebra are conjugate with respect to inner automorphisms
(this does not hold for other simple Lie superalgebras).
The second one is that the element 
$\rho:=\frac{1}{2}(\sum_{\alpha\in\Delta_{\ol{0}}^+}\alpha-
\sum_{\alpha\in\Delta_{\ol{1}}^+}\alpha)$ is equal to zero.

We choose the natural triangular decomposition:
$\fq(n)=\fn^-\oplus\fh\oplus\fn^+$ where $\fh_{\ol{0}}$ consists of the elements
$X_{A,0}$ where $A$ is diagonal, $\fh_{\ol{1}}$ consists of the elements
$X_{0,B}$ where $B$ is diagonal,  and $\fn^+$ (resp., $\fn^-$)
consists of the elements $X_{A,B}$ where $A,B$ are
strictly upper-triangular (resp., lower-triangular). We consider
the induced triangular decompositions of 
$\fsq(n), \fpq(n), \fpsq(n)$. 

The ``naive antiautomorphism'' $x\mapsto x^t$ 
preserves the elements of $\fh$ and
interchanges $\fn^+$ with $\fn^-$.

\subsection{}
\label{notat}
In the standard notation the set of roots of $\fgl(n)=\fq(n)_{\ol{0}}$ 
can be written as
$$\Delta^+=\{\veps_i-\veps_j\}_{1\leq i<j\leq n}$$
and the set of simple roots as $\pi:=\{\veps_1-\veps_2,\ldots,
\veps_{n-1}-\veps_n\}$. Each root space has dimension $(1|1)$.

For $\alpha\in\Delta^+$ let $s_{\alpha}:\fh_{\ol{0}}^*\to\fh_{\ol{0}}^*$ 
be the corresponding reflection: $s_{\veps_i-\veps_j}(\veps_i)=\veps_j$,
$s_{\veps_i-\veps_j}(\veps_k)=\veps_k$ for  $k\not=i,j$.
Denote by $W$ the Weyl group of $\fg_{\ol{0}}$ that is the group
generated by $s_{\alpha}:\alpha\in\Delta^+$.
Recall that $W$ is generated by $s_{\alpha}:\alpha\in\pi$.

The space $\fh_{\ol{0}}^*$ has the standard non-degenerate $W$-invariant bilinear
form: $(\veps_i,\veps_j)=\delta_{ij}$.

\subsubsection{}
\label{hi}
Let $E_{rs}$ be the elementary matrix:
$E_{rs}=(\delta_{ir}\delta_{sj})_{i,j=1}^n$.

The elements
$$h_i:=X_{E_{ii},0}$$
form the standard basis of $\fh_{\ol{0}}$ for $\fg=\fq(n),\fsq(n)$.
We use the notation $h_i$ also
for the image of $h_i$ in the quotient algebras $\fpq(n),\fpsq(n)$.

The elements $H_i:=X_{0,E_{ii}}$ ($i=1,\ldots,n$) form a convenient basis 
of $\fh_{\ol{1}}\subset\fq(n)$; they satisfy the relations
$[H_i,H_j]=2\delta_{ij}h_i$.

\subsubsection{}
\label{halp}
For each positive root $\alpha=\veps_i-\veps_j$ we define
$$\begin{array}{lcl}
h_{\alpha}:=h_i-h_j,\ h_{\ol{\alpha}}:=h_i+h_j,
& & H_{\alpha}:=H_i-H_j,\\
e_{\alpha}:=X_{E_{ij},0},&  & E_{\alpha}:=X_{0,E_{ij}},\\
f_{\alpha}:=X_{E_{ji},0},& & F_{\alpha}:=X_{0,E_{ji}}.
\end{array}$$
All above elements are non-zero in $\fsq(n), \fpq(n), \fpsq(n)$
(since we excluded the cases $\fpq(2),\fpsq(2)$). 

The elements $h_{\alpha}, e_{\alpha},f_{\alpha}$ ($\alpha\in\Delta^+$)
form the standard basis of $\fsl(n)=[\fgl(n),\fgl(n)]$;
the elements $E_{\alpha}$ (resp., $F_{\alpha}$) form the natural
basis of $\fn_{\ol{1}}^+$ (resp., of $\fn^-_{\ol{1}}$) and the elements $H_{\alpha}$
span $\fh_{\ol{1}}\cap\fsq(n)$.

For each $\alpha$ the elements $h_{\alpha}, e_{\alpha},f_{\alpha},
h_{\ol{\alpha}},H_{\alpha},E_{\alpha},F_{\alpha}$ span $\fsq(2)$
 and one has 
$$\begin{array}{c}

[e_{\alpha},f_{\alpha}]=h_{\alpha},\ \
 [E_{\alpha},F_{\alpha}]=h_{\ol{\alpha}},\ \ 
[H_{\alpha},H_{\alpha}]=2h_{\ol{\alpha}}\\

[E_{\alpha},f_{\alpha}]=[e_{\alpha},F_{\alpha}]=H_{\alpha}.
\end{array}$$

\subsubsection{}\label{QQ}
Set
$$Q(\pi):=\sum_{\alpha\in\Delta^+}\mathbb{Z}\alpha,\ \ \ 
 Q^+:=\sum_{\alpha\in\Delta^+}\mathbb{Z}_{\geq 0}\alpha.$$
Define a partial order on $\fh^*_{\ol{0}}$ by $\nu\geq\mu$ iff 
$\nu-\mu\in Q^+$.

\subsection{}
\label{lack}
Set
$$A:=\Sh,\ \ \ R:=\Uh.$$ 
Identify $\cU(\fh_{\ol{0}})$ with $A$. The algebra $R$
is a Clifford algebra over $A$:
it is generated by the odd space $\fh_{\ol{1}}$ endowed by the $A$-valued
symmetric bilinear form $b(H,H')=[H,H']$.
We will describe some properties of this algebra in Appendix.

For $\lambda\in\fh^*_{\ol{0}}$ let $\mathbb{C}(\lambda)$ be the corresponding
one-dimensional $\fh_{\ol{0}}$-module.
Set 
$$\Cl(\lambda):=\Uh\otimes_{\fh_{\ol{0}}}\mathbb{C}_{\lambda}.$$
Clearly, $\Cl(\lambda)$ is isomorphic to a complex
Clifford algebra generated by $\fh_{\ol{1}}$ endowed by the evaluated
symmetric bilinear form $b_{\lambda}(H,H'):=[H,H'](\lambda)$.
Set
$$c(\lambda):=\dim\Ker b_{\lambda}.$$
For $\fg=\fq(n)$, $c(\lambda)$ is the number of zeros
among $h_1(\lambda),\ldots,h_n(\lambda)$.
The complex Clifford algebra $\Cl(\lambda)$ is non-degenerate
iff $c(\lambda)=0$.

Denote by $E(\lambda)$ a simple 
$\Cl(\lambda)$-module 
(up to a shift of grading, such a module is unique--- see~\ref{clpi}).
One has $\dim E(\lambda)=2^{[\frac{\dim\fh_{\ol{1}}+1-c(\lambda)}{2}]}$.

\subsubsection{}\label{Mlam}
Set $\fb:=\fh+\fn^+,\fb^-:=\fh+\fn^-$. 
Endow $\Cl(\lambda)$ with the $\fb$-module structure via the trivial
action of $\fn^+$. Set
$$M(\lambda):=\Ind_{\fb}^{\fg}\Cl(\lambda),
\ \ \ N(\lambda):=\Ind_{\fb}^{\fg}E(\lambda).$$
Clearly, $M(\lambda)$ has a finite filtration with the factors
isomorphic to $N(\lambda)$ up to parity change. We call $M(\lambda)$
a {\em Verma module} and $N(\lambda)$ a {\em Weyl module}.
A Weyl module $N(\lambda)$ has a unique maximal submodule
denoted by $V(\lambda)$.

As a $\fg_{\ol{0}}$-module $N(\lambda)$ has a filtration 
whose factors are $\fg_{\ol{0}}$-Verma modules.
In particular, $N(\lambda)$ has a finite length.

\subsubsection{}
\label{olN}
For a diagonalizable $\fh_{\ol{0}}$-module $N$ and a weight $\mu\in\fh_{\ol{0}}^*$
denote by $N_{\mu}$ the corresponding
weight space. Say that
a module $N$ has the highest weight $\lambda$ if 
$N=\sum_{\mu\leq\lambda} N_{\mu}$ and $N_{\lambda}\not=0$.
If all weight spaces $N_{\mu}$ are finite-dimensional we put
$\ch N:=\sum_{\mu} \dim N_{\mu} e^{\mu}$.

If $N$ has a highest weight 
we denote by $\ol{N}$ the sum of all submodules 
which do not meet the highest weight space of $N$.
One has $V(\lambda)=N(\lambda)/\ol{N}(\lambda)$.

\subsection{Harish-Chandra projection}
\label{HCh}
Denote by $\HC$ the Harish-Chandra projection $\HC:\Ug\to\Uh$
along the decomposition $\Ug=\Uh\oplus(\Ug\fn^++\fn^-\Ug)$.

\subsubsection{}
\begin{lem}{easylemHC}
\begin{enumerate}
\item
If $u\in\Ug$ is the product of $n_-$ elements of $\fn^-$, $n_+$
elements of $\fn^+$ and $n_0$ elements of $\fh$ then
$\deg\HC(u)\leq \min(n_-,n_+)+n_0$.
\item
For $i=1,\ldots,k$ let $x_i$ be an element of weight $\alpha_i$ 
and $y_i$ be an element of weight $-\beta_i$,
where $\alpha_i,\beta_i\in\Delta^+$. If
$\HC(\prod_{i=1}^k x_i\prod_{i=1}^k y_i)$
has degree $k$ then the multisets $\{\alpha_i\}_{i=1}^k$
and $\{\beta_i\}_{i=1}^k$ are equal.
\end{enumerate}
\end{lem}
\begin{proof}
The proof of (i) is an easy induction on $n_-+n_++n_0$. Indeed,
write $u=u'x$ or $u=u'xy_1\ldots y_r$ with 
$x\in (\fh+\fn^+), y_1\ldots y_r\in\fn^-,
u'\in\Ug$. If $u=u'x$  then $\HC(u)=\HC(u')\HC(x)$ and
the assertion follows from by induction.
In the case $u=u'xy_1\ldots y_r$, write
$$u=\pm u'y_1\ldots y_r x+\sum_{i=1}^r\pm u'
y_1\ldots y_{i-1}((\ad x)y_i)y_{i+1}\ldots y_r.$$
As we have already checked  
$\deg\HC(u'y_1\ldots y_r x)\leq \min(n_-,n_+)+n_0$.
The remaining summands
are the products of $n'_-$ elements of $\fn^-$, $n'_+$
elements of $\fn^+$ and $n'_0$ elements of $\fh$ where
$n'_-+n'_++n'_0=n_-+n_++n_0-1$. In all cases
$\min(n'_-,n'_+)+n'_0\leq\min(n_-,n_+)+n_0$ which implies (i);
(ii) easily follows from (i).
\end{proof}

\section{Shapovalov map}
\label{sdefshap}
In this section we construct a Shapovalov map and define
an analogue of Shapovalov determinants. In~\ref{cCO}
we define the main objects and formulate the results of this section.
The proofs are given in~\ref{32}--\ref{IndCl}.

In this section $\fg=\fn^-\oplus\fh\oplus\fn^+$ 
is a classical Lie superalgebra.
If $\fg$ is a $Q$-type superalgebra we keep the above notation.
Otherwise $\fh_{\ol{1}}=0$ and we set $M(\lambda):=N(\lambda)$
for $\lambda\in\fh^*$.

\subsection{Brief description of the main results}
\label{cCO}
View $\Ug$ as a $Q(\pi)$-graded algebra via the adjoint action of 
$\fh_{\ol{0}}$.
Let $\tilde{D}$ be the category of left $Q(\pi)$-graded $\fg$-modules.
Let $Q^-$ be the set of weights of $\cU(\fn^-)$ that is 
$Q^-=-Q^+$. Let $\cD$ (resp., $\cD_+$) 
be the subcategory of $\tilde{D}$ where the objects
are graded $\fg$-modules for which the graded components outside $Q^-$ 
(resp., outside $Q^+$) vanish. Let $\cC$ be the category of left 
$\fh$-modules.
For $K,L\in\cD$  let $\Hom_{\cD}(K,L)$
be the space of degree zero homomorphisms.
For $N\in\cD$ we write $N=\sum_{\nu\in Q^-} N_{\nu}$.
View $M(\lambda)$ as an object of $\cD$: the $Q^-$-grading is defined by
assigning degree zero to the highest weight vectors.

\subsubsection{}\label{Phi+0}
Let $\Phi_0: \cD\to\cC$ be the functor given by $\Phi_0(N)=N_0$.
The functor $\Phi_0$ admits a left adjoint functor $\Ind:\cC\to\cD$
given by $\Ind(L)=\Ind_{\fb}^{\fg} L$ where the action of $\fn^+$ on $L$ 
is supposed to be trivial. It turns out that  $\Phi_0$ also admits
a right adjoint functor which we denote
by $\Coind$. The adjunction morphism $L\to \Ind(L)_0$ is
an isomorphism for any $L\in\cC$;
thus
$$\Hom_{\cD}(\Ind(L),\Coind(L))=\Hom_{\cC}(\Ind(L)_0,L)=
\Hom_{\cC}(L,L).$$

Let $\Xi(L): \Ind(L)\to\Coind(L)$ be the morphism corresponding to
the identity map $\id_L$; in this way we obtain  
a morphism of functors $\Xi:\Ind\to\Coind$.  The following
claim is proven in~\ref{olInd} below.

{\em Claim.} $\Ker \Xi(L)$ is the maximal graded submodule 
of $\Ind(L)$ which does not meet the zero component.

View $R=\Uh$ as  an $\fh$-bimodule.
As a left module $R$ belongs to $\cC$; both $\Ind(R),\Coind(R)$
inherit the right action of $\fh$ and $\Xi(R): \Ind(R)\to\Coind(R)$
is a $\fg$-$\fh$ bimodule map. We call $S:=\Xi(R)$ the {\em Shapovalov map}.
The following claim is proven in~\ref{freyd} below.

{\em Claim.} One has canonical isomorphisms
$\Ind(E)\iso\Ind(R)\otimes_R E,\ \Coind(R)\otimes_R E\iso
\Coind(E)$ identifying
$\Xi(E)$ with $\Xi(R)\otimes_R  \id_E$.

\subsubsection{}\label{cN}
Notice that $\Ind(R)$ has a nice structure: it 
can be identified with $\cU(\fb^-)$ as $\fb^-$-$R$ bimodule.
Set
$$\cM:=\Ind(R).$$
The $\nu$th homogeneous component with respect to $Q^-$-grading takes form
$$\cM_{\nu}=\{v\in \cM|\ hv-vh=-\mu(h)v \text{ for all } h\in\fh_{\ol{0}}\}.$$

We shall now present
convenient realizations of $\Coind(R)$. The proofs are given
in~\ref{freyd} below.

Denote  by $\Hom_{R_r}(-,-)$ the set
of homomorphisms of {\em right} $R$-modules.
Define on $R$ a new bimodule structure $R^{\sigma}$
via $v.r:=(-1)^{p(r)p(v)}\sigma(r)v,
r.v:=(-1)^{p(r)p(v)}v\sigma(r)$ where the dot stands for the new action, 
$\sigma$ is the antiautomorphism introduced in~\ref{sigma},
$r$ is an element of the algebra $R$ and
$v\in R^{\sigma}$. It turns out that $R^{\sigma}$
is isomorphic to $R^*:=\Hom_A(R,A)$ up to a change of parity.
In~\ref{Rsigma} we exhibit various
connections between $R^{\sigma}$ and $\Hom_{R_r}(R,-)$.

Define the functor $\Ind_+:\cC\to\cD_+$ 
similarly to $\Ind$. Let  $\cN$
be the maximal graded subspace of $\Hom_{R_r}(\Ind_+(R),R^{\sigma})$ 
that is 
$$\cN:=\oplus_{\nu\in Q^+} \cN_{-\nu},\ \text{
where } \cN_{-\nu}:=\Hom_{R_r}(\Ind_+(R)_{\nu},R^{\sigma}).$$
Convert the natural structure of $R$-$\fg$-bimodule on $\cN$ to
a $\fg$-$R$-bimodule structure via the antiautomorphism
$\sigma$. Then $\cN\in\cD$ and, moreover, $\cN\cong\Coind(R)$ 
as $\fg$-$R$-bimodules.

For a $\fg$-$R$ bimodule $N$ define a dual module $N^*$ as $\Hom_A(N,A)$
(where $A=\Sh$)
endowed by the $\fg$-$R$ bimodule structure via $\sigma$.
In~\ref{integral} we describe an $A$-homomorphism $\int: R\to A$ 
whose parity is equal to the parity of $\dim\fh_{\ol{1}}$.
For $N$ being a free $R$-module,
the map $\psi\mapsto \int\psi$ provides a map  
$\Hom_{R_r}(N,R^{\sigma})\iso \Pi^{\dim\fh_{\ol{1}}}(N^*)$. 
Both maps are isomorphism of $\fg$-$R$ bimodules
(see~\ref{Rsigma}). Putting $N:=\Ind_+(R)$ we obtain
$\cN\iso\Ind_+(R)^*$
if $\dim\fh_{\ol{1}}$
is even and $\cN\iso\Pi(\Ind_+(R)^*)$
if $\dim\fh_{\ol{1}}$ is odd. Thus we obtain two realizations of $\Coind(R)$:
$\cN$ and $\Pi^{\dim\fh_{\ol{1}}}(\Ind_+(R)^*)$.

\subsubsection{}
The above realizations of $\Coind(R)$ give 
the following realizations of the Shapovalov map:
$S: \Ind(R)\to\Hom_{R_r}(\Ind_+(R),R^{\sigma})$ and
$B:\Ind(R)\to\Pi^{\dim\fh_{\ol{1}}}(\Ind_+(R)^*)$.
Using the natural identification $\Ind(R)=\cU(\fb^-), \Ind_+(R)=\cU(\fb^+)$
we obtain the following formulas:

\begin{equation}\label{Shimap}
S(u_-)(u_+)=(-1)^{p(u_-)p(u_+)}
\HC\bigl(\sigma(u_+)u_-\bigr).
\end{equation}

\begin{equation}\label{Bimap}
B(u_-)(u_+)=(-1)^{p(u_-)p(u_+)}\int\!
\HC\bigl(\sigma(u_+)u_-\bigr).
\end{equation}

\subsubsection{}\label{Snu}
Recall that $S$ is homogeneous of degree zero and write 
$S=\sum_{\nu\in Q^+} S_{\nu}, B=\sum_{\nu\in Q^+} B_{\nu}$
where $S_{\nu}:\cM_{-\nu}\to \cN_{-\nu}$
is the restriction of $S$ and $B_{\nu}$ is defined similarly. Observe that 
$S_{\nu}$ is an $A$-homomorphism between two free $A$-modules
of the same finite rank. 
Thus $\det S_{\nu}\in A=\Sh$ is defined up to an invertible
scalar. If $\fg$ is a contragredient Lie superalgebra,
$\det S_{\nu}$ is called a {\em Shapovalov determinant}.

In~\ref{Sh2n} we show that for a $Q$-type Lie superalgebras
$$\det S_{\nu}=(\Norm S_{\nu})^{2^{\dim\fh_{\ol{1}}}}$$
where $\Norm S_{\nu}\in A$ is a {\em reduced norm} of the operator
$S_{\nu}$. We suggest $\Norm S_{\nu}$ as an analogue of
Shapovalov determinant for $Q$-type superalgebras.

Since $B$ is the composition of $S$ and a $\fg$-$R$ isomorphism,
one has
$$\det S_{\nu}=\det B_{\nu}.$$
The map $B$ is more convenient than $S$ because
 the matrices of $B_{\nu}$ have
entries in the polynomial algebra $A$.
In particular, $\det B_{\nu}$ is equal to the determinant of the
corresponding matrix which we also denote as $B_{\nu}$.

\subsubsection{}\label{xilambda}
Observe that $M(\lambda)=\Ind(\Cl(\lambda))$.
The map $\Xi(\Cl(\lambda))$ is obtained from $\Xi(R)$
by the evaluation at $\lambda$. Its kernel is $\ol{M(\lambda)}$
(see~\ref{olN} for the notation) and it
coincides with $\Ker B(\lambda)$. This gives
\begin{equation}
\label{olNla}
\ol{M(\lambda)}=\{uv_{\lambda}|\ \bigl(\int\HC(u'u)\bigr)(\lambda)=0\ 
\text{ for all } u\in \cU(\fb^+)\}
\end{equation}
where $v_{\lambda}$ is the canonical generator of $M(\lambda)$
(i.e., the image of $1\in R$).

Recall that $M(\lambda)$ has a filtration with factors
isomorphic to $N(\lambda)$ up to the change of parity.
Thus $\ol{M(\lambda)}=0$ if and only if $\ol{N(\lambda)}=0$.
Hence $\ol{N(\lambda)}=0$
if and only if $S(\lambda)$ is injective.

{\em Corollary.} $N(\lambda)$ is simple iff $\Norm S_{\nu}(\lambda)\not=0$
for all $\nu\in Q^+$.

Notice that the matrices of
the evaluated maps $B_{\nu}(\lambda)$ have complex entries.
In particular, the dimension
of the kernel of the evaluated map $B_{\nu}(\lambda)$ is
equal to the corank of the corresponding matrix.

{\em Corollary.} 
$\dim\ol{M(\lambda)}_{\nu}=\corank B_{\nu}(\lambda)$.

\subsubsection{}
Let $\fg$ be a $Q$-type Lie superalgebra.
As we  show later $\cap_{\lambda\in\fh^*_{\ol{0}}} \Ann N(\lambda)=0$.
Since $N(\lambda)$ is a subquotient of $\cM$
one has $\Ann_{\Ug}\cM=0$.
We will use the module $\cM$ for the calculation of the centre
of $\Ug$ in Sect.~\ref{sectcentre}.

\subsection{The proofs of the claims~\ref{Phi+0}}
\label{32}
Retain the notation of~\ref{cCO}. 
It is easy to check that $\Ind$ is left adjoint to $\Phi_0$ and one has
a canonical isomorphism
$$\Ind(L)\iso\Ind(R)\otimes_R  L.$$

Recall that $\cN\in\cD$ is the maximal graded subspace of 
$\Hom_{R_r}(\Ind_+(R),R^{\sigma})$. We will identify $R^{\sigma}$ with $R$
and use the dot to indicate the $R$-bimodule structure on $R^{\sigma}$. Then
$\fg$-$R$-bimodule structure is given by the following formulas

\begin{equation}\label{cN1}
\begin{array}{l}
\psi(xr)=\psi(x).r=(-1)^{p(r)p(\psi(x))}\sigma(r)\psi(x),\\
(g\psi)(x)=(-1)^{p(g)p(\psi)}\psi(\sigma(g)x),\\
(\psi r)(x)=(-1)^{p(r)p(\psi)}\sigma(r).\psi(x)=(-1)^{p(r)p(x)}\psi(x)r.
\end{array}\end{equation}
where $g\in\Ug, r\in R,\psi\in\cN,x\in\Ind_+(R)$.

\subsubsection{}
\begin{prop}{freyd}
\begin{enumerate}
\item
The functor $\Phi_0$ admits a right adjoint functor $\Coind$ 
which is exact.
\item
For any $K\in\cC$ the adjunction map 
$\Coind(K)_0\to K$ is an isomorphism.
\item One has a $\fg$-$R$-bimodule isomorphism $\beta:\cN\iso\Coind(R)$
satisfying $\beta^{-1}(x)(y)=(-1)^{p(x)p(y)}
\sigma(y)x$ for each $x\in\Coind(R)_0=R$
and $y\in\Ind_+(R)_0=R$.

\item One has a canonical isomorphism
$\Coind(R)\otimes_R  L\to\Coind(L)$.

\item Under the above identifications 
$\Xi(L)$ identifies with $\Xi(R)\otimes_R  \id_L$.
\end{enumerate}
\end{prop}
\begin{proof}
(i) The functor $\Phi_0$ is exact. By Freyd Theorem
(see~\cite{mc} Ch. V), $\Phi_0$ admits a  right adjoint functor.
The functor $\Coind$ is left exact as it admits left adjoint. 
It remains to check that an epimorphism
$L\to L'$ induces an epimorphism $\Coind(L)\to \Coind(L')$.

Below we construct a family of modules $X({\nu})\in\cD$ ($\nu\in Q^+$)
with the following properties: 
\begin{itemize}
\item $X(\nu)_0$ is a free $R$-module
and 
\item one has a natural isomorphism
$\Hom_{\cD}(X(\nu),N)\iso N_{-\nu}$ for all $N\in\cD$. 
\end{itemize}
To check the
surjectivity of $\Coind(L)\to \Coind(L')$
it is enough to check it on the homogeneous components that is
to verify the surjectivity of the map
$$\Hom_{\cD}(X(\nu),\Coind(L))\to
\Hom_{\cD}(X(\nu),\Coind(L')).$$ 
The latter amounts to the surjectivity
of $\Hom_{\cC}(X(\nu)_0,L)\to\Hom_{\cC}(X(\nu)_0,L')$ which follows
from the freeness of $X(\nu)_0$.
We construct $X(\nu)\in\cD$ as follows: 
fix a $Q(\pi)$-grading on $\Ug$ by assigning to $1$ degree $-\nu$
and then take the maximal quotient belonging to $\cD$; in other words,
$X(\nu):=\Ug/\sum_{\lambda\not\leq\nu} \Ug\Ug_{\lambda}$
with the grading shifted by $\nu$. One can easily see from PBW theorem
that $X(\nu)_0$ can be identified with $\cU(\fb^+)_{\nu}$.
Hence $X(\nu)_0$ is a free left $R$-module;
the second property is clear. 
This proves (i).

Observe that $X(\nu)$ is a cyclic $\fg$-module and that $X(\nu)_0$
is a finitely generated $R$-module.

(ii)
Identify $\Ind(R)_0$ with $R$ (see~\ref{32}). One has
$$K=\Hom_{\cC}(R,K)=\Hom_{\cD}(\Ind(R),\Coind(K))=
\Hom_{\cC}(R,\Coind(K)_0).$$
Hence $K=\Hom_{\fh}(R,\Coind(K)_0)=\Coind(K)_0$ as required.

(iii) One has an $R$-bimodule isomorphism
 $R\iso \cN_0=\Hom_{R_r}(R,R^{\sigma})$ given by $x\mapsto r_x$
where $r_x(y)=\sigma(x)y$. Observe that the evaluated map
$\Cl(\lambda)\to \Hom_{R_r}(R,\Cl(\lambda))$ is bijective
for any $\lambda$. The inverse map 
$\beta':\cN_0\to R$ induces an $\fg$-$R$ homomorphism
$\beta:\cN\to\Coind(R)$ satisfying $\beta_0=\beta'$.
We show below that
$\beta_{-\nu}\cN_{-\nu}\to\Coind(R)_{-\nu}$ is a morphism of 
free $A$-modules of the same finite rank and that
the evaluated maps $\beta(\lambda)$ 
are injective for all $\lambda$. Thus $\det\beta_{-\nu}$ has nonzero
evaluations at all points and hence is invertible.
Therefore $\beta_{-\nu}$ is bijective for all
$\nu$ and so $\beta$ is bijective as well.

For each $\lambda\in\fh_{\ol{0}}^*$ the evaluated map $\beta(\lambda)$ 
has the source
$\cN(\lambda):=\Hom_{R_r}(\Ind_+(R),\Cl(\lambda)^{\sigma})$.
It is easy to see that any non-zero submodule of $\cN(\lambda)$ 
meets $\cN(\lambda)_0$. Since $\beta(\lambda)_0=\beta'(\lambda)$ 
is injective,
the map $\beta(\lambda)$ is injective as well.

Using (i), one obtains  the natural
isomorphisms of right $R$-modules 
$$\Coind(R)_{-\nu}=\Hom_{\cD}(X(\nu),\Coind(R))=
\Hom_{\cC}(X(\nu)_0,R)=\Hom_{\fh}(\cU(\fb^+)_{\nu},R)=\cN_{-\nu}$$ 
where the right action of $R$ on $\phi\in\Hom(X,Y)$
is induced by the right action of $R$ on $Y$. 
Hence the source and the target of $\beta_{-\nu}$ are free $A$-modules
of the same finite rank.
This establishes (iii).

(iv)  Define a canonical map
\begin{equation}
\label{eq:coind=coind}
 \Coind(R)\otimes_R L\rTo\Coind(L)
\end{equation}
as the one adjoint to the map
$$ (\Coind(R)\otimes_R L)_0\rTo L$$
obtained from the identification $\Coind(R)_0=R$ (see (ii)).

Recall that $\Coind$ is an exact functor (see (i)).
It is easy to show that the map in
(\ref{eq:coind=coind}) is an isomorphism for each $L$ iff
the functor $\Coind$ commutes with infinite direct sums.
To verify the latter observe that
$$\Hom_{\cD}(X(\nu),\Coind(\oplus Y_i))=\Hom_{\cC}(X(\nu)_0,\oplus Y_i)=
\oplus\Hom_{\cC}(X(\nu)_0, Y_i)$$
because $X(\nu)_0$ is a finitely generated $R$-module.
On the other hand, since $X(\nu)$ is finitely generated $\fg$-module one has
$$\Hom_{\cD}(X(\nu),\oplus \Coind(Y_i))=\oplus\Hom_{\cD}(X(\nu),\Coind(Y_i))=
\oplus\Hom_{\cC}(X(\nu)_0, Y_i).$$
This finally yields a bijection  of
the $\nu$th graded components of $\Coind(\oplus Y_i)$
and $\oplus\Coind(Y_i)$ (see the proof of (i)). Now (iv) follows.

(v)
The last claim amounts to checking the diagram below
is commutative for each $L$.
$$
\begin{diagram}
\Ind(L) &\rTo^{\Xi(L)} & \Coind(L) \\
\uTo    &            & \uTo \\ 
 \Ind(R)\otimes_RL  & \rTo^{\Xi(R)\otimes\id} & \Coind(R)\otimes_RL
\end{diagram}
$$
Since the functors $\Coind$ and $X\mapsto X_0$ are adjoint, this
follows from the commutativity of the diagram
$$
\begin{diagram}
\Ind(L)_0 &\rTo & L \\
\uTo    &            & \uTo \\ 
 (\Ind(R)\otimes_RL)_0  & \rTo & (\Coind(R)\otimes_RL)_0
\end{diagram}
$$
which is obvious.
\end{proof}

\subsubsection{}
\begin{prop}{olInd} 
$\Ker \Xi(L)$ is the maximal graded submodule 
of $\Ind(L)$ which does not meet $L$.
\end{prop}
\begin{proof}
A composition
$$X\rTo\Ind(L)\rTo^{\Xi(L)}\Coind(L)$$
is zero if and only if the adjoint composition
$$X_0\rTo\Ind(L)_0=L\rTo\Coind(L)$$
is zero. This implies the statement.
\end{proof}

\subsubsection{}
Let us identify $\Coind(R)$ with $\cN\subset \Hom_{R_r}(\Ind_+(R),R^{\sigma})$ 
via $\beta$ (see~\Prop{freyd}). 
By~\Prop{freyd}, the restriction of
$S:\cM\to\cN$ to $\cM_0=R$ is given by
$x\mapsto r_x$ where $r_x(y)=(-1){p(x)p(y)}\sigma(y)x$ if
$y\in \Ind_+(R)_0=R$ and $r_x(y)=0$ if $y\in \Ind_+(R)_{\nu}$
for $\nu\not=0$. 

Viewing $1$ as an element of $\Ind(R)$ (resp., $\Ind_+(R)$) via
the identification $R=\Ind(R)_0$ (resp., $R=\Ind_+(R)_0$)
write an element of $\Ind(R)$ as $u1$ and an element of $\Ind_+(R)$
as $u'1$ ($u,u'\in\Ug$). 
\subsubsection{}\begin{claim}{}
$$S(u1)(u'1)=(-1)^{p(u)p(u')}\HC(\sigma(u')u)1.$$
\end{claim}
\begin{proof}
Since  $S$ is an even $\fg$-homomorphism one has
$$\begin{array}{ll}
S(u1)(u'1)&=(-1)^{p(u)p(u')}(\sigma(u')S)(u1)(1)=
(-1)^{p(u)p(u')}S(\sigma(u')u_-u1)(1)\\
& =(-1)^{p(u)p(u')}S(\HC(\sigma(u')u)1)(1)=(-1)^{p(u)p(u')}
\HC(\sigma(u')u)1.
\end{array}$$
\end{proof}

Identifying $\cM$ with $\cU(\fb^-)$ and $\cN$ with $\cU(\fb^+)$
we obtain the formulas~(\ref{Shimap}),~(\ref{Bimap}).

\subsection{Shapovalov determinants}
Let $S_{\nu}:\cM_{-\nu}\to\cN_{-\nu}$ be the restriction of 
$S$ to $\cM_{-\nu}$.
Viewed as left $A$-modules $\cM_{-\nu}$ and $\cN_{-\nu}$ are free
of the same finite rank. Thus $\det S_{\nu}$ is an element of
$A$ defined up to the multiplication by an 
invertible  element, i.e. by an element of $\mathbb{C}^*$.

If $\fg$ is a contragredient Lie superalgebra, $\det S_{\nu}$ is called a
{\em Shapovalov determinant}.

\subsubsection{Reduced norms}\label{Sh2n}
Let $\fg$ be a $Q$-type Lie superalgebra.

Set $\tilde{R}:=R\otimes R$ where the tensor product means
that of the graded algebras. View $\tilde{R}$ as a non-graded algebra.
Notice that $\tilde{R}$ is a Clifford
$A$-algebra whose evaluation at a generic point is
isomorphic to the matrix algebra $\Mat(k,\mathbb{C})$
where $k=2^{\dim\fh_{\ol{1}}}$.
As we will explain in~\ref{rednorm}, for any $\tilde{R}$-module $L$
which is free over $A$ there exists a unique map
 $\Norm:\End_{\tilde{R}}(L)\to A$ which satisfies the properties
$$\Norm(\id)=1,\ \ \Norm(\psi\psi')=\Norm\psi\Norm\psi',\ \
\det\psi=(\Norm\psi)^k.$$

Convert $R$-bimodules $\cM_{-\nu}$ and $\cN_{-\nu}$ to left (non-graded)
$\tilde{R}$-modules via the antiautomorphism $\sigma$;
denote these modules by $X$ and $Y$ respectively.
In Sect.~\ref{XY} we will construct an isomorphism $\Psi:X\to Y$.
This allows us to define a map $\Norm:\Hom_{\tilde{R}}(X,Y)\to A$
by setting $\Norm \psi:=\Norm(\Psi^{-1}\psi)$. One has
$$\det\psi=(\Norm\psi)^{2^{\dim\fh_{\ol{1}}}}.$$

Since $S_{\nu}$ is an even homomorphism of $R$-bimodules,
it can be viewed as an element of $\Hom_{\tilde{R}}(X,Y)$.
We call $\Norm S_{\nu}$ a Shapovalov determinant.

\subsubsection{}
The results presented in~\ref{xilambda}
immediately follow from the above.

\subsection{Applications to Verma modules $M(\lambda)$.}
\label{IndCl}
Recall that $\Cl(\lambda)$ has a filtration with the factors
of the form $E(\lambda),\Pi(E(\lambda))$. This filtration induces
a filtration on $M(\lambda)$
which has the same number of factors and each factor
is either $N(\lambda)$ or $\Pi(N(\lambda))$. In its turn,
$\ol{M(\lambda)}$ admits a filtration with the factors
of the form $\ol{N(\lambda)}$ or $\Pi(\ol{N(\lambda)})$
and the number of factors is not greater than one for the previous 
filtrations.  Therefore
$$
\ch \ol{M(\lambda)}=m\ch\ol{N(\lambda)}, \text{ for some } 1\leq
m\leq \frac{2^{\dim\fh_{\ol{1}}}}{\dim E(\lambda)}.$$

The example $\fg=\fsq(2)$, $\lambda=0$ (see~\ref{sq2cl0})
illustrates
that the numbers of factors in the above filtrations
can be not equal. Indeed, the module
$N(0)$ has the filtration of length two:
$N(0)/M(0)\cong\Pi(M(0))$, but $\ol{N(0)}=\ol{M}(0)$.
We have $\dim \fh_{\ol{1}}=1,\dim E(\lambda)=1$ and $m=1$.

\subsubsection{}\label{clambda}
If $c(\lambda)=0$ (see~\ref{lack} for the notation), 
we conclude from~\ref{clpi} below that for $\dim\fh_{\ol{1}}$ being even
$$\begin{array}{l}
M(\lambda)\cong N(\lambda)^{\oplus s}\oplus
\Pi(N(\lambda))^{\oplus s},\\
\ol{M}^{\oplus}(\lambda)\cong \ol{M}(\lambda)^{\oplus s}\oplus
\Pi(\ol{M}(\lambda))^{\oplus s},
\end{array}$$
and for $\dim\fh_{\ol{1}}$ being odd
$$\begin{array}{l}
N(\lambda)\cong\Pi(N(\lambda)),\\
M(\lambda)\cong N(\lambda)^{\oplus s},\\
\ol{M(\lambda)}\cong \ol{N(\lambda)}^{\oplus s}
\end{array}$$
where $s:=2^{[\frac{\dim\fh_{\ol{1}}-1}{2}]}$.

\subsubsection{}
Combining~\ref{xilambda} and the above analysis of filtrations we conclude
\begin{equation}
\label{corankt}
\begin{array}{l}
\dim \ol{M(\lambda)}_{\lambda-\nu}=\corank B_{\nu}(\lambda),\\
\corank B_{\nu}(\lambda)/r\leq
\dim\ol{N(\lambda)}_{\lambda-\nu}\leq \corank B_{\nu}(\lambda),\\
\dim V(\lambda)_{\lambda-\nu}=\rank B_{\nu}(\lambda)/r \text{ if } 
c(\lambda)=0
\end{array}\end{equation}
where $r=2^{\dim\fh_{\ol{1}}}/\dim E(\lambda)$. Recall that
the condition $c(\lambda)=0$ simply means that $\Cl(\lambda)$
is a non-degenerate Clifford algebra; the
example~\ref{sq2cl0} shows that
this condition is essential in the last formula.

\subsubsection{}
\begin{cor}{corshap}
\begin{enumerate}
\item
$\ol{M(\lambda)}_{\lambda-\nu}=0\ \Longleftrightarrow\ 
\ol{N(\lambda)}_{\lambda-\nu}=0\ \Longleftrightarrow\ 
\Norm S_{\nu}\not=0$,
\item 
$N(\lambda)$ is simple iff $\Norm S_{\nu}\not=0$ for all $\nu\in Q^+$.
\end{enumerate}
\end{cor}

\subsubsection{Example: $\nu=0$} 
In this case $Y=R$ and the map 
$S_0: R\to\Hom_A(R,A)$ coincides with $\alpha$. From the formula~(\ref{intHJ})
we see that the matrix $B_0$ written with respect to an appropriate basis
have zero entries everywhere except the secondary diagonal; the entries
of this diagonal are equal to $\pm 1$. Thus $\det B_0=\pm 1$
and $\Norm S_0=1$.

\section{On Shapovalov determinants for various algebras}
\label{conn}
\subsection{Assumptions}
The construction described in Sect.~\ref{sdefshap} is applicable
to a large class of Lie superalgebras admitting
``nice triangular decompositions''. The main assumption is
that $[\fh_{\ol{0}},\fh]=0$.

Let $\fg=\fn^-+\fh+\fn^+,\fg'=(\fn^-)'+\fh'+(\fn^+)'$ 
be Lie superalgebras satisfying the assumption. Let
$\psi:\fg\to\fg'$ be a homomorphism 
such that the restriction of $\psi$ gives bijections 
$\fn^{\pm}\iso(\fn^{\pm})'$.
 In this section we show that $\Norm S'_{\nu}=\psi(\Norm S_{\nu})$;
notice that $\det S'_{\nu}\not=\psi(\det S_{\nu})$
if $\dim\fh_{\ol{1}}\not=\dim\fh'_{\ol{1}}$.

\subsubsection{}
Retain notations of Sect.~\ref{sdefshap}; set $R':=\cU(\fh')$,
$\cM':=\Ind(R')$ and so on. Extend $\psi$ to the homomorphism
$\Ug\to\cU(\fg')$. The module $\cM'$ inherits a structure
of $\fg$-$R$ bimodule and it lies in the category $\cD$.
The restriction of $\psi$ to $R$ induces $\fg$-$R$-maps
$\cM\to \cM'$ and $(\cN)'\to\cN$. The following diagram is
commutative 
$$
\begin{diagram}
\cM &\rTo^{S} & \cN \\
\dTo &        &\uTo\\
\cM' &\rTo^{S'} & (\cN)'.
\end{diagram}
$$

In particular, $\cM', (\cN)'$ inherit  $Q^-$-gradings and
$S'=\sum_{\nu\in Q^+} S'_{\nu}$ ($Q^+,Q^-\subset Q$ 
where $Q$ stands for the  root lattice corresponding to $\fg$).

\subsection{}
\begin{thm}{clcon}
One has  $\Norm S'_{\nu}=\psi(\Norm S_{\nu})$ for each $\nu$.
\end{thm}
\begin{proof}
Let $H_1,\ldots, H_n$ (resp., $H_1',\ldots,H_m'$) be a basis of $\fh_{\ol{1}}$
(resp., of $\fh_{\ol{1}}'$) and $\psi(H_i)=H_i'$ for $i=1,\ldots,s$,
$\psi(H_i)=0$ for $=s+1,n$.
Normalize $\int,\int'$ is such a way that
$\int H_1\ldots H_n=\int' H'_1\ldots H'_m=1$.
For $J\subset \{1,\ldots,n\}$ 
set $H_J:=\prod_{j\in J}H_j$ ($H_{\emptyset}=1$)
and define $H'_J$ similarly. Clearly, $H_J$
form a free basis of $\Uh$ over $\Sh$.
Fix bases $\{x_1,\ldots,x_r\}$
in $\cU(\fn^-)_{-\nu}$ and $\{y_1,\ldots,y_r\}$
in $\cU(\fn^+)_{\nu}$. Then
the products $x_iH_J$ form a basis in $\cU(\fb^-)$
and  the products $y_iH_J$ form a basis in $\cU(\fb^+)$. 
Consider Shapovalov matrices 
$B_{\nu}, B'_{\nu}$ corresponding to this choice
of bases. More precisely one has
$$B_{\nu}=\bigl(b_{(i,I; j,J)}\bigr),\ \ 
b_{(i,I; j,J)}=(-1)^{p(x_iH_I)p(y_jH_J)}\int\!
\HC\bigl(\sigma(x_iH_I)y_jH_J\bigr)$$
and the similar formulas for $B'_{\nu}$.

Consider three cases: 

(i) the restrictions of $\psi$ 
gives bijection $\fh_{\ol{1}}\iso \fh'_{\ol{1}}$;

(ii) the restrictions of $\psi$ 
gives a monomorphism $\fh_{\ol{1}}\to \fh'_{\ol{1}}$ and
     $\dim\fh'_{\ol{1}}-\dim\fh_{\ol{1}}=1$;

(iii) the restrictions of $\psi$ gives an epimorphism 
$\fh_{\ol{1}}\to \fh'_{\ol{1}}$ and
     $\dim\fh_{\ol{1}}-\dim\fh'_{\ol{1}}=1$.

In the first case $\Ker\psi\subset\fh_{\ol{0}}$ and so
the matrix $B'_{\nu}$ is the evaluation
of the matrix $B_{\nu}$.
In particular, $\det B'_{\nu}=\psi(\det B_{\nu})$
and so $\Norm S_{\nu}=\psi(\Norm  S'_{\nu})$.

In the second case $n=s=m-1$. One has
$\int' \psi(a)H_m=\int a$ for all $a\in\cU(\fh)$.
Then $\int'\psi(a)\psi(b)H_m=\int ab=\pm\int' H_m\psi(a)\psi(b))$
and $\int'\psi(a)\psi(b)=\int' H_m\psi(a)\psi(b)H_m=0$ 
for all $a,b\in\Ug$.
Therefore 
$$B'_{\nu}=\begin{pmatrix}
0 & B_{\nu}\\
\pm B_{\nu} &0 
\end{pmatrix}$$
Thus $\det B'_{\nu}=\pm (\det B_{\nu})^2$
and so $\Norm S_{\nu}=\Norm  S'_{\nu}$.

The third case $m=s=n-1$ is similar to the second one.

Finally, remark that general $\psi$ can be presented
as the composition of maps of the form (i)--(iii).
Indeed, by the assumption on $\fg$ one has $[\fh_{\ol{0}},\fh]=0$. Since
$\Ker\psi$ is a subspace of $\fh$ and  an ideal,
$\Ker\psi\cap\fh_{\ol{0}}$ is an ideal as well.
Thus $\psi=\psi_2\circ\psi_1$ where $\psi_1$
is of the form (i) and the restriction of $\psi_2$  
gives a bijection $\fh_{\ol{0}}\iso\fh'_{\ol{0}}$. Therefore in is enough to show
that $\psi$ satisfying $\Ker\psi\cap\fh_{\ol{0}}=0$ 
an be presented
as the composition of maps of the form (ii), (iii).
It suffices to construct a chain of ideals 
$I_1\subset I_2\subset\ldots \subset I_k=\Ker\psi\subset\fg$ 
satisfying $\dim I_j=j$ and a chain of subalgebras 
$\im\psi=\fp_{0}\subset\fp_1\subset\ldots \subset \fp_t=\fg'$
satisfying $\dim\fp_i=\dim\fp_{i-1}+1$.
Recall that $\Ker\psi\subset \fh_{\ol{1}}$ and so all
elements of $\Ker\psi$ are central and so any subspace of $\Ker\psi$
is an ideal. Similarly, the assumption $[\fh'_{\ol{0}},\fh']=0$ implies that
the sum $\im\psi+X$ is a subalgebra for any $X\subset\fh_{\ol{1}}$.
The statement follows.
\end{proof}

\subsection{Shapovalov determinants for various $Q$-type algebras}
\label{coon}
Recall that $\fpq(n),\fpsq(n)$ are quotients
of $\fq(n)$ and $\fsq(n)$ respectively, by a one-dimensional centre
$\mathbb{C}h_0$.
As a consequence, the Shapovalov determinants for 
$\fpq(n),\fpsq(n)$ are evaluations at $h_0=0$ of the corresponding
Shapovalov determinants for $\fq(n),\fsq(n)$.

Let $S_{\nu}$ be the Shapovalov map for $\fq(n)$
and $S'_{\nu}$ be the Shapovalov map for $\fsq(n)$.
The natural embedding $\fsq(n)\to\fq(n)$ gives
$\Norm S_{\nu}=\Norm  S'_{\nu}$.

\section{Comparing $\cM$ to $\cN$.}
\label{XY}

The modules $\cM$ and $\cN$ are not isomorphic as $R$-bimodules.
In this section we show that they become (non-canonically) isomorphic 
when considered as non-graded $R\otimes R$-modules.
This fact is used in the definition of reduced norms in~\ref{Sh2n}.

\subsection{$\cM$ and $\cN$ as $R$-bimodules}
Set $\tilde{R}:=R\otimes R$
where the tensor product sign stands for the tensor
product of graded algebras.

Any $R$-bimodule $V$ can be viewed as a left $\tilde{R}$ module
via the antiautomorphism $\sigma$:
$$(u\otimes 1)v:=uv,\ \ (1\otimes u)v:=(-1)^{p(v)p(u)}v\sigma(u).$$

View $\tilde{R}$ as a non-graded algebra; we will prove
that $\cM\cong \cN$ as (non-graded) $\tilde{R}$-modules.

\subsubsection{}
Recall that as $R$-bimodules $\cM=\cU(\fb^-)$ 
and  $\cN=\sum_{\nu\in Q^+} \Hom_{R_r}(\cU(\fb^+)_{\nu},R^{\sigma})$.
As  modules over $1\otimes R$ both $\cM,\cN$ can be decomposed as
$$\cM=(1\otimes R)\otimes \cU(\fn^-) ,\ \ \ 
\cN=(1\otimes R)\otimes  \cU(\fn^+)^{\#}$$
where $\cU(\fn^+)^{\#}$ is the image of a natural embedding 
$\Hom(\sum_{\nu\in Q^+}\cU(\fn^+)_{\nu},\mathbb{C})\to \cN$.
It is easy to see that both $\cU(\fn^-)$ and $\cU(\fn^+)^{\#}$ 
are $\ad\fh$-stable.

We construct in~\ref{X'Y'} below an isomorphism 
$\Psi':\cU(\fn^-)\iso \cU(\fn^+)^{\#}$
of non-graded $\ad\fh$-modules.

Finally, the following~\Lem{Psi'} implies that an isomorphism
of $\ad\fh$ modules $\cU(\fn^-)$ and $\cU(\fn^+)^{\#}$ extends canonically 
to an $\tilde{R}$-isomorphism $\cM\iso\cN$.

\subsubsection{}
\begin{lem}{Psi'}
Let $\fp$ be a Lie superalgebra, $U:=\cU(\fp)$ and let 
$\Delta:U\to U\otimes U$ be the comultiplication.
Let $X,Y$ be non-graded modules over $U\otimes U$ which, as
modules over $1\otimes U\subset U\otimes U$, have the form
$$X=(1\otimes U)\otimes X',\ \ Y=(1\otimes U)\otimes Y'$$
where $X',Y'$ are $\Delta(\fp)$-stable subspaces.

Then any homomorphism of $\Delta(\fp)$-modules
 $\Psi':X'\to Y'$ can be uniquely extended to
a $U\otimes U$-homomorphism $\Psi:X\to Y$.
\end{lem}

Proof is straightforward.

\subsubsection{}
Let us rewrite the above lemma in terms of bimodules.

Let  $\tilde{X},\tilde{Y}$ be $U$-bimodules and
$X,Y$ are corresponding $U\otimes U$-modules
(defined via the antiautomorphism $\sigma$) viewed
as non-graded modules. Notice that the action of $\Delta(\fp)$
on $X, Y$ corresponds to the action of $\ad\fp$ on
$\tilde{X},\tilde{Y}$.

Assume that $\tilde{X}=X'\otimes U,\ \tilde{Y}=Y'\otimes U$
as right modules where $X',Y'$ are $\ad\fp$-stable. \Lem{Psi'}
says that
any non-graded $\ad\fp$-homomorphism 
 $\Psi':X'\to Y'$ can be uniquely extended to
a $U\otimes U$-homomorphism $\Psi:X\to Y$.

Thus, we have proven the following

\subsubsection{}
\begin{cor}{corXY}
 $\cM$ and $\cN$ are isomorphic as non-graded $\tilde{R}$-modules.
\end{cor}

\subsection{An isomorphism $\cU(\fn^-)\to\cU(\fn^+)^{\#}$}
\label{X'Y'}
In this subsection we will construct a 
non-graded $\fh$-isomorphism $\cU(\fn^-)\iso \cU(\fn^+)^{\#}$.
All $\fh$-module structures in this subsection are induced by the adjoint
action. 

\subsubsection{}\label{negr}
Denote by $\cU(\fg_{\alpha})^{\#}$ the maximal 
$\fh$-invariant submodule of $\Hom(\cU(\fg_{\alpha}),\mathbb{C})$.
We claim that $\cU(\fg_{-\alpha})$ and 
$\cU(\fg_{\alpha})^{\#}$
are isomorphic as non-graded $\fh$-modules.
Indeed, $\cU(\fg_{-\alpha})$ and 
$\cU(\fg_{\alpha})^{\#}$ are isomorphic as $\fh_{\ol{0}}$-modules and
their weights are of the form $-r\alpha$; their zero weight
spaces are even trivial $\fh$-modules
and for $r>0$ the  corresponding
subspaces are simple two-dimensional $\fh$-modules.
The algebra $\fh$ has at most two non-isomorphic simple
modules of a given $\fh_{\ol{0}}$-weight, these modules are $E(-r\alpha)$
and $\Pi(E(-r\alpha))$ (see~\ref{clpi}).
Thus as non-graded $\fh$-modules $\cU(\fg_{-\alpha})$ and 
$\cU(\fg_{\alpha})^{\#}$ are isomorphic (see~\Rem{Xalpha} for 
more details).

\subsubsection{}
The product in $\Ug$ gives the following isomorphisms of graded $\fh$-modules:
$$\cU(\fn^-)\cong\otimes_{\alpha\in\Delta^+}\cU(\fg_{-\alpha}),\ \ \ \ 
\cU(\fn^+)\cong\otimes_{\alpha\in\Delta^+}\cU(\fg_{\alpha})$$
where $\otimes$ stands for the graded tensor product.
Then $\cU(\fn^+)^{\#}$ 
is isomorphic to the maximal $\fh$-invariant submodule
of the graded tensor product 
$\Hom(\otimes_{\alpha\in\Delta^+}\cU(\fg_{\alpha}),\mathbb{C})=
\otimes_{\alpha\in\Delta^+}\Hom(\cU(\fg_{\alpha}),\mathbb{C})$
that is 
$$\cU(\fn^+)^{\#}\cong \otimes_{\alpha\in\Delta^+}\cU(\fg_{\alpha})^{\#}.$$
Now the existence of a non-graded isomorphism $\cU(\fn^+)\cong \cU(\fn^+)^{\#}$
follows from~\ref{negr}.

\subsubsection{}
\begin{rem}{Xalpha}
Let us explain why an isomorphism between $\cU(\fg_{-\alpha})$ and 
$\cU(\fg_{\alpha})^{\#}$ is non-graded.
As $\fh$-modules, $\cU(\fg_{-\alpha})\cong \sum_{r=0}^{\infty} 
\cS^r(\fg_{-\alpha})$
and $\cU(\fg_{\alpha})^{\#}=\sum_{r=0}^{\infty} \cS^r(\fg_{\alpha}^*)$
where $\cS^r$ stands for the $r$th symmetric power.
It is easy to check that there exists an isomorphism
$\psi_1:\fg_{-\alpha}\iso \Pi(\fg_{\alpha}^*)$.
This induces an isomorphism
 $\psi_r: \cS^r(\fg_{-\alpha})\iso \Pi\bigl(\cS^r(\fg_{\alpha}^*)\bigr)$
for all $r>0$. 
On the other hand, $\cS^0(\fg_{-\alpha})\not\cong\Pi\bigl(
\cS^0(\fg_{\alpha}^*)\bigr)$ since 
$\cS^0(\fg_{-\alpha}), \cS^0(\fg_{\alpha}^*)$ 
are trivial even $\fh$-modules. Therefore
$\sum_{r>0} \psi_r: \sum_{r=1}^{\infty} 
\cS^r(\fg_{-\alpha})\iso \Pi\bigl(\sum_{r=1}^{\infty} 
\cS^r(\fg_{\alpha}^*)\bigr)$
can not be extended to an isomorphism $\cU(\fg_{-\alpha})\to\Pi
\bigl(\cU(\fg_{\alpha})^{\#}\bigr)$.
\end{rem}

\section{The case $\fg=\fsq(2)$}
\label{sectq2}
In this section we write down the Shapovalov
matrices for the simplest case $\fg:=\fsq(2)$.

The simple modules $V(\lambda)=N(\lambda)/\ol{N(\lambda)}$ 
for $\fsq(2)$ were described in~\cite{pe}.
We explicitly  calculate the kernels of Shapovalov maps
which gives both $\ol{M(\lambda)}$ and $\ol{N(\lambda)}$.

The even part $\fg_{\ol{0}}=\fgl(2)$ has the only root $\alpha$;
the elements $e_{\alpha},f_{\alpha},h_{\alpha}$ form the standards
basis of $\fsl(2)\subset\fgl(2)$ and $h':=h_{\ol{\alpha}}$ 
is a central element in $\fgl(2)$. 
The space $\fh_{\ol{1}}$ is spanned by an element $H$ satisfying $H^2=h'$.

We will omit the lower index $\alpha$ 
writing $f$ instead of $f_{\alpha}$ and
etc.

\subsection{Verma and Weyl modules}
Observe that $\Uh$ has rank two over $\Sh$ and $c(\lambda)=1$ if 
$h'(\lambda)=0$, $c(\lambda)=0$ if $h'(\lambda)\not=0$.
The superalgebra $\Cl(\lambda)$ 
is a simple module over itself iff $c(\lambda)=0$.
Therefore $M(\lambda)=N(\lambda)\cong \Pi(N(\lambda))$ 
if $h'(\lambda)\not=0$. If $\lambda$ is such that $h'(\lambda)=0$
then $M(\lambda)$ has a 
submodule $N(\lambda)$ 
(having an odd highest weight vector) whose quotient is $\Pi(N(\lambda))$.
Notice that for $n>0$ 
the weight space $\cU(\fn^-)_{n\alpha}$ is two-dimensional.
This gives the following formulae for non-graded characters
$$\ch M(\lambda)=e^{\lambda}(2+4(e^{-\alpha}+e^{-2\alpha}+\ldots))$$
and
$$
\ch N(\lambda)=\left\{
\begin{array}{ll} e^{\lambda}(1+2(e^{-\alpha}+e^{-2\alpha}+\ldots)), 
& h'(\lambda)=0,\\
e^{\lambda}(2+4(e^{-\alpha}+e^{-2\alpha}+\ldots))& h'(\lambda)\not=0.
\end{array}\right.$$

\subsection{Shapovalov matrices}
\label{Salpha}
Set $R:=\Uh, A:=\Sh$. 
The elements $1,H$ form a free $A$-basis of $R$; 
one has $\int\! H=1,\int\! 1=0$. 

The weight space  $\cU(\fb^-)_0$ coincides with $R$.
The matrix $B_0$ written with respect to the basis $1,H$ takes form
$$B_0=
\begin{pmatrix}
\ \ 0\ & 1\\
-1\ &0\\
\end{pmatrix}.$$

For $k\geq 0$ the elements $f^{k+1},f^kFH,f^{k+1}H,f^kF$ 
form a free $A$-basis of $\cU(\fb^-)_{(k+1)\alpha}$.
Using the formulas of~\ref{q2ltrm} we see that the matrix $B_{(k+1)\alpha}$
written in this basis takes form
$$B_{(k+1)\alpha}=
\begin{pmatrix}
\ \ 0\ &0\ &(-1)^{k+1}(k+1)hh_k\ & (-1)^{k+1}(k+1)h_k\\
\ \ 0\ &0\ &(-1)^k(k+1)h'h_k\ & (-1)^kh'h_k\\
(-1)^k(k+1)hh_k\ & (-1)^k(k+1)h'h_k\ & 0\ &0\\
(-1)^{k+1}(k+1)h_k\ & (-1)^{k+1}h'h_k\ & 0\ &0
\end{pmatrix}$$
where 
$$h_0:=1,\ \ \ h_k:=k!(h-1)\ldots(h-k) \text{ for }k>0.$$

Therefore $\det B_{(k+1)\alpha}=(h'h_kh_{k+1})^2$ 
and so $\Norm S_{(k+1)\alpha}=\mathbb{C}^*h'h_kh_{k+1}$ that is
$$\Norm S_{(k+1)\alpha}=h_{\ol{\alpha}}(h_{\alpha}-1)^2\ldots(h_{\alpha}-k)^2
(h_{\alpha}-(k+1)).$$

\subsection{Structure of $M(\lambda)$ as $\fg_{\ol{0}}$-module}
\label{fdim}
Writing $\lambda=(x,a)$ where $x:=h(\lambda), a:=h'(\lambda)$, we see
that $N(\lambda)$ is not simple iff $x\in\mathbb{Z}_{>0}$ or $a=0$.
Let us describe $V(\lambda)=N(\lambda)/\ol{N(\lambda)}$ 
as an $\fsl(2)$-module.
Denote by $M_{\fsl(2)}(x)$ (resp., $V_{\fsl(2)}(x)$) a Verma (resp., simple)
$\fsl(2)$-module of the highest weight $x$.

Recall that the evaluated
Shapovalov matrices $B_{\nu}(\lambda)$ correspond to the restriction
of Shapovalov map $S(\lambda)$ to $M(\lambda)_{\lambda-\nu}$.
Let $E(\lambda)$ be a simple $\Cl(\lambda)$-module.
By~\Claim{olInd}, $\ol{N(\lambda)}$ coincides with the kernel
of Shapovalov map $\Xi(E(\lambda))$. Since
$E(\lambda)$ is a submodule of $\Cl(\lambda)$,
$N(\lambda)$ is a submodule of 
$M(\lambda)$
and $\Xi(E(\lambda))$ is the restriction of $S(\lambda)$
to $N(\lambda)$.

Identify $M(\lambda)$ with $\cU(\fb^-)$. The submodule
$N(\lambda)$ is given by the formula
$$N((x,a))=\left\{\begin{array}{ll}
\cU(\fb^-), & a\not=0,\\
\cU(\fn^-)H, & a=0.
\end{array}\right.$$

\subsubsection{}
\label{craq2}
Take  $\lambda=(m,a)$ such that $a\not=0$ and $m\in\mathbb{Z}_{>0}$.
The Shapovalov matrix $B_{k\alpha}$ has rank $4$ if $0<k<m$,
$B_{m\alpha}$ has rank $2$ and $B_{k\alpha}$
has rank $0$ for $k>m$. Therefore
$V(1,a)=V_{\fsl(2)}(1)\oplus V_{\fsl(2)}(1)$
and 
$
V(m,a)=V_{\fsl(2)}(m)^{\oplus 2}\oplus V_{\fsl(2)}(m-2)^
{\oplus 2}\ \ \text{for } m>1.
$

\subsubsection{}
Observe that
$$N((x,0))\cong M_{\fsl(2)}(x)\oplus M_{\fsl(2)}(x-2)\ \text{ for } x\not=0.$$
Take $\lambda=(x,0)$ where $x\not\in\mathbb{Z}_{\geq 0}$. 
The Shapovalov matrices
$B_{(k+1)\alpha}$ ($k\geq 0$) have rank two and 
$\Ker B_{(k+1)\alpha}$ is spanned by $f^kFH,f^k(Fh-fH)$.
Therefore  $\Ker B_{(k+1)\alpha}\cap \cU(\fn^-)H$ 
is spanned by $f^kFH$ and thus
$$\ol{M}((x,0))\cong M_{\fsl(2)}(x-2), \ \ 
V((x,0))\cong M_{\fsl(2)}(x)\cong V_{\fsl(2)}(x)
\ \ \ \text{ for } x\not\in \mathbb{Z}_{\geq 0}.$$

\subsubsection{}
Take $\lambda=(m,0)$ where $m\in\mathbb{Z}_{>0}$.
For $0\leq k <m$ the Shapovalov matrix
$B_{(k+1)\alpha}$  has rank two and 
$\Ker B_{(k+1)\alpha}$ is spanned by $f^kFH,f^k(Fh-fH)$;
thus $\Ker B_{(k+1)\alpha}\cap \cU(\fn^-)H$ 
is spanned by $f^kFH$. 
The Shapovalov matrices
$B_{(k+1)\alpha}$ ($k\geq m$) are equal to zero. Therefore
$$\ol{M}((m,0))\cong M_{\fsl(2)}(m-2)\oplus M_{\fsl(2)}(-m-2), \ \ \ 
V((m,0))\cong V_{\fsl(2)}(m).$$

\subsubsection{}
\label{sq2cl0}
Finally, take $\lambda=(0,0)$.
The module $V(0)$ is one-dimensional.
The matrices $B_{(k+1)\alpha}$ ($k\geq 0$) have rank $2$: 
the kernels are spanned by
$f^kFH,f^{k+1}H$. Thus $\ol{N(0)}=\ol{M(0)}$.

\subsubsection{}
\label{findimq2}
Summarizing we obtain
\begin{equation}
\label{Vxa}
V((x,a))\cong\left\{\begin{array}{lll}
V_{\fsl(2)}(1)\oplus V_{\fsl(2)}(1), & x=1, & a\not=0,\\
V_{\fsl(2)}(x)^{\oplus 2}\oplus V_{\fsl(2)}(x-2)^
{\oplus 2}, & x\in\mathbb{Z}_{>0}, x\not=1, & a\not=0,\\
V_{\fsl(2)}(x), & & a=0
\end{array}\right.
\end{equation}
and $V((x,a))=N((x,a))$ for $x\not\in\mathbb{Z}_{>0}$ and
$a\not=0$.

\subsection{}
\begin{cor}{corshq2qn}
Let $\fg$ be a $Q$-type Lie superalgebra and $\alpha_i$ be a simple root.
\begin{enumerate}
\item
If $\lambda\in\fh_{\ol{0}}^*$ is such that  
$m:=\lambda(h_{\alpha_i})\in\mathbb{Z}_{>0}$
then  $N(\lambda)$ has a primitive vector
of the weight $s_{\alpha}(\lambda)=\lambda-m\alpha_i$.

\item
If $\lambda\in\fh_{\ol{0}}^*$ is such that  $\lambda(h_{\ol{\alpha}_i})=0$
then  $N(\lambda)$ has a primitive vector
of the weight $\lambda-\alpha_i$.
\end{enumerate}
\end{cor}
\begin{proof}
Fix a  simple root $\alpha:=\alpha_i$. 
The superalgebra
$\fp:=\fg_{\alpha}+\fg_{-\alpha}+[\fg_{\alpha},\fg_{-\alpha}]$
is isomorphic to $\fsq(2)$. Let $U^-({\alpha})$ 
be a subalgebra of $\Ug$ generated by $f_{\alpha},F_{\alpha},H_{\alpha}$.
Let $v\in N(\lambda)$ be a highest weight vector.
Clearly, $N':=U^-(\alpha)v$ is a Weyl $\fp$-module. For $\lambda$
satisfying the assumption (i) (resp., (ii))
the formula~(\ref{Vxa}) shows the existence of a non-zero vector
$v'\in N'_{-m\alpha}$ (resp., $v'\in N'_{-\alpha}$)
such that $E_{\alpha}v'=e_{\alpha}v'=0$.
For any simple root $\alpha_j\not=\alpha$ one has 
$E_{\alpha_j}N'=e_{\alpha_j}N'=0$
because all vectors in $N'$ have weights of the form $\lambda-r\alpha_i$
and  $N(\lambda)$ does not have vectors of weight 
$\lambda-r\alpha_i+\alpha_j$. Since
$\{E_{\beta}, e_{\beta}: \beta\in\pi\}$ generate $\fn^+$,
we conclude that $v'$ is primitive.
\end{proof}

\subsection{Useful formulas} 

In~\ref{Salpha}   we used the following formulas.

\begin{lem}{q2ltrm}
$$\begin{array}{ll}
(i) & \HC(e^kf^k)=k!h(h-1)\ldots(h-(k-1)), \\
(ii) & \HC(Ee^kf^kF)=k!h'(h-1)(h-2)\ldots(h-k), \\
(iii) & \HC(Ee^{k}f^{k+1})=\HC(e^{k+1}f^kF)=(k+1)!(h-1)\ldots(h-k)H, \\
\end{array}$$
\end{lem}
\begin{proof}
By~\Lem{easylemHC}, $\HC(e^kf^k)$ has degree
$k$. Note that the term $e^kf^k$ annihilates
$V_{\fsl(2)}(m)$ for $m=0,1,\ldots,k-1$. Thus, up to a constant,
$\HC(e^kf^k)=h(h-1)\ldots(h-(k-1))$. Let $v_1$ be a highest weight vector
of $V_{\fsl(2)}(1)$ and $v_1^{\otimes k}$ be the corresponding vector
in $V_{\fsl(2)}(1)^{\otimes k}$.
 It is easy to check that 
$e^kf^k(v_1^{\otimes k})=k!^2(v_1^{\otimes k})$. Since $v_1^{\otimes k}$
is a primitive vector of weight $k$ this gives $\HC(e^kf^k)(k)=k!^2$
and (i) follows.

 For (ii) set $a_k:=\HC(e^kf^k)$.
Modulo the left ideal
$\Ug e^2$ one can write $e^kf^k=a_k+fb_ke$ where 
$b_k\in\mathbb{C}[h]$ ($b_1=1$). Then 
$a_{k+1}=\HC(e^{k+1}f^{k+1})=ea_kf+efb_kef$ that is $h^2b_k=a_{k+1}-h\zeta(a_k)$
where $\zeta:\mathbb{C}[h]\to \mathbb{C}[h]$ is an algebra homomorphism 
given by $h\mapsto h-2$.
Therefore $b_k=k!k(h-2)\ldots(h-k)$. 
Now $\HC(Ee^kf^kF)=\HC(E(a_k+fb_ke)F)=h'\zeta(a_k)+Hb_kH=h'\zeta(a_k)+h'b_k$
which gives (ii). Finally, for (iii) one has
$\HC(e^{k+1}f^kF)=\HC(e(a_k+fb_ke)F)=H\zeta(a_k)+hb_kH=(\zeta(a_k)+hb_k)H$
and similarly $\HC(Ee^{k}f^{k+1})=(\zeta(a_k)+hb_k)H$.
\end{proof}

\section{The leading term of Shapovalov determinant}
\label{leadterm}
Recall that $\det S_{\nu}=\det B_{\nu}\in\Sh$  for each $\nu\in Q^+$.
In this section we calculate the leading term
of the polynomial $\det S_{\nu}$.

\subsection{The main result of the section}
\label{partition}
The Kostant partition function $\tau:\ Q\to \mathbb{Z}_{\geq 0}$
is defined by the formula
$$\ch {\cU}({\mathfrak n}^-)=\displaystyle\prod_{\alpha\in\Delta^+_{\ol{1}}}
(1+e^{-\alpha})\displaystyle\prod_{\alpha\in\Delta^+_{\ol{0}}}(1-e^{-\alpha})^{-1}=:
\displaystyle\sum_{\eta\in Q} \tau(\eta)e^{-\eta}.$$
Note that $\tau(Q\setminus Q^+)=0$.

\subsubsection{}
\begin{defn}{} A vector $\mathbf{k}=
\{ k_{\alpha}, k_{\ol{\alpha}}\}_{\alpha\in
\Delta^+}$ is called a {\em partition of $\nu\in Q^+$} if 
$$\nu=\displaystyle\sum_{\alpha\in \Delta^+}
(k_{\alpha}+k_{\ol{\alpha}})\alpha;\ \ \ 
k_{\alpha}\in \mathbb{Z}_{\geq 0},\ 
k_{\ol{\alpha}}\in \{0,1\}\text{ for all }\alpha\in \Delta^+.$$
\end{defn}
Denote by ${\cP}(\nu)$ the set of all partitions of $\nu$. 
Clearly, $|{\cP}(\nu)|=\tau(\nu)$.

For each $\alpha\in\Delta^+$ set
$$\tau_{{\alpha}}(\nu):=
|\{\mathbf{k}\in \cP(\nu)|\ k_{\ol{\alpha}}=0\}|.$$

In this section we will show that for suitable bases one has the

\subsubsection{}
\begin{claim}{corestimsh}
The leading term of the polynomial $\Norm S_{\nu}$ is equal to
$$\prod_{\alpha\in\Delta^+}h_{\alpha}^
{\sum_{m=1}^{\infty}\tau(\nu-m\alpha)}
h_{\ol{\alpha}}^{\tau_{{\alpha}}(\nu-\alpha)}
\bigr.$$
\end{claim}

In particular, $\det S_{\nu}\not=0$ for all $\nu\in Q^+$. 
By~\Cor{corshap}, this is
equivalent to the existence of a simple Weyl module.

\Claim{corestimsh} is proven in~\ref{outi}--\ref{66} below.

Using~\Lem{estimshap2} we obtain the following useful formula
\begin{equation}\label{kappa}
\deg\Norm S_{\nu}=\sum_{\mathbf{m}\in\cP(\nu)}|\mathbf{m}|.
\end{equation}
where $|\mathbf{m}|:=\sum_{\alpha\in\Delta^+} m_{\alpha}+m_{\ol{\alpha}}$.

\subsection{Outline of the proof}\label{outi}
We reduce a computation of leading term of Shapovalov determinants 
to the case ${\fsq}(2)$.  This is done in several steps described below.

As it was explained in~\ref{coon}, Shapovalov determinants
for various $Q$-type algebras can be expressed via
Shapovalov determinants for $\fsq(n)$. In particular,
it is enough to prove~\Claim{corestimsh}
for $\fg:=\fsq(n)$.  The proof for $\fsq(n)$
goes as follows.
Define a filtration $F$ on $\fg$ by setting
$$F^{0}(\fg)=0,\ F^1(\fg)=\fn^-+\fn^++\fh_{\ol{1}},\ F^i(\fg)=\fg
\text{ for }i>1,\ \ \dot{\fg}:=\gr_F\fg.$$

Denote by $S$ (resp., $\dot{S}$) the Shapovalov map for $\fg$ 
(resp., for $\dot{\fg}$).
As we will show later (\ref{chuh}), $\Norm \dot{S}_{\nu}$
 is either zero or  equal to the leading term of
$\Norm {S}_{\nu}$.

Now let 
$${\fg}^{\sqcap}=\mathop{\prod}
\limits_{\alpha\in\Delta^+}\dot{\fs}^{(\alpha)}$$
where $\dot{\fs}^{(\alpha)}=\gr_F\fsq(2)$.
Observe that $[\dot{\fg},\dot{\fg}]=\dot{\fh}$ where $\dot{\fh}:=\gr_F\fh$
and that $\dot{\fh}\cong\fh$ as Lie algebras.
As a consequence, $\dot{\fg}$ is the quotient of
${\fg}^{\sqcap}\to\dot{\fg}$ by an ideal lying in
${\fh}^{\sqcap}$.
By~\Thm{clcon}, $\Norm \dot{S}_{\nu}=\psi(\Norm {S}^{\sqcap}_{\nu})$ 
where ${S}^{\sqcap}_{\nu}$ is the Shapovalov map for
${\fg}^{\sqcap}$. Hence the leading term of
$\Norm {S}_{\nu}$ is equal to $\psi(\Norm {S}^{\sqcap}_{\nu})$
provided  the latter is non-zero.
Since ${\fg}^{\sqcap}$ is the direct product of copies
of $\dot{\fsq}(2)$, a computation of
$\Norm {S}^{\sqcap}_{\nu}$
is reduced to the case $\dot{\fsq}(2)$--- see~\ref{approx}.

\subsection{The algebras $\dot{\fg}$}
\label{dotfg}
Retain notation of~\ref{outi} and
extend the filtration $F$ to $\Ug$.
Then  $\cU(\dot{\fg})$  identifies  with the associated graded of
$\Ug$; denote by $\dot{x}$ the image of $x\in\Ug$ or
$x\subset \Ug$ in $\cU(\dot{\fg})$.

The algebra $\dot{\fg}=\gr_F\fg$ admits
a decomposition
$\dot{\fg}=\dot{\fn}^-+\dot{\fh}+\dot{\fn}^+$.
Construct a Shapovalov map $\dot{S}$ via
the above decomposition.
The algebra $\cU(\dot{\fn}^{+})$ inherit the gradings by $Q^+$
and this grading leads to a decomposition
$\dot{S}=\sum \dot{S}_{\nu}$. In this way we obtain the 
Shapovalov matrices $\dot{B}_{\nu}$ and the determinants
$\det \dot{B}_{\nu}=(\Norm \dot{S}_{\nu})^{2^{\dim\fh_{\ol{1}}}}$.

Let us explain why $\Norm {S}'_{\nu}$
is either zero or  equal to the leading term of
$\Norm {S}_{\nu}$.

\subsubsection{}\label{chuh}
Recall that the Rees algebra
$\tilde{\fg}:=\oplus_{k=0}^{\infty} \vareps^kF^k(\fg)$
is a Lie superalgebra over $\mathbb{C}[\vareps]$
with the bracket induced by the bracket in $\fg$;
as a Lie superalgebra over $\mathbb{C}$ $\tilde{\fg}$
is graded with $k$th homogeneous component  $\vareps^k F^k(\fg)$.

The evaluation at $\vareps=0$ is canonically isomorphic to
$\gr_F\fg$ and the evaluation at $\vareps=1$ is canonically isomorphic to
$\fg$:
$\tilde{\fg}/\vareps\tilde{\fg}\cong\gr_F\fg=\dot{\fg}$, 
$\tilde{\fg}/(\vareps-1)\tilde{\fg}\cong\fg$.

The algebra $\tilde{\fg}$ inherits a triangular decomposition
where $\tilde{\fn}^{\pm}:=
\oplus_{k=0}^{\infty} \vareps^k (F^k(\fg)\cap\fn^{\pm})$, 
$\tilde{\fh}:=
\oplus_{k=0}^{\infty} \vareps^k (F^k(\fg)\cap\fh)$. Define
the Harish-Chandra projection $\tilde{\HC}:\cU(\tilde{\fg})\to
\cU(\tilde{\fh})$ and the Shapovalov map $\tilde{S}$.
The algebras $\cU(\tilde{\fn}^{\pm})$ inherit the gradings by $Q^{\pm}$
and these gradings leads to a decomposition
$\tilde{S}=\sum \tilde{S}_{\nu}$. 
In this way we obtain the 
Shapovalov matrices $\tilde{B}_{\nu}$ and the determinants
$\det \tilde{B}_{\nu}=(\Norm \tilde{S}_{\nu})^{2^{\dim\fh_{\ol{1}}}}$.

The evaluation of $\tilde{B}_{\nu}$
at $\vareps=0$ (resp., at $\vareps=1$) is the Shapovalov matrix 
$\dot{B}_{\nu}$ for $\dot{\fg}$ (resp., $B_{\nu}$ for $\fg$).

\subsubsection{}
For each $h\in\fh_{\ol{0}}$ one has
$\tilde{h}:=\vareps^2h\in\tilde{\fh}$; the reduction modulo
$\vareps$ identifies  $\tilde{h}$ with $\dot{h}\in\dot{\fh}$
and the reduction modulo $\vareps-1$ identifies  $\tilde{h}$ with $h$.
Observe that the entries of Shapovalov matrices are polynomials
homogeneous in $\vareps$.
As a consequence, the leading term of $B_{\nu}$ is equal to
 $\det \dot{B}_{\nu}$ providing that the latter is non-zero
via the obvious identification $h\mapsto \dot{h}$ of 
$\fh_{\ol{0}}$ with $\dot{\fh}_{\ol{0}}$.

For instance, for $\fg=\fsl(2)$, one has
$\tilde{B}_{2\alpha}=2(\vareps^2 h)^2-2\vareps^4 h=2\tilde{h}^2-
2\vareps^2\tilde{h}$; the identification 
$\tilde{\fg}/\vareps\tilde{\fg}$ with $\fg$ gives
$B_{2\alpha}=2h^2-2h$ and
the identification $\tilde{\fg}/(\vareps-1)\tilde{\fg}$ with $\dot{\fg}$
gives $\dot{B}_{2\alpha}=2\dot{h}^2$.

\subsection{The algebra ${\fg}^{\sqcap}$}\label{q2tensor}
For each $\alpha$ in $\Delta^+$
let $\fs^{(\alpha)}$ be a subalgebra of $\fg$ generated
by $\fg_{\pm\alpha}$ and let $\dot{\fs}^{(\alpha)}$ be its
image in $\dot\fg$.
Clearly, $\fs^{(\alpha)}\cong\fsq(2)$.

Consider a Lie algebra ${\fg}^{\sqcap}:=\mathop{\prod}
\limits_{\alpha\in\Delta^+}\dot{\fs}^{(\alpha)}$
(the direct product of Lie superalgebras). 
Set $\fh^{\sqcap}:=\prod_{\alpha\in\Delta^+}\dot{\fh}^{(\alpha)}$
where $\dot{\fh}^{(\alpha)}\subset \dot{\fs}^{(\alpha)}$ is the image of 
$\fh\cap\fs^{(\alpha)}$.

\subsubsection{}
Consider the triangular-type decomposition 
${\fg}^{\sqcap}=({\fn}^{\sqcap})^-+\fh^{\sqcap}+({\fn}^{\sqcap})^+$
where 
$({\fn}^{\sqcap})^{\pm}:=\prod_{\alpha\in\Delta^+}\fg_{\pm\alpha}$.
Let $S^{\sqcap}$ be the Shapovalov map; write 
$S^{\sqcap}=\sum_{\nu\in Q^+} S^{\sqcap}_{\nu}$ as above.
In~\ref{approx} we will prove the following formula
\begin{equation}
\label{formdotS}
\Norm S^{\sqcap}_{\nu}=\prod_{\alpha\in\Delta^+}\dot{h}_{\alpha}^
{\sum_{m=1}^{\infty}\tau(\nu-m\alpha)}
\dot{h}_{\ol{\alpha}}^{\tau_{{\alpha}}(\nu-\alpha)}
\end{equation}
where $\dot{h}_{\alpha}, \dot{h}_{\ol{\alpha}}$ stand for the images 
of these elements in $\dot{\fs}^{(\alpha)}\subset {\fg}^{\sqcap}$.

\subsubsection{}\label{psii}
Recall that $\fs^{(\alpha)}\cong\fsq(2)$
has a $(2|1)$-dimensional Cartan algebra with a basis
$h_{\alpha},h_{\ol{\alpha}},H_{\alpha}$.
Observe that $({\fn}^{\sqcap})^{\pm}\cong \dot{\fn}^{\pm}$ 
since they are commutative
Lie superalgebras. This implies
$\dot{\fg}\cong {\fg}^{\sqcap}/I$ where
$$I:=\spn\{\dot{h}_{\alpha+\beta}-\dot{h}_{\alpha}-\dot{h}_{\beta},\ 
\dot{h}_{\ol{\alpha+\beta}}-\dot{h}_{\ol{\alpha}}-\dot{h}_{\ol{\beta}},\ 
\dot{H}_{\alpha+\beta}-\dot{H}_{\alpha}-\dot{H}_{\beta}|\ 
\alpha,\beta,\alpha+\beta\in\Delta^+\}.$$

Now using~\Thm{clcon}, we obtain~\Claim{corestimsh} from
the formula~(\ref{formdotS}).

\subsection{Proof of~(\ref{formdotS})}
\label{approx}
Define the equivalence relation on $\cP(\nu)$ by setting
$$\mathbf{m}\approx\mathbf{m}'\ \Longleftrightarrow\ 
\forall \alpha\in\Delta^+\ \ 
m_{\alpha}+m_{\ol{\alpha}}=m'_{\alpha}+m'_{\ol{\alpha}}.$$
Denote by $\ol{\cP}(\nu)$ the set of of equivalence
classes in $\cP(\nu)$. For $\mathbf{k}\in\ol{\cP}(\nu)$
set $k_{\alpha}=m_{\alpha}+m_{\ol{\alpha}}$ where
$\mathbf{m}\in\cP(\nu)$ belongs to the class $\mathbf{k}$.
For $\mathbf{m}\in\cP(\nu)$ set
$$\supp\mathbf{m}:=
\{\alpha\in\Delta^+:\ m_{\alpha}+m_{\ol{\alpha}}\not=0\}.$$
Then $\supp$ is well-defined for $\mathbf{k}\in\ol{\cP}(\nu)$ and
$\supp\mathbf{k}:=\{\alpha\in\Delta^+:\ k_{\alpha}\not=0\}$.

\subsubsection{}
Recall that ${\fg}^{\sqcap}:=
\prod_{\alpha\in\Delta^+}\dot{\fs}^{(\alpha)}$. Let
$S^{(\alpha)}:\cM^{(\alpha)}\to (\cN)^{(\alpha)}$ be the Shapovalov map for
$\dot{\fs}^{(\alpha)}$. It
is easy to see that $\cM^{\sqcap}=\otimes_{\alpha\in\Delta^+}
\cM^{(\alpha)}$, $({\cN})^{\sqcap}=\otimes_{\alpha\in\Delta^+}
(\cN)^{(\alpha)}$ and $S^{\sqcap}=\otimes_{\alpha\in\Delta^+}S^{(\alpha)}$.
Choose  the integral $\int$ on $\fg^{\sqcap}$ to be
$\otimes_{\alpha\in\Delta^+} \int^{(\alpha)}$ where  $\int^{(\alpha)}$
is the integral for $\dot{\fs}^{(\alpha)}$. Recall that
$B$ is the composition of $S$ and a map induced by $\int$. We obtain
$B^{\sqcap}=\otimes_{\alpha\in\Delta^+}B^{(\alpha)}$
where $B^{\sqcap}$ (resp., $B^{(\alpha)}$)
is the map $B$ for ${\fg}^{\sqcap}$ (resp., for $\dot{\fs}^{(\alpha)}$).

Recall that $\cM^{\sqcap}$ (resp., $\cM^{(\alpha)}$) is a free module
over a polynomial algebra $\cU(\fh^{\sqcap}_{\ol{0}})$
(resp., over $\cU(\dot{\fh}^{(\alpha)})$). If $L$ is 
a free submodule of $\cM^{\sqcap}$ (resp., of  $\cM^{(\alpha)}$) we will use
the notation ``$\rnk L$'' for the rank over this polynomial algebra.

\subsubsection{}
One has
$$\cM^{\sqcap}_{-\nu}=\oplus_{\mathbf{k}\in \ol{\cP}(\nu)}
\cM^{\sqcap}_{\mathbf{k}},\
\text{ where }
\cM^{\sqcap}_{\mathbf{k}}:=\otimes_{\alpha\in\Delta^+}
\cM^{(\alpha)}_{-k_{\alpha}\alpha}$$
with a similar formula  for
$({\cN})^{\sqcap}_{-\nu}$. Let  $B^{\sqcap}_{\mathbf{k}}$ be
the restriction of $B^{\sqcap}$ to $\cM^{\sqcap}_{\mathbf{k}}$.
One has
\begin{equation}\label{Bak}
B^{\sqcap}_{\mathbf{k}}=\otimes_{\alpha\in\Delta^+}B^{(\alpha)}_
{-k_{\alpha}\alpha}
\end{equation}
where $B^{(\alpha)}_{r\alpha}$ is the restriction of 
$B^{(\alpha)}$ to ${\cM}^{(\alpha)}_{-r\alpha}$.
Since $B^{(\alpha)}$ maps
$\cM^{(\alpha)}_{-r\alpha}$ to $(\cN)^{(\alpha)}_{-r\alpha}$, 
$B^{\sqcap}_{\mathbf{k}}$ maps $\cM^{\sqcap}_{\mathbf{k}}$ to 
${\cM^{\sqcap}}^{\#}_{\mathbf{k}}$. As a consequence,
\begin{equation}\label{Dak}
\det B^{\sqcap}_{\nu}=\prod_{\mathbf{k}\in \ol{\cP}(\nu)}
\det B^{\sqcap}_{\mathbf{k}}.\end{equation}

\subsubsection{}
Fix $\mathbf{k}\in \ol{\cP}(\nu)$ and
let us compute $\det B^{\sqcap}_{\mathbf{k}}$ 
using the decomposition~(\ref{Bak}).
Recall that for $\psi_i\in\End(V_i)$ one has
$\det(\otimes\psi_i)=\prod(\det\psi_i)^{n/n_i}$ 
where $n_i:=\dim V_i, n:=\prod n_i=\dim \otimes V_i$. 
The module $\cM^{(\alpha)}_{-r\alpha}$ has rank $2$ for $r=0$  and 
rank $4$ for $r>0$ (see~\ref{Salpha}).
Therefore
$$\begin{array}{rl}
\det B^{\sqcap}_{\mathbf{k}}&=\prod_{\alpha\in\supp\mathbf{k}} 
(\det B^{(\alpha)}_{k_{\alpha}\alpha})^{r(\mathbf{k})/4},\\
&\ \ \ \ \
 \text{where } r(\mathbf{k}):=\rnk\cM^{\sqcap}_{\mathbf{k}}.\end{array}
$$
Recall that $\dot{\fs}^{(\alpha)}=\gr_F\fsq(2)$. By~\ref{dotfg}, the entries
of $B^{(\alpha)}_{(k+1)\alpha}$ are the leading terms
of the entries of Shapovalov matrix $B_{k\alpha}$ for $\fsq(2)$.
Shapovalov matrices $B_{k\alpha}$ were computed in~\ref{Salpha}; 
using the formulas there we get
$$\det B^{(\alpha)}_0=1,\ \ \det B^{(\alpha)}_{k_{\alpha}\alpha}=
(k_{\alpha}-1)!^2 \dot{h}_{\ol{\alpha}}^2\dot{h}_{\alpha}^{4k_{\alpha}-2}
\text{  for } k_{\alpha}>0.$$

Substituting our formulas in~(\ref{Dak}) we obtain, up to a non-zero scalar,
\begin{equation}
\label{d(alpha)}
\begin{array}{rl}
\det B^{\sqcap}_{\nu}&=\mathop{\prod}\limits_{\mathbf{k}\in \ol{\cP}(\nu)}
\mathop{\prod}\limits_{\alpha\in\Delta^+}  
\dot{h}_{\ol{\alpha}}^{d(\ol{\alpha})}\dot{h}_{\alpha}^{d(\alpha)},\\
& \text{where }
d(\ol{\alpha}):={1\over 2}
\sum_{\mathbf{k}\in \ol{\cP}(\nu): \alpha\in\supp\mathbf{k}}
\rnk\cM^{\sqcap}_{\mathbf{k}},\\
& \ \ \ \ \ \ \ \ \  
d(\alpha):=\sum_{\mathbf{k}\in \ol{\cP}(\nu):: \alpha\in\supp\mathbf{k}}
\frac{4k_{\alpha}-2}{4}\rnk\cM^{\sqcap}_{\mathbf{k}}.
\end{array}
\end{equation}

We will simplify the expressions for $d(\ol{\alpha}),d(\alpha)$
in~\ref{jumps} below.

\subsection{The multiplicities $d(\alpha),d(\ol{\alpha})$}
\label{jumps}
\subsubsection{}
For a partition $\mathbf{m}\in {\cP}(\nu)$ view
$\mathbf{f}^{\mathbf{m}}:=
\prod f_{\alpha}^{m_{\alpha}}F_{\alpha}^{m_{\ol{\alpha}}}$
as an element of $\fg^{\sqcap}$. Identify $\cM^{\sqcap}$
with $\cU({\fh}^{\sqcap}+({\fn}^{\sqcap})^-)$ and for each 
$\mathbf{m}\in {\cP}(\nu)$ denote by $\cM^{\sqcap}_{\mathbf{m}}$
the space $\cU({\fh}^{\sqcap})\mathbf{f}^{\mathbf{m}}\subset \cM^{\sqcap}$.
Recall that $\ol{\cP}(\nu)={\cP}(\nu)/\approx$.
For any ${\mathbf{k}}\in\ol{\cP}(\nu)$ one has
$$\cM^{\sqcap}_{{\mathbf{k}}}=\mathop{\oplus}
\limits_{\mathbf{m}\in {\cP}(\nu),\ 
\mathbf{m}\in {\mathbf{k}}}
\cM^{\sqcap}_{\mathbf{m}}.$$ 
Recall that $\cU({\fh}^{\sqcap})$
is a Clifford algebra 
over the polynomial algebra $\cU({\fh}^{\sqcap})_{\ol{0}})$
and so
$\cM^{\sqcap}_{\mathbf{m}}$ is a free module of rank $2^N$
over this algebra where $N:=\dim {\fh}^{\sqcap}_{\ol{1}}=|\Delta^+|$.

\subsubsection{}
Let us simplify the expression
$d(\ol{\alpha})={1\over 2}\rnk\sum_{\mathbf{k}\in \ol{\cP}(\nu): 
\alpha\in \supp \mathbf{k}}
\cM^{\sqcap}_{\mathbf{k}}$ obtained in~(\ref{d(alpha)}).
One has
$$d(\ol{\alpha})=2^{N-1}
|\{\mathbf{m}\in {\cP}(\nu): m_{\alpha}+m_{\ol{\alpha}}\not=0\}|
$$
because
$\sum_{\mathbf{k}\in \ol{\cP}(\nu): 
\alpha\in \supp \mathbf{k}}\cM^{\sqcap}_{\mathbf{k}}=
\sum_{\mathbf{m}\in {\cP}(\nu): \alpha\in \supp \mathbf{m}}
\cM^{\sqcap}_{\mathbf{m}}$. Observe that
$$|\{\mathbf{m}\in {\cP}(\nu): m_{\alpha}+m_{\ol{\alpha}}\not=0\}|=
2|\{\mathbf{m}\in {\cP}(\nu):m_{\ol{\alpha}}\not=0\}|=
2\sum_{\mathbf{m}\in {\cP}(\nu)}m_{\ol{\alpha}}.$$
Using~\Lem{estimshap2} (ii) we obtain
\begin{equation}
\label{dol}
d(\ol{\alpha})=2^{N}\sum_{\mathbf{m}\in {\cP}(\nu)}m_{\ol{\alpha}}=
2^{N}\tau_{{\alpha}}(\nu-\alpha).
\end{equation}

\subsubsection{}
Let us simplify the term $d(\alpha)$. One can rewrite~(\ref{d(alpha)}) as
$$d(\alpha)+d(\ol{\alpha})=\sum_{\mathbf{k}\in \ol{\cP}(\nu)}
k_{\alpha}\rnk\cM^{\sqcap}_{\mathbf{k}}=
\sum_{\mathbf{m}\in {\cP}(\nu)}
(m_{\alpha}+m_{\ol{\alpha}})\rnk\cM^{\sqcap}_{\mathbf{m}}=
2^N\sum_{\mathbf{m}\in {\cP}(\nu)}(m_{\alpha}+m_{\ol{\alpha}}).$$
Using~\Lem{estimshap2} (i) we get
\begin{equation}
\label{dal}
d(\alpha)=2^N\sum_{r=1}^{\infty}\tau(\nu-r\alpha).
\end{equation}

\subsubsection{}
Finally, recalling that $N=\dim\fh^{\sqcap}_{\ol{1}}$ and
substituting~(\ref{dol}),\ref{dal}) into~(\ref{d(alpha)})
we obtain the formula~(\ref{formdotS}).

\subsubsection{}\label{66}
Retain notation of~\ref{partition}.

\begin{lem}{estimshap2}
For all $\alpha\in {\Delta}^+$ one has 
\begin{enumerate}
\item 
$\ \ \ \displaystyle\sum_{\mathbf{k}\in \cP(\nu)}k_{\alpha}=
\displaystyle\sum_{m=1}^{\infty}\tau(\nu-m\alpha);$
\item 
$\ \ \ \displaystyle\sum_{\mathbf{k}\in \cP(\nu)}k_{\ol{\alpha}}=
\tau(\nu)-\tau_{\alpha}(\nu)=\tau_{{\alpha}}(\nu-\alpha)$
\end{enumerate}
\end{lem}
\begin{proof}  For $\nu\not\geq\alpha$ the assertions obviously 
hold since both sides
of both equations are equal to zero. 
Fix $\nu\geq\alpha$ and assume that (i) holds for all $\mu<\nu$.
The map $\mathbf{k}\mapsto (\mathbf{k}-\alpha)$ 
gives a bijection from the set $\{\mathbf{k}\in \cP(\nu)|\ 
k_{\alpha}\not=0\}$ onto $\cP(\nu-\alpha)$. Therefore
$$\begin{array}{lcl}
\displaystyle\sum_{\mathbf{k}\in \cP(\nu)}k_{\alpha}
&=&\displaystyle\sum_{\mathbf{k}\in \cP(\nu-\alpha)}(k_{\alpha}+1)
=\displaystyle\sum_{m=1}^{\infty}\tau(\nu-\alpha-m\alpha)
+|\cP(\nu-\alpha)|\\
&=&\displaystyle\sum_{m=2}^{\infty}\tau(\nu-m\alpha)+
\tau(\nu-\alpha)=\displaystyle\sum_{m=1}^{\infty}\tau(\nu-m\alpha)
\end{array}$$
and (i) follows.
The map $\mathbf{k}\mapsto
(\mathbf{k}-\alpha)$ gives a bijection 
$$\{\mathbf{k}\in \cP(\nu)|\ 
k_{\ol{\alpha}}=1\}\iso 
\{\mathbf{k}\in \cP(\nu-\alpha)|\ k_{\ol{\alpha}}=0\}.$$
This gives (ii).
\end{proof}

\section{Jantzen filtration}
\label{sectjan}
The Jantzen filtration and sum formula were described by Jantzen in~\cite{ja}
for a Verma module over a semisimple Lie algebra.
In this section we adapt this construction  for $Q$-type superalgebras.

Throughout this section we assume that all Shapovalov determinants
are non-zero polynomials that is $\det B_{\nu}\not=0$
for all $\nu$.

We define the Jantzen filtration in~\ref{constjan}.

\subsection{Main properties of Jantzen filtration.}
\label{notjan}
The construction of Jantzen filtration depends on a ``generic vector''
$\rho'\in\fh_{\ol{0}}^*$ 
 satisfying the following property

{\em (J1) the hypersurfaces } $\det B_{\nu}=0$ 
{\em do not contain straight lines parallel to } $\rho'$.

In other words, $\forall\lambda\in\fh_{\ol{0}}^*\ \ \exists c\in\mathbb{C}\
\text{such that} \det B_{\nu}(\lambda+c\rho')\not=0$.
It is not hard to show that the condition $\det B_{\nu}\not=0$
for all $\nu$ ensures the existence of $\rho'$ satisfying (a).
For semisimple Lie algebras one can take $\rho':=\rho$.

For each $\lambda\in\fh_{\ol{0}}^*$
the Jantzen filtration $\{F^r(M(\lambda))\}$ is a finite decreasing filtration
by $\fg$-submodule; one has
$$F^0(M(\lambda))=M(\lambda),\ \ \ F^1(M(\lambda))=\ol{M(\lambda)}.$$

\subsubsection{}
Define the order $l$ of zero of a polynomial $q\in\Sh$ at $\mu\in\fh_{\ol{0}}^*$ 
as follows:
$l=0$ if $q(\mu)\not=0$; $l=1$ if $q(\mu)=0$ and there exists
a non-zero partial derivative
${\frac{\partial q}{\partial x}}(\mu)\not=0$ and so on. 

Denote by $m_{\nu}(\lambda,\rho')$  the order of zero of 
$\det B_{\nu}(\lambda+x\rho')\in\mathbb{C}[x]$ at $x=0$ and
by $m_{\nu}(\lambda)$  the  order of zero of $\det B_{\nu}$. 
Let $\Gamma$ be the set of irreducible components of 
hypersurfaces $\det B_{\nu}=0$
and $d_{\gamma}(\nu)$ be the multiplicity of $\gamma$ in $\det B_{\nu}=0$
(thus $d_{\gamma}(\nu)$ is a non-negative integer). Then, if all
$\gamma\in\Gamma$ are smooth at $\lambda$ one has
$$m_{\nu}(\lambda)=\sum_{\gamma\in\Gamma:\lambda\in\gamma} d_{\gamma}(\nu).$$
Note that the formula remains valid if $\lambda\not\in\gamma$ for
$\gamma\in\Gamma$. 

Clearly,
$$\corank B_{\nu}(\lambda)\leq m_{\nu}(\lambda)\leq m_{\nu}(\lambda,\rho').$$
Moreover
$$m_{\nu}(\lambda,\rho')=m_{\nu}(\lambda)$$ 
if the vector $\rho'$
is transversal to the hypersurface $\det B_{\nu}=0$ at point $\lambda$ i.e.,
it is transversal to all irreducible components passing through $\lambda$.

\subsubsection{}
The following property is proven in~\cite{ja} 
\begin{equation}\label{mlambdarho}
m_{\nu}(\lambda,\rho')=\sum_{r=1}^{\infty} \dim
F^r(M(\lambda))_{\lambda-\nu}.
\end{equation}

Assume that $\lambda$ is such that

{\em (J2) the vector }$\rho'$
{\em is transversal to the hypersurfaces }$\det B_{\nu}=0$
{\em at point} $\lambda$.

Then the following ``sum formula'' holds
\begin{equation}\label{mlambda}
\sum_{\gamma\in\Gamma:\lambda\in\gamma} d_{\gamma}(\nu)=
\sum_{r=1}^{\infty} \dim F^r(M(\lambda))_{\lambda-\nu}.
\end{equation}

It is not always possible to find $\rho'$ such that (J2) holds
for all $\lambda$:
for instance, for $\det B_{\nu}=h_1^2+h_2^2-1$ there is no
$\rho'$ with this property. 
However, one can always choose $\rho'$ 
transversal to the hypersurface $\det B_{\nu}=0$ at all points
$\lambda\in\fh_{\ol{0}}^*\setminus X$
where $X$ is a set of codimension two. Thus the sum formula
holds for a generic point of each hypersurface $\det B_{\nu}=0$.

\subsubsection{}
\begin{rem}{ordzero}
Since $F^1(M(\lambda))=\ol{M(\lambda)}$ one has
$\corank B_{\nu}(\lambda)=\dim F^1(M(\lambda))_{\lambda-\nu}$.
Using~(\ref{mlambdarho}) we get
$$F^2(M(\lambda))=0\ \ \Longleftrightarrow\ \ 
m_{\nu}(\lambda,\rho)=\corank B_{\nu}(\lambda)\ \ \forall \nu\in Q^+.$$
\end{rem}

\subsubsection{}
Later on we show that 
the Shapovalov determinants for $Q$-type Lie superalgebras
admit linear factorizations. In other words,
 the union of
the hypersurfaces $\det B_{\nu}=0$ is the union of 
countably many hyperplanes.
We choose $\rho'$ which is not parallel to these hyperplanes.
Clearly, $\rho'$ satisfies the condition (J2) for all $\lambda$. 
Hence the sum formula~(\ref{mlambda})
holds for all $\lambda\in\fh_{\ol{0}}^*$.

\subsubsection{}
Interesting questions are whether the Jantzen filtration depends on $\rho'$
and whether for given $\rho'$ the filtration induces a unique filtration
of $N(\lambda)$. For semisimple Lie algebras the Jantzen filtration
does not depend on a ``generic vector'', see~\cite{bb}, 5.3.1.

\subsection{A construction of Jantzen filtration.}
\label{constjan}
Let $x$ be an indeterminate, $L$ be the local ring ${{\mathbb C}[x]}_{(x)}$, 
and $F$ be its field of fractions. Endow $L$ and $F$ with the trivial 
$\mathbb{Z}_2$-grading:
$L_{\overline{1}}=F_{\overline{1}}=0$.

We shall extend the scalars of "our favorite objects" from $\mathbb{C}$ 
to $L$ and to $F$.
For a $\mathbb{C}$-vector superspace $V$ denote by $V_L$ the
$L$-module $L\otimes V$ and by $V_F$ the $F$-vector superspace
$F\otimes V$.
We identify ${\cU}({\mathfrak g}_L)$ with the superalgebra
${\cU}({\mathfrak g})_L$. Retain notation of~\ref{cCO} and
define the Shapovalov map $\Xi_L$.
Clearly, $\cM_L=\Ind(R_L)$. 
For any $\mu\in(\fh_{\ol{0}}^*)_L$
define a $\fg_L$-$\fh_L$ bimodules
$M_L(\mu):=\Ind(\Cl(\mu))$ and $\Coind(\Cl(\mu))$
Observe that a Shapovalov matrix for $\fg_L$ written
with respect to a bases lying in $\fg$ coincides with
the Shapovalov matrix for $\fg$ written
with respect to the same bases. Consequently,
the Shapovalov determinants $\det B_{\nu}\in \Sh$ viewed
as elements of the algebra $\Sh_L$ coincide with the Shapovalov
determinants $\det B_{\nu}$ constructed for ${\cU}({\mathfrak g}_L)$.

The vector space ${\mathfrak g}_F$
admits the natural structure of a Lie superalgebra and
${\cU}({\mathfrak g}_F):={\cU}({\mathfrak g})_F$ 
is its enveloping superalgebra. For any $\mu\in (\fh_{\ol{0}}^*)_L$
the localized module $M_L(\mu)\otimes_L F$
is naturally isomorphic to the $\fg_F$-module
$M_F(\mu)$ where $\mu$ is viewed as an element of
$(\fh_{\ol{0}}^*)_F$ via the natural embedding $(\fh_{\ol{0}}^*)_L\hra (\fh_{\ol{0}}^*)_F$.

\subsubsection{}
\label{Fmnu}
Fix $\lambda\in\fh^*_{\ol{0}}$ and set 
$$M:=M_L(\lambda+x\rho'),\ \ 
M':=\Coind(\Cl(\lambda+x\rho')),\ \
\varphi:=S(\lambda+x\rho').$$
Thus $\phi: M\to M'$ is the evaluation of the Shapovalov map at 
$\lambda+x\rho'$.

Define a decreasing $\mathbb{Z}$-filtration on 
$M$ by setting
$F^r(M)=M$ for
$r\leq 0$ and 
$$F^r(M):=\{v| \varphi(v)\in (x^r)M'\}\ \text{ for } r>0.$$
Notice that each term $F^r(M)$ is a
$\fg_L$-$\fh_L$ bisubmodule because $\varphi$
is a $\fg_L$-$\fh_L$ homomorphism by~\ref{IndCl}. 
Due to the condition (J1) of~\ref{notjan}
$\det B_{\nu}(\lambda+x\rho')\not=0$ for all $\nu\in Q^+$.
Thus $\Ker\varphi=0$ and so
$$\displaystyle\bigcap_{r=0}^{\infty} F^r(M)=0.$$

\subsubsection{}
\label{janlambda}
Observe that $M(\lambda)=M/(xM)$
and 
$\varphi/(x\varphi)=S(\lambda): M(\lambda)
\to\Coind (\Cl(\lambda))$.
Let $F^r(M(\lambda))$
be the image of $F^r(M)$. We get a decreasing $\fg$-$\fh$ filtration on 
$M(\lambda)$ with the property
$$\displaystyle\bigcap_{r=0}^{\infty} F^r(M(\lambda))=0.$$

The filtration is finite since $M(\lambda)$ has a finite length.
 
One has
$$\begin{array}{l}
F^0(M(\lambda))=M(\lambda),\\
F^1(M(\lambda))=\Ker S(\lambda)=\ol{M(\lambda)}.
\end{array}$$

The filtration $F^r(M(\lambda))$ is an analogue 
of the Jantzen filtration for the module $M(\lambda)$.

\subsection{Example: The case $\fg=\fsq(2)$}
\label{jntz}
Retain notation of~\ref{fdim}. In~\ref{fdim} we described 
$\ol{N(\lambda)}$ and $\ol{M(\lambda)}=F^1(M(\lambda))$.
Below we describe the Jantzen filtration on $M(\lambda)$.

For $\lambda=(x,a)$ where $x\not\in\mathbb{Z}_{>0}$ and $a\not=0$
the module $M(\lambda)=N(\lambda)$ is simple.

For $\lambda=(m,a)$ where $m\in\mathbb{Z}_{>0}$ and $a\not=0$
the module $M(\lambda)=N(\lambda)$ has length two and its
Jantzen filtration has length two: $F^2(M(\lambda))=0$.
More precisely, $M(\lambda)$ has a unique non-trivial
submodule $N(s_{\alpha}\lambda)=M(s_{\alpha}\lambda)=F^1(M(\lambda))$.

For $\lambda=(x,0)$ where $x\not\in\mathbb{Z}_{>0}$
one has $F^2(M(\lambda))=0$.
The module $F^1(M(\lambda))$ has a basis $f^iFH,f^i(Fh-fH)$ 
with $i\geq 0$ and so $F^1(M(\lambda))
\cong M(\lambda-\alpha)/\ol{M(\lambda-\alpha)}$.
The module $N(\lambda)$ has length two:
its unique non-trivial submodule is $\ol{N(\lambda)}\cong
V(\lambda-\alpha)$. 
The module
$M(\lambda)$ has length four because it admits a submodule
isomorphic to $N(\lambda)$ with the quotient isomorphic to $\Pi(N(\lambda))$.

For $\lambda=(m,0)$ where $m\in\mathbb{Z}_{>0}$ one has
$F^3(M(\lambda))=0$. The term
$F^1(M(\lambda))$ contains the weight spaces
$M(\lambda)_{\mu}$ for $\mu<s_{\alpha}\lambda$ and the vectors
$f^iFH,f^i(Fh-fH)$. The term 
$F^2(M(\lambda))$ is spanned by $f^{m+k}FH$ ($k\geq 0$) and
thus $F^2(M(\lambda))\cong V(s_{\alpha}\lambda)$.

\subsection{The sum formulas}
\label{lemJan}
The formula(\ref{mlambdarho}) is an immediate consequence of the 
following fact proven in~\cite{ja}.
Let $L$ be the local ring ${\mathbb C}[x]_{(x)}$ and  
let $N, N'$ be free $L$-modules of a finite rank $r$.
Let $\varphi: N\to N'$ be an injective linear map.
Define a decreasing filtration on $N$ by setting
$$F^j(N):=\{v\in N|\ \varphi(v)\in (x^j)N'\}.$$
The claim is that the sum
$\sum_{j=1}^{\infty}\dim \bigl(F^j(N)/(F^j(N)\cap xN)\bigr)$
is equal to
the order of zero of  $\det D$ at the origin where
$D\in \Mat_{r}(L)$ is a matrix corresponding to $\varphi$.
Observe that for different choice of free bases in $N,N'$
the determinants of corresponding matrices differ by 
the multiplication on an invertible scalar and so have equal
orders of zero at the origin.
It is easy to see that we can choose free bases
$v_1,\ldots,v_r$ in $N$ and $v'_1,\ldots,v'_r$ in $N'$ in such a way
 that $\varphi(v_i)=t^{s_i}v'_i$ for some
$s_i\in\mathbb{Z}_{\geq 0}$. Let $D$ be the matrix of $\varphi$
with respect to these bases.
The order of zero of  $\det D$ at the origin is
$\sum_{i=1}^r s_i$ and $\dim F^j(N)/(F^j(N)\cap xN)=|\{i: s_i\geq j\}|$.
Since $\sum_i s_i=\sum_{j=1}^{\infty}|\{i: s_i\geq j\}|$, the claim results.

\section{ The anticentre}
\label{sectantie}
For a finite dimensional Lie superalgebra $\fp$
the anticentre $\Ap$ can be defined as the set of invariants
of $\Up$ with respect to a twisted adjoint action:
$\Ap:=\Up^{\ad'\fp}$ where $\ad'$ is given by the formula
$$(\ad' g)u=gu-(-1)^{p(g)(p(u)+1)} ug.$$
We see that the odd elements of the anticentre $\Ap$
commute with all elements of $\Up$ and the even elements of $\Ap$
commute with the even elements of $\Up$ and anticommute
with the odd ones. Clearly, the product of two anticentral elements
is central.

In this section we describe $\Ag$, see~\Thm{TUg}. This provides us a bunch
of central elements (see~\Cor{corZUg}). As it was indicated
in~\ref{linfactint}, the central elements are useful
for the proof of linear factorizability
of Shapovalov determinants, see
Sect.~\ref{shapf} below.

\subsection{Schur's lemma}
Recall that Schur's lemma for Lie superalgebras takes
the following form: for a simple $\fp$-module $V=V_{\ol{0}}\oplus V_{\ol{1}}$ 
either $\End(V)^{\ad\fp}=\mathbb{C}\id$ or
$\End(V)^{\ad\fp}=\mathbb{C}\id\oplus \mathbb{C}\upsilon$ 
where $\upsilon$ is odd and satisfies $\upsilon^2=\id$.

Let $\theta:V\to V$ be the map $v\mapsto (-1)^{p(v)}v$. It is easy to see that
the action of $z$ on a simple $\fp$-module $V$ is proportional to
\begin{equation}
\label{tZp}
\begin{array}{ll}
\id, & \text { if } z\in\Zp \text { and } z \text { is even, }\\
0, & \text { if } z\in\Zp \text { and } z \text { is odd, }\\
\theta, & \text { if } z\in\Ap \text { and } z \text { is even, }\\
\upsilon\theta, & \text { if } z\in\Ap \text { and } z \text { is odd.}
\end{array}\end{equation}

\subsubsection{}
For $\fp$ being a basic classical or $Q$-type Lie superalgebra 
the formula~(\ref{tZp}) holds for Weyl modules as well.
This follows from the fact that the action of $z$ on a Weyl module
is determined by its action  on the highest weight space 
which is a simple $\fh$-module.

\subsection{}
\label{antie}
Assume that $\fp$ satisfies the following condition
$$\bigwedge\! ^{\top}\fp_{\ol{1}} \text{ is a trivial $\fp_{\ol{0}}$-module.}\eqno{(*)}$$
Then the anticentre $\Ap$ admits the following description
(see~\cite{ghost}, Sect.3).

Let $\mathbb{C}v$ be the even trivial $\fp_{\ol{0}}$-module. The induced module
$\Ind_{\fp_{\ol{0}}}^{\fp}\mathbb{C}v$  contains a unique trivial 
$\fp$-submodule; let $u_0\in\Up$ be such that 
$u_0v$  generates this submodule. The map 
$$\vartheta:z\mapsto (\ad' u_0)(z)$$ 
is a linear isomorphism from $\cZ(\fp_{\ol{0}})$ to $\Ap$.
 The map $\vartheta$ is defined up to
a multiplicative scalar: if $u_0,u'_0$ are such that $\fp(u_0v)=\fp(u'_0v)=0$
and $\vartheta,\vartheta'$ are the corresponding isomorphisms then
$\vartheta=c\vartheta'$ for some $c\in\mathbb{C}^*$.
One has $\gr(\vartheta(z))\in\bigwedge^{\top}\fp_{\ol{1}}\gr z$.
In particular, the anticentre is pure even if $\dim\fp_{\ol{1}}$ is even
 and is pure odd otherwise.

\subsubsection{}
Define a filtration $\cF_{1/2}$ of $\Up$ by letting the odd elements of $\fp$
have degree $1/2$ and the even elements of $\fp$ have degree $1$.
Denote by $\deg_{1/2}u$ the degree of $u\in\Up$ with respect to $\cF_{1/2}$.
Observe that  $\deg_{1/2}u=\deg u$ for $u\in\cU(\fp_{\ol{0}})$.

We claim that
\begin{equation}
\label{deg12phi}
\deg_{1/2}\vartheta(z)\leq\frac{\dim\fp_{\ol{1}}}{2}+\deg z.
\end{equation}
Indeed, let $X$ be a basis of $\fp_{\ol{1}}$.
The module $\Ind_{\fp_{\ol{0}}}^{\fp}\mathbb{C}v$ has a basis of the form
$\{u_kv\}$ ($k=1,\ldots, 2^{\dim\fp_{\ol{1}}}$) where each $u_k$ is a product
of distinct elements of $X$ and so
$\deg_{1/2}u_k\leq\frac{\dim\fp_{\ol{1}}}{2}$.
One can choose $u_0$ to be a linear combination of $u_k$. Then
$\deg_{1/2}u_0\leq\frac{\dim\fp_{\ol{1}}}{2}$ and hence~(\ref{deg12phi}).

\subsubsection{Element $T_{\fp}$.}
\label{Tp}
By~\ref{antie}, $\Ap$ contains a unique (up
to a scalar) element $T_{\fp}$ such that
$\gr T_{\fp}\in\bigwedge^{\top}\fp_{\ol{1}}$. One has $T_{\fp}=(\ad' u_0)(1)$ and so
$\deg_{1/2} T_{\fp}\leq \frac{\dim\fp_{\ol{1}}}{2}$.

\subsection{}
Let $\fg$ be a $Q$-type Lie superalgebra. The algebras
$\fp=\fg,\fh$ satisfy the condition (*). In~\ref{ZAh} we show that
$\cA(\fh)=\Sh T_{\fh}$  where $T_{\fh}$ is given by the formula~(\ref{Tfh}).
We describe $T_{\fg}$ and $\Ag$ in~\Thm{TUg} below.

\subsubsection{}
The following proposition is proven in~\cite{lm} (Cor. D):

\begin{prop}{ann}
For any Zariski dense subset $\Omega$ of $\fh_{\ol{0}}^*$ one has
$$\bigcap_{\lambda\in \Omega}\Ann_{\Ug}N(\lambda)=0.$$
\end{prop}

The proof in~\cite{lm} is based on the similar assertion for
semisimple Lie algebras. The assertion for semisimple Lie algebras
can be deduced (see, for instance,~\cite{jbook}, 7.1.9)
from the fact that the determinants
of all Shapovalov forms are not equal to zero.
The same reasoning as in~\cite{jbook}, 7.1.9 works 
in our case:~\Prop{ann} can be easily deduced from
the inequalities $\det B_{\nu}\not=0$ (for all 
$\nu\in Q^+$) obtained in~\ref{corestimsh}.

\subsubsection{}
\begin{lem}{HChAg}
\begin{enumerate}
\item
The restriction of $\HC$ to $\Zg$ is injective and its image
lies in $\Sh^W$. In particular, $\Zg$ is pure even.
\item
The restriction of $\HC$ to $\Ag$ is injective and its image
lies in $\cA({\fh})$.
\end{enumerate}
\end{lem}
\begin{proof} Obviously,  $\HC(\Zg)\subset\cZ(\Uh), \HC(\Ag)\subset\cA(\fh)$.
It is easy to see that $\cZ(\Uh)=\Sh$.
Observe that $z\in \Zg\cup\Ag$ kills a Weyl module $N(\lambda)$
iff $z$ kills its highest weight space $N(\lambda)_{\lambda}$.
For  $v\in N(\lambda)_{\lambda}$ one has
$zv=\HC(z)v$. Combining the above observation with~\Prop{ann}
we get $\HC(z)\not=0$ for $z\not=0$. Hence (ii). 
The inclusion
$\HC(\Zg)\subset \Sh^W$ follows from~\Cor{corshq2qn} (i) and
the fact that $zv=\HC(z)(\lambda)v$ for $v$ being a primitive 
vector of weight $\lambda$.
\end{proof}

\subsection{}
Recall that $T_{\fg}$ is defined up to an invertible scalar.

\begin{thm}{TUg}
\begin{enumerate}
\item We can choose $T_{\fg}$ such that
$$\HC(T_{\fg})=
T_{\fh}\prod_{\alpha\in\Delta_{\ol{0}}^+} h_{\ol{\alpha}}.$$

\item
The restriction of $\HC$ to $\Ag$ is a linear isomorphism
$\Ag\to \Sh^W\HC(T_{\fg})$.
\end{enumerate}
\end{thm}
\begin{proof}
Let $\vartheta:\cZ(\fg_{\ol{0}})\iso\Ag$ be the map described in~\ref{antie}.
Recall that $\HC(\Ag)\subset \Ah=\Sh T_{\fh}$. 
This gives an injective linear map $\cZ(\fg_{\ol{0}})\to\Sh$ 
$$z\mapsto u_z\ \ \text{ such that }
\HC(\vartheta(z))=T_{\fh}u_z.$$ 

Retain notation of~\ref{hi}. In~\ref{Thq} we obtain the formula
\begin{equation}
\label{tfh}
t_{\fh}:=T^2_{\fh}=\left\{\begin{array}{ll}
\pm h_1\ldots h_n & \text{ for } \fg=\fq(n),\fpq(n)\\
\pm\sum h_1\ldots \hat{h_i}\ldots h_n
& \text{ for } \fg=\fsq(n),\fpsq(n).
\end{array}\right.
\end{equation}

Let us estimate $\deg u_z$.
Observe that $\deg_{1/2} u_z=\deg u_z$ because $u_z\in\Sh$.
Taking into account that $\gr_{1/2}\Uh\cong\Uh$ and that
$u_z\in\Sh$ is a non-zero divisor in $\Uh$, we conclude that
$\gr_{1/2} u_z$ is a non-zero divisor in $\gr_{1/2}\Uh$ and so
$\deg_{1/2} (u_zT_{\fh})=\deg_{1/2} u_z+\deg_{1/2}T_{\fh}$.
By~(\ref{deg12phi}), $\deg_{1/2}T_{\fh}\leq\frac{\dim\fh_{\ol{1}}}{2}$.
Since $t_{\fh}=T_{\fh}^2\in \Sh$ and $\deg t_{\fh}=\dim\fh_{\ol{1}}$  one has
$\deg_{1/2}t_{\fh}=\dim\fh_{\ol{1}}$. Therefore
$\deg_{1/2}T_{\fh}=\frac{\dim\fh_{\ol{1}}}{2}$.
Using~(\ref{deg12phi}) we obtain the following estimation
\begin{equation}
\label{usual}
\deg u_z\leq \frac{\dim\fg_{\ol{1}}}{2}+\deg z-\deg_{1/2}T_{\fh}=
\dim\fn^+_{\ol{1}}+\deg z.
\end{equation}
 
We claim that $u_z$ is divisible by $h_{\ol{\alpha}}$ for any root
$\alpha\in\Delta_{\ol{0}}^+$. Since
$t_{\fh}$ is not divisible by $h_{\ol{\alpha}}$, it suffices to show
that $p:=\HC(\vartheta(z))^2=t_{\fh}u_z^2$ is 
divisible by $h_{\ol{\alpha}}$ for any root
$\alpha\in\Delta_{\ol{0}}^+$. Since $\vartheta(z)^2\in\Zg$ \Cor{HChAg} gives
$p\in\Sh^W$. As a consequence, it is enough to verify
that $p$ is 
divisible by $h_{\ol{\alpha}}$ for any $\alpha\in\pi$.
Fix $\alpha\in\pi$. In the notation of~\ref{hi} one has
$h_{\ol{\alpha}}=h_i+h_{i+1}$ for some $i$.
Take $\lambda\in\fh_{\ol{0}}^*$ such that $h_i(\lambda)=h_{i+1}(\lambda)=0$.
Observe that $t_{\fh}(\lambda)=0$ and so
$p(\lambda)=0$. 
For any $c\in\mathbb{C}$ one has $h_{\ol{\alpha}}(\lambda-c\alpha)=0$
and therefore, by~\Cor{corshq2qn} (ii),  
$N(\lambda-c\alpha)$ has a primitive vector
of the weight $\lambda-(c+1)\alpha$. This gives 
$p(\lambda-c\alpha)=p(\lambda-(c+1)\alpha)$ because $\vartheta(z)^2$ is central.
Now $p(\lambda)=0$ forces
$p(\lambda-c\alpha)=0$ for all $c\in\mathbb{C}$. Any
 $\lambda'\in\fh_{\ol{0}}^*$ satisfying
$h_{\ol{\alpha}}(\lambda')=0$ takes form
$\lambda'=\lambda''-c\alpha$ where $c=-\lambda'(h_i)$ and
$\lambda''$ is such that $h_i(\lambda)=h_{i+1}(\lambda)=0$. 
Thus $p(\lambda')=0$ for all $\lambda'$ 
satisfying $h_{\ol{\alpha}}(\lambda')=0$.
This means that $p$ is divisible by $h_{\ol{\alpha}}$ and the claim follows.

Hence
$u_z=\prod_{\alpha\in\Delta_{\ol{0}}^+} h_{\ol{\alpha}}u'_z$.
By~(\ref{usual}), $\deg u'_z\leq \deg z$.
In particular, for $z\in\mathbb{C}$ one has $\deg u'_z=0$ that is
$u'_z\in\mathbb{C}$ and thus
$u_1=\prod_{\alpha\in\Delta_{\ol{0}}^+} h_{\ol{\alpha}}$ up to a scalar.
Since $T_{\fg}=\vartheta(1)$ up to an invertible scalar,
this implies (i).

Now we have $\HC(\vartheta(z))=\HC(T_{\fg})u'_z$.
One can easily deduce from~\Cor{corshq2qn} (i) that
$u'_z\in\Sh^W$.
Define the map $\vartheta':\cZ(\fg_{\ol{0}})\to\Sh^W$ by $z\mapsto u'_z$.
Obviously, $\vartheta'$ is a linear injective map. 
Since $\dim\{z\in \cZ(\fg_{\ol{0}})|\ \deg z=m\}=\dim\{s\in\Sh^W|\ \deg s=m\}$
and $\deg u'_z\leq \deg z$, the map $\psi$ is surjective.  This proves (ii).
\end{proof}

\subsubsection{}\label{tg}
We obtain
\begin{equation}
\label{tfg}
t_{\fg}:=\HC(T_{\fg}^2)=t_{\fh}
\bigl(\prod_{\alpha\in\Delta_{\ol{0}}^+} h_{\ol{\alpha}}\bigr)^2.
\end{equation}
where $t_{\fh}=T_{\fh}^2$ is given by the formula~(\ref{tfh}).

Since the product of two anticentral elements is central,
\Thm{TUg} implies the following corollary.

\subsubsection{}
\begin{cor}{corZUg}
The restriction of $\HC$ induces an algebra monomorphism
$\Zg\to\Sh^W$ 
whose image contains $t_{\fg}\Sh^W$.
\end{cor}

\section{Computation of Shapovalov determinants}
\label{shapf}
In this section we calculate Shapovalov determinants, see~\Thm{thmdetsh}.
Then, in~\ref{coval},
we study Weyl modules $N(\lambda)$ for $\lambda$ having the smallest
possible degeneracy. Results of~\ref{coval} are used for the calculation
of $\Zg$ in Sect.~\ref{sectcentre}.

\subsection{}
\label{detsh}
\begin{thm}{thmdetsh}
Up to a non-zero scalar,
$$\det B_{\nu}=
\bigl(\prod_{\alpha\in\Delta_{\ol{0}}^+} (h_{\ol{\alpha}})^{
\tau_{{\alpha}}(\nu-\alpha)}
\prod_{\alpha\in\Delta_{\ol{0}}^+}
\prod_{r\geq 1}(h_{\alpha}-r)^{\tau (\nu-r\alpha)}\bigr)^{2^{\dim\fh_{\ol{1}}}}.$$
\end{thm}

\begin{cor}{cordetsh}
$$\Norm S_{\nu}=\prod_{\alpha\in\Delta_{\ol{0}}^+} (h_{\ol{\alpha}})^{
\tau_{{\alpha}}(\nu-\alpha)}
\prod_{\alpha\in\Delta_{\ol{0}}^+}
\prod_{r\geq 1}(h_{\alpha}-r)^{\tau (\nu-r\alpha)}.$$
\end{cor}

This theorem is proven in~\ref{zeros},\ref{setR} below.
First in~\ref{zeros} we deduce from~\Cor{corZUg} and~\Claim{corestimsh}
that $\det B_{\nu}$  admits a linear factorization. Then using the sum 
formula~(\ref{mlambda}) we compute the multiplicities
of linear factors of $\det B_{\nu}$
in~\ref{setR}. 

\subsubsection{}\label{gamm}
For $h\in\fh_{\ol{0}},c\in {\mathbb C}$  set
$$\gamma_{h}:=\{\zeta\in\fh^*_{\ol{0}}|\ h(\zeta)=0\},\ \ \
\gamma_{h,c}:=\{\zeta\in\fh^*_{\ol{0}}|\ h(\zeta)=c\}.$$

Set
$$\Gamma:=\{\gamma_{\ol{\alpha}},
\gamma_{\alpha,r}:\alpha\in\Delta_{\ol{0}}^+, r=1,2,\ldots\}.$$

Say that $\lambda\in\fh^*_{\ol{0}}$ is {\em regular} 
if $\lambda\not\in\gamma$ for any $\gamma\in\Gamma$.
Note that $N(\lambda)$ is simple iff $\lambda$ is regular.

Say that $\lambda\in\fh^*_{\ol{0}}$ is {\em subregular} 
if it belongs to exactly one of the hyperplanes from $\Gamma$.

\Thm{thmdetsh} implies that $N(\lambda)$ is simple iff $\lambda$ is regular.

In~\ref{coval} we will show that the Jantzen filtration of $M(\lambda)$
has length  two if $\lambda$ is subregular. This fact is used
in the computation of the centre (see Sect.~\ref{sectcentre}).

\subsection{Proof of~\Thm{thmdetsh}: set of zeros}
\label{zeros}
In this subsection we show that $\det B_{\nu}$ admits 
a linear factorization of the form~(\ref{det22}).

\subsubsection{}
\label{shzero1}
Fix $\nu\in Q^+$. Assume that $\lambda\in\fh_{\ol{0}}^*$ is such
that $t_{\fg}(\lambda)\not=0$ and $\det B_{\nu}(\lambda)=0$. 
By~\Cor{corshap}, $\det B_{\nu}(\lambda)=0$
iff $\ol{N(\lambda)}_{\lambda-\nu}\not=0$.
The last implies the existence of a primitive vector in
$N(\lambda)_{\lambda-\mu}$ for some $0<\mu \leq \nu$.
Then $\HC(z)(\lambda)=\HC(z)(\lambda-\mu)$
for any $z\in\Zg$. By~\Cor{corZUg},
$\HC(\Zg)\supset t_{\fg}\Sh^W$. Now 
$t_{\fg}(\lambda)\not=0$ gives $\lambda-\mu=w\lambda$ for some $w\in W$.
In particular, $(\lambda-\mu,\lambda-\mu)=(\lambda,\lambda)$
where $(-,-)$
stands for the standard $W$-invariant bilinear form on $\fh_{\ol{0}}^*$.
We conclude that 
$(\lambda,\mu)-{\frac{1}{2}}(\mu,\mu)=0$ for some $0<\mu \leq \nu$.

Using the formula 
$t_{\fg}=t_{\fh}\prod_{\alpha\in\Delta^+}h_{\ol{\alpha}}^2$ we get
$$\det B_{\nu}(\lambda)=0 \ \Longrightarrow t_{\fh}(\lambda)=0\ \text{ or }
\lambda\in 
\bigcup_{\gamma\in\tilde{\Gamma}}\gamma$$
where $\tilde{\Gamma}$ is the following set of hyperplanes
$$\tilde{\Gamma}:=
\bigcup_{\mu\in Q^+,\mu\not=0}
\{\zeta\in\fh_{\ol{0}}^*|\ (\zeta,\mu)-{\frac{1}{2}}(\mu,\mu)=0\}
\cup \bigcup_{\alpha\in\Delta^+} \{\zeta\in\fh_{\ol{0}}^*|\ 
h_{\ol{\alpha}}(\zeta)=0\}.$$

As a consequence, the hypersurface $\det B_{\nu}=0$ 
is the union of some hyperplanes from $\tilde{\Gamma}$
and some irreducible components of the hypersurface $t_{\fh}=0$.
Consider the case $\fg=\fsq(n),\fpsq(n)$ where $n\geq 3$.
Then
$t_{\fh}=\pm\sum h_1\ldots \hat{h_i}\ldots h_n$
and so the hypersurface $t_{\fh}=0$ is irreducible.
From the formula for leading
term of $\det B_{\nu}$ given in~\ref{corestimsh} 
we see that the hypersurface $t_{\fh}=0$ is not a component
 of the hypersurface $\det B_{\nu}=0$.
For $\fg=\fsq(2)$ the  hypersurface $t_{\fh}=0$
coincides with the hyperplane $h_{\ol{\alpha}}=0$.
Now consider the case $\fg=\fq(n),\fpq(n)$. Then
$t_{\fh}=\pm h_1\ldots \ldots h_n$
and from~\ref{corestimsh} we see that the hyperplanes $h_i=0$
are not components of the hypersurface $\det B_{\nu}=0$.
Finally we get
$$\det B_{\nu}(\lambda)=0 \ \Longrightarrow
\lambda\in \displaystyle\bigcup_{\gamma\in\tilde{\Gamma}}\gamma.$$

\subsubsection{}
\label{shzero3}
Let us analyze the set of zeros  of $\det B_{\nu}$ more carefully.

Let $\mu\in\fh_{\ol{0}}^*$ be such that the hyperplane 
$\gamma:=\{\zeta\in\fh_{\ol{0}}^*|\ (\zeta,\mu)-{\frac{1}{2}}(\mu,\mu)=0\}$
is a component of the hypersurface $\det B_{\nu}=0$ and is not a component of
the hypersurface $t_{\fg}=0$.
By~\ref{shzero1} for any $\zeta\in\gamma$ satisfying $t_{\fg}(\zeta)\not=0$
there exists $w\in W$
such that $\zeta-\mu=w\zeta$; in other words,  
$\gamma\subset\cup_{w\in W} X_w\cup\{\zeta|\ t_{\fg}(\zeta)=0\}$
where  $X_w:=\{\eta\in\fh_{\ol{0}}^*|\ 
\eta-\mu=w\eta\}$. Since each $X_w$ is a proper linear subspace
of $\fh^*_{\ol{0}}$ we get
$\gamma=X_w$ for some $w\in W$. Writing $\gamma=\mu/2+\gamma'$
where $\gamma':=\{\zeta\in\fh_{\ol{0}}^*|\ (\zeta,\mu)=0\}$
we obtain $w(\mu/2+\zeta)=-\mu/2+\zeta$ for any
$\zeta\in\gamma'$. In particular, 
$w\mu=-\mu$ and so $w\zeta=\zeta$ for all $\zeta\in\gamma'$.
Hence $w$ is the reflection with respect to $\mu$. 
Since the only reflections in 
the symmetric group $S_n$ are transpositions  $(ij)$ 
which correspond to the reflections with respect to the roots,
we conclude that $\mu=r\alpha$
for some $r\in\mathbb{C},\alpha\in\Delta^+$. Taking into
account that $0<\mu\leq\nu$, one obtains $r\in\mathbb{Z}_+$. Consequently,
$$\{\zeta\in\fh_{\ol{0}}^*|\ (\zeta,\mu)-{\frac{1}{2}}(\mu,\mu)=0\}=
\{\zeta\in\fh_{\ol{0}}^*|\ (\zeta,\alpha)-r=0\}.$$

\subsubsection{}\label{det2}
Summarizing~\ref{shzero1} and~\ref{shzero3} we conclude that
up to a non-zero scalar
\begin{equation}
\label{det22}
\det B_{\nu}=
\displaystyle{\prod_{\alpha\in\Delta_{\ol{0}}^+}} h_{\ol{\alpha}}^{d'_{\alpha}(\nu)}
\displaystyle{\prod_{\alpha\in\Delta_{\ol{0}}^+,r\geq 1}}
(h_{\alpha}-r)^{d_{r\alpha}(\nu)}.
\end{equation}
Comparing the above formula with the formula for leading term of 
$\det B_{\nu}$ (see~\ref{corestimsh}) we obtain
$$d'_{\alpha}(\nu)=2^{\dim\fh_{\ol{1}}}\tau_{\alpha}(\nu-\alpha),\ \ \ 
\sum_{r\geq 1} d_{r\alpha}(\nu)=2^{\dim\fh_{\ol{1}}}
\sum_{m\geq 1}\tau(\nu-m\alpha).$$
for all $\alpha\in\Delta_{\ol{0}}^+$.

\subsection{Proof of~\Thm{thmdetsh}: computation of multiplicities.}
\label{setR}
In this subsection we compute the multiplicities $d_{r\alpha}(\nu)$.
Fix $\alpha\in\Delta_{\ol{0}}^+, r\in \mathbb{Z}_{>0}$ and set
$\gamma:=\gamma_{\alpha,r}$. Let $\breve{\gamma}$ be the set of subregular
points $\lambda$ in $\gamma$ satisfying $t_{\fh}(\lambda)\not=0$.
Observe that  $\breve{\gamma}$ is dense in  $\gamma$.

Let $s_{\alpha}\in W$ be the reflection with respect to $\alpha$.
For $\lambda\in\gamma$ one has $s_{\alpha}\lambda=\lambda-r\alpha$.

\subsubsection{}
\begin{lem}{genpt}
If $\lambda\in \breve{\gamma}$ then
\begin{enumerate}
\item
$\Hom_{\fg}(N(\mu), N(\lambda))=0$ for all 
$\mu\not=\lambda, s_{\alpha}\lambda=\lambda-r\alpha$;
\item
$N(s_{\alpha}\lambda)$ is simple and $t_{\fh}(s_{\alpha}\lambda)\not=0$.
\end{enumerate}
\end{lem}
\begin{proof}
Suppose that $\Hom_{\fg}(N(\mu), N(\lambda))\not=0$ and $\mu\not=\lambda$. 
Then $\lambda-\mu\in Q^+$ and
moreover $\mu\in W\lambda$ by~\Cor{corZUg} since $t_{\fg}(\lambda)\not=0$. 
Take $w\in W$ such that $\mu=w\lambda$. There exist linearly independent
positive roots $\beta_1,\ldots,\beta_k$ such that 
$w=s_{\beta_1}\ldots s_{\beta_k}$ (see, for instance,~\cite{jbook}, A.1.18).
One has $\lambda-w\lambda=(\lambda,\beta_k)\beta_k+
(\lambda,s_{\beta_k}\beta_{k-1})\beta_{k-1}+\ldots\in Q^+$.
Recall that $(\lambda,\beta)\in {\mathbb Z}$ for $\beta\in\Delta^+$ forces
$\beta=\alpha$. Therefore $k=1$ and $\beta_1=\alpha$. 
Thus $\mu=s_{\alpha}\lambda$ and (i) follows.

Since $t_{\fg}$ is $W$-invariant one has $t_{\fg}(s_{\alpha}\lambda)\not=0$.
Now to verify the simplicity of $N(s_{\alpha}\lambda)$, we need to check
that for any $\beta\in\Delta^+$ the value
$h_{\beta}(s_{\alpha}\lambda)=(s_{\alpha}\beta,\lambda)$
is not a positive integer. If $s_{\alpha}\beta\in\Delta^+$ then 
$(s_{\alpha}\beta,\lambda)\not\in \mathbb{Z}_{>0}$ 
since $\lambda$ is subregular. One has
 $(s_{\alpha}\alpha,\lambda)=-r\not\in \mathbb{Z}_{>0}$. Finally, 
if $s_{\alpha}\beta\not\in\Delta^+$ and $\alpha\not=\beta$
then $s_{\alpha}\beta=\beta-m\alpha$ for some $m>0$ and
so $(s_{\alpha}\beta,\lambda)=(\beta,\lambda)-mr\not\in \mathbb{Z}_{>0}$
because $(\beta,\lambda)\not\in \mathbb{Z}_{>0}$.
\end{proof}

\subsubsection{}
\label{chum}
Take $\lambda\in \breve{\gamma}$ and consider the Jantzen filtration on 
$M(\lambda)$. Since $\gamma=\gamma_{\alpha,r}$ is the only element of $\Gamma$
 which contains $\lambda$,
the sum formula~(\ref{mlambda}) gives
\begin{equation}
\label{m3}
d_{r\alpha}(\nu)=\sum_{j\geq 1} \dim F^j(M(\lambda))_{\lambda-\nu}
\end{equation}
for any $\nu\in Q^+$.
\Lem{genpt} implies that any proper submodule $N$ of $N(\lambda)$
satisfies $\ch N=i\ch M(s_{\alpha}\lambda)$ for some $i\geq 0$. 
Therefore for each $j\geq 1$ we have
$$\ch F^j(M(\lambda))=k_j\ch N(\lambda-r\alpha)$$
for some $k_j\geq 0$. Putting $m_r:=\sum_{j\geq 1} k_j$
one obtains
$$d_{r\alpha}(\nu)=m_r\tau (\nu-r\alpha)$$
for all $\nu\in Q^+$.

\subsubsection{}
By~\ref{det2} we have $\sum_{r\geq 1} d_{r\alpha}(\nu)=2^{\dim\fh_{\ol{1}}}
\sum_{j\geq 1}\tau(\nu-j\alpha)$
that is
$$\sum_{r\geq 1}m_r\tau (\nu-r\alpha)=
2^{\dim\fh_{\ol{1}}}\sum_{j\geq 1}\tau(\nu-j\alpha).$$
Let us show that $m_r=2^{\dim\fh_{\ol{1}}}$ by induction on $r$. Indeed,
substituting $\nu:=\alpha$ we get $m_1=2^{\dim\fh_{\ol{1}}}$. Now assuming that 
$m_1=\ldots=m_{i-1}=2^{\dim\fh_{\ol{1}}}$ and putting $\nu=i\alpha$ one concludes
$m_i=2^{\dim\fh_{\ol{1}}}$ as required. Hence
$$d_{r\alpha}(\nu)=2^{\dim\fh_{\ol{1}}}\tau (\nu-r\alpha)$$
for all $\alpha\in\Delta_{\ol{0}}^+, r\in \mathbb{Z}_{>0}$, and $\nu\in Q^+$.

This completes the proof of~\Thm{thmdetsh}.

\subsection{The Jantzen filtration in subregular points.}
\label{coval}
In this subsection we show that the Jantzen filtration
in the points having the smallest possible degeneracy
has length two.

Retain notation of~\ref{gamm}. Set
$$V^{\oplus}(\mu):=M(\mu)/\ol{M(\mu)}.$$

Recall that, by~\ref{clla}, Verma modules $M(\lambda)$ are the direct sum of
Weyl modules if $\lambda$ are such that $t_{\fh}(\lambda)\not=0$.

\subsubsection{}
\label{remolM}
For $\gamma=\gamma_{\alpha,r}$,
let $\breve{\gamma}$ be the set of subregular points of $\gamma$
satisfying $t_{\fh}(\lambda)\not=0$ (as in~\ref{setR}).
Clearly, $\breve{\gamma}$ is dense in $\gamma$.

Recall that
 $\fg\not=\fpsq(2)$ throughout the paper.
For $\gamma=\gamma_{\ol{\alpha}}$
and $\fg\not=\fsq(2)$, let $\breve{\gamma}$
be the maximal set of subregular points of $\gamma$
which does not meet the hypersurface
$t_{\fh}=0$ and is invariant under the shift by $\alpha$:
$\lambda\in\breve{\gamma}\ \Longrightarrow\ 
\lambda-\alpha\in\breve{\gamma}$.
The condition $\fg\not=\fsq(2),\fpsq(2)$ ensures that
the hypersurface $t_{\fh}=0$ does not contain hyperplanes parallel
to $\gamma$. As a consequence, $\breve{\gamma}_{\ol{\alpha}}$
is dense in $\gamma_{\ol{\alpha}}$.

For $\fg=\fsq(2)$ let $\breve{\gamma}_{\ol{\alpha}}$ be the set of subregular
points of $\gamma_{\ol{\alpha}}$.

\subsubsection{}
\begin{prop}{length}
For any $\lambda\in\breve{\gamma}$ one has
\begin{enumerate}
\item
A Jordan-H\"older filtration of module $N(\lambda)$ has length two.
\item 
The Jantzen filtration of $M(\lambda)$ has length two.
\item $\ol{M}(\lambda)=V^{\oplus}(\tilde{\lambda})$ where
$\tilde{\lambda}=\lambda-r\alpha$ if $\gamma=\gamma_{\alpha,r}$
and $\tilde{\lambda}=\lambda-\alpha$ if $\gamma=\gamma_{\ol{\alpha}}$.
\end{enumerate}
\end{prop}
\begin{proof}
The case $\fg=\fsq(2)$ was treated in~\ref{jntz}.
Suppose $\fg\not=\fsq(2)$.

From the formula for Shapovalov determinants we know
that $N(\lambda)$ is not simple if $\lambda\in\gamma$.
Since $\lambda$ is subregular, 
$\dim\ol{N(\lambda)}_{\tilde{\lambda}}\not=0$.
Let $E(\mu)$ be a simple $\Cl(\lambda)$-module ($E(\mu)=N(\mu)_{\mu}$).
One has $\dim\ol{N(\lambda)}_{\tilde{\lambda}}=ke$
where $e:=\dim E(\tilde{\lambda})$ and $k$ is a positive integer.
By the construction of $\breve{\gamma}$ one has $t_{\fh}(\lambda),
t_{\fh}(\tilde{\lambda})\not=0$. 
Therefore $\Cl(\lambda),\Cl(\tilde{\lambda})$ are non-degenerate.
Then $M(\lambda)$  is the direct sum of copies
of $N(\lambda)$ and of $\Pi(N(\lambda))$ and
$M(\tilde{\lambda})$ has the similar form
with the same number of summands.
Then $\dim \ol{M(\lambda)}_{\tilde{\lambda}}=2^{\dim\fh_{\ol{1}}}k$.
However $\det B_{r\alpha}$ at $\lambda$ has zero of order $2^{\dim\fh_{\ol{1}}}$
so $\dim\ol{M(\lambda)}_{\lambda-r\alpha}\leq 2^{\dim\fh_{\ol{1}}}$.
Therefore $k=1$ so $\ol{N(\lambda)}_{\tilde{\lambda}}$ is a simple
$\fh$-module and 
$\ol{M(\lambda)}_{\tilde{\lambda}}\cong \Cl(\tilde{\lambda})$
as $\fh$-module.
In particular, $\ol{M(\lambda)}$ has a subquotient isomorphic 
to $V^{\oplus}(\tilde{\lambda})$.
Observe that $\ol{M(\lambda)}\cong V^{\oplus}(\lambda-\alpha)$ iff
$\ol{N(\lambda)}$ is simple. Hence (i) and (iii) become equivalent.

Consider the case $\gamma=\gamma_{\alpha,r}$. 
Combining~\Lem{genpt} 
 with the above conclusions
we obtain  (i), (iii) and moreover that
$\ol{M(\lambda)}\cong M(\lambda-r\alpha)$.
The sum formula gives 
$$\sum_{j=1}^{\infty}\ch  F^j(M(\lambda))=
\sum_{\mu}e^{\lambda-r\alpha-\mu}\tau(\mu)=\ch M(\lambda-r\alpha)=\ch 
 F^1(M(\lambda)).$$
The equality  $F^1(M(\lambda))\cong M(\lambda-r\alpha)$ forces
$F^j(M(\lambda))=0$ for $j>1$. This completes the proof for
$\gamma=\gamma_{\alpha,r}$.

Fix $\gamma:=\gamma_{\ol{\alpha}}$ and
take $\lambda\in \breve{\gamma}$. 
Let us prove (iii).
Observe that for
a subregular point $\lambda'\in\gamma$  the sum formula
gives
\begin{equation}\label{tsss}
\sum_{j=1}^{\infty}\ch  F^j(M(\lambda'))=
\sum_{\mu}e^{\lambda'-\alpha-\mu}\tau_{\alpha}(\mu)
\end{equation}
and therefore
$$\begin{array}{l}\ch \ol{M(\lambda)}\leq 
\sum_{\mu}e^{\lambda-\alpha-\mu}\tau_{\alpha}(\mu),\\
\ch \ol{M(\lambda-\alpha)}\leq 
\sum_{\mu}e^{\lambda-2\alpha-\mu}\tau_{\alpha}(\mu)=
\sum_{\mu}e^{\lambda-\alpha-\mu}\tau_{\alpha}(\mu-\alpha).
\end{array}$$
where $\sum a_{\mu} e^{\mu}\geq\sum b_{\mu} e^{\mu}$ means that
$a_{\mu}\geq b_{\mu}$ for all $\mu$.
Using the formula $\tau(\nu)=\tau_{\alpha}(\nu)+\tau_{\alpha}(\nu-\alpha)$
we get
$$\ch \ol{M(\lambda-\alpha)}+\ch \ol{M(\lambda)}\leq\ch M(\lambda-\alpha).$$

On the other hand,
since $V^{\oplus}(\lambda-\alpha)$ is a subquotient of $\ol{M(\lambda)}$,
we obtain
$$\ch \ol{M(\lambda)}\geq \ch V^{\oplus}(\lambda-\alpha)
=\ch M(\lambda-\alpha)-\ch \ol{M(\lambda-\alpha)}.$$
Comparing the above inequalities, we conclude
\begin{equation}\label{chopl}
\ch \ol{M(\lambda)}=\ch V^{\oplus}(\lambda-\alpha)=
\sum_{\mu}e^{\lambda-\alpha-\mu}\tau_{\alpha}(\mu)
\end{equation}
so $\ol{M(\lambda)}\cong V^{\oplus}(\lambda-\alpha)$ as required.
This proves (iii) and (i).
Finally, combining~(\ref{chopl}) with~(\ref{tsss}) 
we obtain $F^j(M(\lambda))=0$ for $j>1$.
\end{proof}

\subsubsection{}
\begin{cor}{thmcrk}
Assume that $\lambda\in\fh^*_{\ol{0}}$ is subregular.
Then the Jantzen filtration of $M(\lambda)$ has length two:
$F^2(M(\lambda))=0$.
\end{cor}
\begin{proof}
Recall that the Jantzen filtration of $M(\lambda)$ has length two
iff, for any $\nu$ the order of zero of $\det B_{\nu}$ at point $\lambda$
is equal to the corank of $B_{\nu}$ at point $\lambda$.
Fix a hyperplane $\gamma\in\Gamma$ and take a subregular
point $\lambda\in\gamma$.
Observe that for all subregular points of  $\gamma\in\Gamma$ 
the order of zero of $\det B_{\nu}$
is the same number, say, $r(\nu)$. From~\Prop{length} (ii), the corank of 
$B_{\nu}$ at point $\lambda$ is at least $r(\nu)$. Since the corank
does not exceed the order of zero, it is equal to $r(\nu)$.
The assertion follows.
\end{proof}

\subsubsection{}
Take $\gamma\in\Gamma$ and set $\nu:=r\alpha$ if $\gamma=\gamma_{\alpha,r}$
and $\nu:=\alpha$ if $\gamma=\gamma_{\ol{\alpha}}$.
By~\Thm{thmcrk}, $M(\lambda)$ has $2^{\dim\fh_{\ol{1}}}$
linearly independent primitive vectors of the weight $\lambda-\nu$
if $\lambda\in\gamma$ is a subregular point.
It is easy to see that  $M(\lambda)$ has
$k$ linearly independent primitive vectors of a given
weight iff a certain matrix with entries in $\Sh$ at the point $\lambda$
has corank equal to $k$.
Since the set of subregular points is dense in $\gamma$, 
$M(\lambda')$ has at least $2^{\dim\fh_{\ol{1}}}$
linearly independent primitive vectors of the weight $\lambda'-\nu$
for any $\lambda'\in\gamma$.

\subsubsection{}
\begin{cor}{primv}
The module $N(\lambda)$ is simple iff $h_{\ol{\alpha}}(\lambda)\not=0$
and $h_{\alpha}(\lambda)\not\in\mathbb{Z}_{>0}$ 
for all $\alpha\in\Delta_{\ol{0}}^+$.

If $h_{\ol{\alpha}}(\lambda)=0$ then $\ol{N(\lambda)}$ has a subquotient
isomorphic to $V(\lambda-\alpha)$ or to $\Pi(V(\lambda-\alpha))$.

If $r:=h_{\alpha}(\lambda)\in\mathbb{Z}_{>0}$ 
then $\ol{N(\lambda)}$ has a subquotient
isomorphic to $V(s_{\alpha}\lambda)$ or to $\Pi(V(s_{\alpha}\lambda))$.
\end{cor}

\section{On the centre of completion $\hat{U}$ of $\Ug$}
\label{secthat}
A classical theorem by Chevalley states  that
for a semisimple Lie algebra the restriction of a non-zero $\fg$-invariant
regular function on $\fg$ to $\fh$ is non-zero, that is
for $f\in\cS(\fg^*)^{\fg}$, $f\not=0$ forces $f_{\fh}\not=0$.
This theorem was generalized by A. Sergeev (see~\cite{s3}) 
to all finite-dimensional Lie superalgebras.
If $\fg$ admits an even non-degenerate invariant bilinear form, then
$\fg\cong\fg^*$ as $\fg$-modules. Since $\gr z\in\cS(\fg)^{\fg}$
for $z\in\Zg$ the Chevalley theorem can be reformulated as follows:
$\deg\HC(z)=\deg(z)$  for any $z\in\Zg$.
In this section we will prove a similar statement
for a certain completion $\hat{U}$ of  $\Ug$ where
$\fg$ is a Kac-Moody superalgebra
with a symmetrizable Cartan matrix or $Q$-type Lie superalgebra.
For a  finite-dimensional case this implies $\cZ(\hat{U})=\Zg$.

The centre of  $\hat{U}$ for $\fg$ being a Kac-Moody superalgebra
with a symmetrizable Cartan matrix
was described by Kac in~\cite{klapl}.
We elucidate his approach to $Q$-type algebras in
Section~\ref{sectcentre}; this will give us a description of
$\Zg$.

We denote by $\fg$ a $Q$-type Lie superalgebra. However,
all results of this section are valid for a Kac-Moody superalgebra
with a symmetrizable Cartan matrix (see~\Rem{KMo} for details).

\subsection{A construction $\hat{U}$}
Define on $\Ug$ a topology where a basis of neighborhoods
of zero are left ideals $J(\nu):=\Ug\cU(\fn^+)_{\geq\nu}$
for $\nu\in Q(\pi)$. 
Clearly, $J(\nu)=\Ug$ for $\nu\leq 0$ and $J(\nu)\subset J(\nu')$
if $\nu>\nu'$.
Let $\hat{U}$ be the completion of $\Ug$ with respect
to this topology. It is easy to see 
that the structure of associative algebra and the adjoint action of $\fg$
 can be uniquely extended to $\hat{U}$.
Clearly, $\cZ(\hat{U})=\hat{U}^{\fg}$.

Let $N$ be a $\fg$-module with locally nilpotent action of $\fn^+$ {\footnote 
{this means that $\dim \cU(\fn^+)v <\infty$ for any $v\in N$.}}.
Then the action of $\fg$ can be canonically extended to the action
of $\hat{U}$. In particular, 
the $\fg$--$\fh$ bimodule  $\cM$ introduced in~\ref{cN}
can be viewed as a $\hat{U}$--$\fh$ bimodule.

\subsubsection{}\label{J'}
Set 
$$J'(\nu):=\cU(\fb^-)\cU(\fn^+)_{\nu}.$$
Then, by PBW theorem, 
$J(\nu)=\sum_{\mu\geq \nu} J'(\mu)$ and
\begin{equation}
\label{dprod}
\Ug=\oplus_{\nu\in Q^+}J'(\nu),\ \ \hat{U}=\prod_{\nu\in Q^+}J'(\nu).
\end{equation}
We write the elements of $\hat{U}$ 
in the form $u=\sum_{\nu\in Q^+} u_{\nu}$ where $u_{\nu}\in J'(\nu)$.

\subsubsection{}
Put $\fb^-:=\fn^-+\fh$.
Observe  that $\Uh$ and $\cU(\fb^-)$ are 
closed with respect to our topology.
 One has
$\hat{U}=\Uh\oplus(\fn^-\cU(\fb^-)+\hat{U}\fn^+)$.
Extend the Harish-Chandra projection to $\hat{U}$ along
the above decomposition.

Take $u=\sum_{\nu\in Q^+} u_{\nu}$ where $u_{\nu}\in J'(\nu)$.
Then $\HC(u)=\HC(u_0)$  and $\HC(u)=u_0$ if $u$ has weight zero, because
$J'(0)=\cU(\fb^-)$.

\subsubsection{}
The canonical filtration $\cF$ on $\Ug$ can be naturally extended
to $\hat{U}$, however, the resulted filtration is not exhausting:
$\cup_{r=0}^{\infty}\cF^r(\hat{U})\not=\hat{U}$.
Let $\deg u$ be the degree of $u\in\hat{U}$ with respect to this filtration
($\deg u=\infty$ if $u\not\in \cup_{r=0}^{\infty}\cF^r(\hat{U})$).

\subsection{Results}
The following theorem was suggested to the author by J.~Bernstein.
\subsubsection{}
\begin{thm}{bern}
For any $z\in\cZ(\hat{U})$ one has
$$\deg z=\deg\HC(z).$$ 
In other words,
$z\in\cZ(\hat{U})$ takes form $z=\sum_{\nu\in Q^+} z({\nu})$
where $z(\nu)\in J'(\nu)$ is such that $\deg z(\nu)\leq \deg z(0)$.
\end{thm}

In particular, $\HC$ induces an algebra embedding $\cZ(\hat{U})\to\cU(\fh)$.
The image of the embedding lie in the centre of $\cU(\fh)$ which is
$\cU(\fh_{\ol{0}})$. We describe this image in~\Thm{thmkac} below.

\subsubsection{}
\begin{cor}{corbern}
One has $\cZ(\hat{U})=\Zg$.
\end{cor}
\begin{proof}
From~\Thm{bern}, 
$z(\nu)\in\cF^r(\Ug)$ where $r:=\deg z(0)$. Then
$z(\nu)=0$ if $\nu$ is ``sufficiently large'' that is
$\nu\not\in\Delta^{(r)}$ where
$\Delta^{(r)}:=\{\sum_{i=1}^s\alpha_i|\ \alpha_i\in\Delta^+, s\leq r\}$.
Since $\Delta$ is finite, $\Delta^{(r)}$ is also finite and thus
$z$ is a finite sum of elements of $\Ug$. Hence $z\in\Ug$.
\end{proof}

\subsubsection{}
View $\cM$ as a $\hat{U}$-$\fh$ bimodule.

\begin{prop}{ua}
\begin{enumerate}
\item
$\cM$ is a faithful $\hat{U}$-module.
\item
For $z\in\hat{U}$ one has
$$z\in\cZ(\hat{U})\ \Longleftrightarrow\ zv=v\HC(z)\ \text{ for any }
v\in\cM.$$
\end{enumerate}
\end{prop}

In other words, we have an embedding of $\fg$-modules $\hat{U}\to\End(\cM)$
and the image of a central element is an endomorphism induced
by the right action.

\subsection{Proof of~\Thm{bern}}
Take a non-zero $z\in\hat{U}^{\fg}$  and write
$z=\sum_{\nu\in Q^+} z_{\nu}$
where $z_{\nu}\in J'(\nu)$. Set $r:=\deg z_0$.
Let us prove the inequality
\begin{equation}
\label{z0}
\deg z_{\nu}\leq r
\end{equation}
by induction on $\htt(\nu)$
where the function height on $Q^+$ is given by
$\htt(\sum_{\alpha_i\in\pi} k_i\alpha_i)=\sum k_i$.
The inequality trivially holds for $\htt(\nu)=0$.

\subsubsection{}
Suppose that 
the inequality~(\ref{z0}) holds for all $\nu$ such that $\htt(\nu)\leq k$ and
let us show that~(\ref{z0}) holds for all weights of height $k+1$.
Take $\mu$ of height $k+1$ and write $\mu=\nu+\alpha_0$
where $\alpha_0$ is a simple root and $\nu\in Q^+$ is such 
that $\htt(\nu)=k$.

Denote by $p_{\eta}$ the projection $\hat{U}\to J'(\eta)$
with respect to the decomposition~(\ref{dprod}). 
Since  $z\in\hat{U}^{\fg}$ one has $[e_{\beta},z]=0$ 
for any simple root $\beta$ and any $e_{\beta}\in\fg_{\beta}$. Therefore
$$0=p_{\mu}
\bigl([e_{\beta}, z]\bigr)=p_{\mu}\bigl(
[e_{\beta},z_{\mu-\beta}]\bigr)+p_{\mu}\bigl([e_{\beta},z_{\mu}]\bigr).$$
The canonical filtration is $\ad\fg$-stable
 and $p_{\eta}$-stable; thus, by  the induction hypothesis,
the first summand has degree at most $r$. Hence 
\begin{equation}
\label{Fr-1}
p_{\mu}([e_{\beta},z_{\mu}])\in \cF^r(\Ug).
\end{equation}
for all simple roots $\beta$.
We need to deduce that $z_{\mu}\in \cF^r(\Ug)$.

\subsubsection{}
Fix a basis 
$\{\mathbf{e}^{\mathbf{k}}\}$ in $\cU(\fn^+)_{\mu}$, 
$\{\mathbf{h}^{\mathbf{s}}\}$ in $\cU(\fh)$ and
write $z_{\mu}=\sum_{\mathbf{k}, \mathbf{s}}
{a}_{\mathbf{k},\mathbf{s}}{\mathbf{h}}^{\mathbf{s}}
\mathbf{e}^{\mathbf{k}}$ where $a_{\mathbf{k},\mathbf{s}}$ are elements of 
$\cU(\fn^-)_{-\mu}$. One has
\begin{equation}
\label{Fr-2}
p_{\mu}([e_{\beta},z_{\mu}])=
\sum_{\mathbf{k},\mathbf{s}}[e_{\beta},{a}_{\mathbf{k},\mathbf{s}}]
{\mathbf{h}}^{\mathbf{s}}\mathbf{e}^{\mathbf{k}}\in \cF^r(\Ug).
\end{equation}
by~(\ref{Fr-1}). Recall that $\cS(\fg)=\cS(\fb^-)\otimes\cS(\fn^+)$.
As a consequence, for any  linearly independent elements
$\{s^+_i\}\subset \cS(\fn^+)$ and any elements
$\{s^-_i\}\subset \cS(\fb^-)$ the degree of
$\sum_i s^-_is^+_i$ is the maximum of $\deg s^-_i+\deg s^+_i$.
Denoting the degree of $\mathbf{e}^{\mathbf{k}}$ by  $|\mathbf{k}|$ 
we obtain from~(\ref{Fr-2}) 
\begin{equation}
\label{Fr-3}
\deg \sum_{\mathbf{s}}[e_{\beta},{a}_{\mathbf{k},\mathbf{s}}]
{\mathbf{h}}^{\mathbf{s}}
\leq r-|\mathbf{k}|
\end{equation}
for any $\mathbf{k}$.
Identify $\fb^-$ with $\fg/\fn^+$.
Then $\sum_{\mathbf{s}}
[e_{\beta},{a}_{\mathbf{k},\mathbf{s}}]{\mathbf{h}}^{\mathbf{s}}$ identifies
with $(\ad e_{\beta})\sum_{\mathbf{s}}
{a}_{\mathbf{k},\mathbf{s}}{\mathbf{h}}^{\mathbf{s}}$ where the adjoint
action $\fn^+$ on $\fg/\fn^+$ is induced by the usual adjoint
action. 

Suppose that $\deg z_{\mu}>r$. Then
\begin{equation}
\label{Fr-4}
\deg\sum_{\mathbf{s}}
{a}_{\mathbf{k},\mathbf{s}}{\mathbf{h}}^{\mathbf{s}}> r-|\mathbf{k}|
\end{equation}
for some $\mathbf{k}$.
In the light of~\Lem{coad}, 
the image of $\sum_{\mathbf{s}}
{a}_{\mathbf{k},\mathbf{s}}{\mathbf{h}}^{\mathbf{s}}$ in $\cS(\fg/\fn^+)$
is not $\fn^+$-invariant and thus there exists a simple root $\beta$
such that
$\deg \sum_{\mathbf{s}}
[e_{\beta},{a}_{\mathbf{k},\mathbf{s}}]{\mathbf{h}}^{\mathbf{s}}=\deg
\sum_{\mathbf{s}}
{a}_{\mathbf{k},\mathbf{s}}{\mathbf{h}}^{\mathbf{s}}$.
Comparing~(\ref{Fr-3}) with~(\ref{Fr-4}) we get a contradiction.
The statement follows.
\qed

\subsection{}
Put $\fn:=\fn^+,\fb:=\fb^+$.

\begin{lem}{coad}
View $\fg/\fn$ as $\fn$-module via the adjoint action. Then
$\cS(\fg/\fn)^{\fn}=\cS(\fh)$ where the embedding $\fh\to \fg/\fn$
is induces by the isomorphism $\fb/\fn\iso\fh$.
\end{lem}
\subsubsection{}
\begin{rem}{}
Assume that $\fg$ is finite-dimensional Lie algebra. Then
$\fg/\fn$ identifies with $\fb^*$ via the invariant bilinear form and
$\cS(\fg/\fn)$ identifies with the set of regular functions
on $\fb$,
and $\cS(\fg/\fn)^{\fn}$ identifies with the invariant functions.
The formula $\cS(\fg/\fn)^{\fn}=\cS(\fb/\fn)$
means that for any invariant function $\phi$ one has
$\phi(h+n)=\phi(h)$ for any $h\in\fh,n\in\fn$.
This is a standard fact.
Indeed, let $\phi$ be an invariant function and $N$ be the Lie group
corresponding to $\fn$.
If $h$ is a generic element of $\fh$ then the orbit
$N.h$ is dense in $h+\fn$ and thus $\phi(h+n)=\phi(h)$
for generic $h$. Since the set of generic elements
is dense in $\fh$,  $\phi(h+n)=\phi(h)$ for any $h\in\fh,n\in\fn$.
\end{rem}

{\em Proof of~\Lem{coad}.}
Identify the image of $\fg/\fn$ in $\cS(\fg/\fn)$
with $\fb^-$ and $\cS(\fg/\fn)$ with $\cS(\fb^-)$.
 For an algebra $S$ and its subspaces $X,Y$
denote by $XY$ the span of $xy, x\in X,y\in Y$.

\subsubsection{}
\label{stepi}
First, let us check that for any $\alpha\in\Delta^+$ one has
$\cS(\fg_{-\alpha}+\fh)^{\fn}=\cS(\fh)$.
Let $\fg$ be a $Q$-type.
Taking  a natural basis (see Sect.~\ref{sectq2})
in $\fg_{\alpha}+\fg_{-\alpha}+[\fg_{\alpha},\fg_{-\alpha}]\cong\fsq(2)$ 
we obtain $(\ad e)(u)=h\frac{\partial{u}}{\partial{f}}+
H\frac{\partial{u}}{\partial{F}}$
and $(\ad E)(u)=H\frac{\partial{u}}{\partial{f}}
+h'\frac{\partial{u}}{\partial{F}}$
where $H^2=h'$. Thus $(\ad e)(u)=(\ad E)(u)=0$ forces
$\frac{\partial{u}}{\partial{f}}=\frac{\partial{u}}{\partial{F}}=0$ 
that is $u\in\cS(\fh)$
as required.

\subsubsection{}
Now let us verify the statement of lemma.
Extend the partial order on the root lattice $Q(\pi)$ to
a total order 
compatible with the addition that is
$$\beta>\beta'\ \Longrightarrow\ \beta+\alpha>\beta'+\alpha$$
(for instance, take an embedding $Q(\pi)$ into $\mathbb{R}$).
Let $\Delta^+$ be the set (not the multiset)
of positive roots and let 
$\cP:=\sum_{\alpha\in\Delta^+}\mathbb{Z}_{\geq 0}\alpha$
be the positive lattice generated by $\Delta^+$.
For $\mathbf{k}\in\cP$ set $X_{\mathbf{k}}:=\prod_{\alpha\in\Delta^+}
\fg_{-\alpha}^{k_\alpha}\cS(\fh)\subset \cS(\fb^-)$. Then
$\cS(\fb^-)=\sum_{\mathbf{k}\in\cP} X_{\mathbf{k}}$; let
$p_{\mathbf{k}}$ be the projection with respect to this decomposition.
For $u\in\cS(\fb^-)$ set $\supp u:=\{\mathbf{k}|\ p_{\mathbf{k}}(u)\not=0\}$.

Define a lexicographic order on $\cP$ by putting
$\mathbf{k}>\mathbf{m}$ if for some $\beta\in\Delta^+$
one has $k_{\beta}>m_{\beta}$ and  $k_{\alpha}=m_{\alpha}$
for all $\alpha>\beta$. Suppose that $u\in\cS(\fb^-)^{\fn}$
is such that $u\not\in\cS(\fh)$.
Let $\mathbf{k}$ be the maximal element in $\supp u$ and let
$\alpha\in\Delta^+$ be the minimal root satisfying $k_{\alpha}\not=0$.
Put $\mathbf{k}':=\mathbf{k}-\alpha$ and let us compute 
$p_{\mathbf{k}'}(\ad e)(u)$ for $e\in\fg_{\alpha}$.
One easily sees that 
$$p_{\mathbf{k}'}(\ad e)(u)=
p_{\mathbf{k}'}(\ad e)\bigl(p_{\mathbf{k}}(u)\bigr).$$
Fix a basis $\{f_i\}$ in $\fg_{-\alpha}^{k_{\alpha}}$ and write
$$p_{\mathbf{k}}(u)=\sum f_i a_i,\ a_i\in X_{\mathbf{k}-k_{\alpha}\alpha}.$$
Then for any $e\in\fg_{\alpha}$
$$0=p_{\mathbf{k}'}(\ad e)\bigl(p_{\mathbf{k}}(u)\bigr)=
\sum (\ad e)(f_i) a_i.$$
One has $\cS(\fb^-)=\cS(\fg_{-\alpha})\otimes\cS(\fh)\otimes \cS'$ 
where $\cS':=\otimes_{\beta\not=\alpha} \cS(\fg_{-\beta})$.
Notice that $(\ad e)(f_i)\in\cS(\fg_{-\alpha})\cS(\fh)$
and $a_i\in \cS(\fh)\cS'$.
Thus $\sum (\ad e)(f_i) a_i=0$ forces $\sum (\ad e)(f_i) b_i=0$
for some non-zero $b_i\in\cS(\fh)$ which is impossible by~\ref{stepi}.
\qed

\subsection{Proof of~\Prop{ua}}
For (i) take a non-zero $u\in \Ug_{\lambda}$
and write $u=\sum u_{\nu}$ where $u_{\nu}\in J'(\nu)$.
Let $\mu\in Q^+$ be a minimal element satisfying $u_{\mu}\not=0$;
note that $\mu\geq \lambda$.
Retain notation of~\ref{partition} and fix a PBW basis
$\mathbf{f}^{\mathbf{k}}, \mathbf{k}\in{\cP}(\mu-\lambda)$
in $\cU(\fn^-)_{\lambda-\mu}$. Write
$u_{\mu}=\sum_{\mathbf{k}\in{\cP}(\mu-\lambda)}\mathbf{f}^{\mathbf{k}}
a_{\mathbf{k}}$ where $a_{\mathbf{k}}$ are elements of $\cU(\fb^+)$.
Identify $\cU(\fb^-)$ and $\cM$ as $\fb^-$-modules. 
Thanks to the minimality of $\mu$, for $v\in\cM_{-\mu}$ one has
$uv=u_{\mu}v$. Then
$$u_{\mu}v=\sum_{\mathbf{k}\in{\cP}(\mu-\lambda)}\mathbf{f}^{\mathbf{k}}
a_{\mathbf{k}}v=\sum_{\mathbf{k}}\mathbf{f}^{\mathbf{k}}
\HC(a_{\mathbf{k}}v).$$
Take $\mathbf{k}$ such that $a_{\mathbf{k}}\not=0$.
Since the restriction of Shapovalov map to $\cM_{\mu}$ is non-degenerate,
there exists $v\in \cM_{-\mu}$ satisfying 
$\HC(a_{\mathbf{k}}v)\not=0$. Hence $uv=u_{\mu}v\not=0$ 
and this gives (i).

For (ii) recall that $\cM$ is generated by the image of $1\in\Uh$
which is annihilated by $\fn^+$. This gives the implication 
$$z\in\cZ(\hat{U})\ \Longrightarrow\ zv=v\HC(z)\ \text{ for any }
v\in\cM.$$
The inverse implication follows from (i) and the fact that
$\cM$ is a $\hat{U}$-$\fh$ bimodule and so
the map $v\mapsto v\HC(z)$ lies in $\End_{\hat{U}}(\cM)$.
\qed

\subsection{}
\begin{rem}{KMo}
All constructions and results of this section are valid for
a Kac-Moody superalgebra with a symmetrizable Cartan matrix.
In particular,~\Thm{bern} gives $\deg z=\deg\HC(z)$
for any $z\in\cZ(\hat{U})$ and thus $\HC$ induces an algebra embedding
$\cZ(\hat{U})\to\cU(\fh)$.
The image $\HC\bigl(\cZ(\hat{U})\bigr)$
was described by V.~Kac in~\cite{klapl}, Remark 3 and Section 8.
\Cor{corbern} gives
$$\cZ(\hat{U})=\cZ(\fg)$$
if $\fg$ is a finite-dimensional Kac-Moody (contragredient) superalgebra.

All proofs except for the proof of~\Lem{coad} work for the Kac-Moody case. 
The only difference occurs in the proof of the formula
$\cS(\fg_{-\alpha}+\fh)^{\fn}=\cS(\fh)$,
see~\ref{stepi}.
If $\fg$ is a Kac-Moody superalgebra with a symmetrizable Cartan matrix,
this step can be done as follows. The algebra $\fg$ admits
a non-degenerate invariant bilinear form $(-,-)$
and there exists $h\in\fh$ such that $[e,f]=(f,e)h$
for any $f\in \fg_{-\alpha},e\in\fg_{\alpha}$. 
Let $\{e_i\}\subset\fg_{\alpha}$ and $\{f_i\}\subset\fg_{-\alpha}$ 
form dual bases with respect to $(-,-)$. Viewing
$u\in\cS(\fg_{-\alpha}+\fh)$ as a polynomial in $\{f_i\}$ 
we obtain $(\ad e_i)(u)=h\frac{\partial{u}}{\partial{f_i}}$.
Thus $(\ad e_i)(u)=0$ for all $i$ gives $u\in\cS(\fh)$ as required.
\end{rem}

\section{The centre of a $Q$-type Lie superalgebra}
\label{sectcentre}
In this section we describe the centre of a $Q$-type Lie superalgebra
(see~\Thm{thmkac}). The central elements
correspond to the polynomials in $\Sh$ which have the same values
at $\lambda$ and $\lambda'$ provided that  $\lambda$ is subregular
and $\lambda'$ is the highest weight of $\ol{M(\lambda)}$.

We also show that
$\cZ(\fq(n))=\cZ(\fsq(n))$ and $\cZ(\fpq(n))=\cZ(\fpsq(n))$ 
(see~\Cor{centcon}).
Notice that $\cZ(\fq(n))$ was described in~\cite{sq}.

Throughout  the section $\fg$ is a $Q$-type Lie superalgebra and 
$\fg\not=\fpq(2),\fpsq(2)$. 

\subsection{}
\begin{thm}{thmkac}
Let $\fg$ be a $Q$-type Lie superalgebra, $\fg\not=\fpq(2),\fpsq(2)$.
The restriction of $\HC$ to $\Zg$ is an algebra isomorphism $\Zg\iso Z$
where $Z$ is the set of $W$-invariant polynomial functions on $\fh_{\ol{0}}^*$
which are constant along each straight line  parallel to a root
$\alpha$ and lying in the hyperplane $h_{\ol{\alpha}}(\lambda)=0$. 
In other words,
$$Z:=\Sh^W\cap\bigcap_{\alpha\in\Delta} Z_{\alpha},$$
where
$$Z_{\alpha}:=\{f\in\Sh|\ 
h_{\ol{\alpha}}(\lambda)=0\ \Longrightarrow\ f(\lambda)=f(\lambda-c\alpha)
\ \forall c\in\mathbb{C}\}.$$
\end{thm}

We prove this theorem using a modification of Kac
method presented in~\cite{klapl}. Actually we prove that
$\HC(\cZ(\hat{U}))=Z$ (see Section~\ref{secthat} for the definition
of $\hat{U}$) and then use the equality $\cZ(\hat{U})=\Zg$
obtained in~\Cor{corbern}.

\subsection{Proof of~\Thm{thmkac}}
By~\Thm{bern} the restriction of $\HC$ to $\cZ(\hat{U})$ is injective.
For $z\in\cZ(\hat{U})$ the image $\HC(z)$ is central in $\Uh$ so
$\HC(z)\in\Sh$. By~\Prop{ua}, 
$z$ acts on $M(\lambda)$ by $\HC(z)(\lambda)$.
Using~\Cor{primv} we conclude
that $\HC(\cZ(\hat{U}))\subseteq Z$. 
To prove the  opposite inclusion $\HC(\cZ(\hat{U}))\supseteq Z$ 
we assign to each $\phi\in Z$
an element $z\in\hat{U}$ with the property: $\HC(z)=\phi$ and
$zv=v\phi$ for all $v\in \cM$. By~\Prop{ua}, $z$ is central.

We construct the element $z=\sum z_{\nu}$ by a recursive procedure;
the summands $z_{\nu}\in J'(\nu)$ (see~\ref{J'} for the notation) 
are chosen to fulfill the condition
$\sum_{\mu\leq\nu} z_{\mu}v=v\phi$ for all  $v\in \cM_{-\nu}$.
For any $\nu\in Q^+$ and any $v\in \cM_{-\nu}$ one has
 $zv=\sum_{\mu\leq\nu} z_{\mu}v$. Hence $zv=v\phi$ for all $v\in\cM$.

Putting $z_{<\nu}:=\sum_{\mu<\nu} z_{\mu}$ we can rewrite the
above condition as
\begin{equation}
\label{pr*}
z_{\nu}v=v\phi-z_{<\nu}v, \ \forall v\in \cM_{-\nu}.
\end{equation}
In the rest of the proof we show the existence
of $z_{\nu}$ satisfying~(\ref{pr*}).

\subsubsection{}
\label{iota-mu}
The term $z_{\nu}$ lies in $\cU(\fb^-)_{-\nu}\otimes_{\fh}\cU(\fb^+)_{\nu}$.
Identify $\cU(\fb^-)$ with $\Ind(R)=\cM$ 
 and $\cU(\fb^+)$ with $\Ind_+(R)$.
Under these identifications, $z_{\nu}$ lies in 
$\cM_{-\nu}\otimes_{R}\Ind_+(R)_{\nu}$.
The action of $\cM_{-\nu}\otimes_{R}\Ind_+(R)_{\nu}$
on $\cM_{-\nu}$ takes form
\begin{equation}
\label{S*}
(a\otimes b)v=a\HC(\sigma(b)v)=a S^*(b)(v)
\end{equation}
where $S^*: \Ind_+(R)\to \Hom_{R_r}(\cM, R^{\sigma})$
is induced by the Shapovalov map
$S: \cM\to \Hom_{R_r}(\Ind_+(R), R^{\sigma})$.
Consider the chain of homomorphisms
$$\begin{diagram}
\cM\otimes_{R}\Ind_+(R)& \rTo^{\id\otimes S^*} &
\cM\otimes_{R_r} \Hom_{R_r}(\cM, R^{\sigma}) & \rTo^{\iota} & \End_{R_r}(\cM)
\end{diagram}$$
where $\iota$ is the natural map ($\iota(v\otimes f)(v'):=vf(v')$).
Let $\psi$ be the composed map $\psi:=\iota\circ(\id\otimes S^*)$
and  $\psi_{\nu}: \cM_{-\nu}\otimes_{R}\Ind_+(R)_{\nu}\to 
\End_{R_r}(\cM_{-\nu})$ be the restriction of $\psi$.
In the light of~(\ref{S*}),
an element $x \in \cM_{-\nu}\otimes_{\fh}\Ind_+(R)_{\nu}$
acts on $v\in \cM_{-\nu}$ by the formula
$$x v=\psi_{\nu}(x)(v).$$

Thus the existence of $z_{\nu}$ satisfying~(\ref{pr*}) is equivalent to
the inclusion $C\in \im\psi_{\nu}$ where
$C\in\End (\cM_{-\nu})$  is given by 
$Cv=v\phi-z_{<\nu}v$;
notice that $C\in \End_{R_r}(\cM_{-\nu})$ because $\phi\in\Sh$
belongs to the centre of $R=\Uh$.
The condition $C\in \im\psi_{\nu}$ can be rewritten
as $\im C^*\subset \im S^*_{\nu}$
where $C^*\in\End_{\fh}(\Hom_{R_r}(\cM_{-\nu}, R^{\sigma}))$
is transpose to $C$. Thus it remains to verify the
inclusion 
\begin{equation}\label{C*S*}
\im C^*\subset \im S^*_{\nu},\ \text{ where } Cv=v\phi-z_{<\nu}v.
\end{equation}

\subsubsection{}
Since $S^*_{\nu}$ and $C^*$ are linear maps, 
the inclusion~(\ref{C*S*}) 
is equivalent to the linear equation $C^*=S^*_{\nu}X$
over the polynomial algebra $\Sh$. Rewrite $C^*=S^*_{\nu}X$
as $YS_{\nu}=C$ for $Y:=X^*$.

The equation $YS_{\nu}=C$ has a solution over the field of fractions of $\Sh$,
because $S_{\nu}$ is a monomorphism between free  $\Sh$-modules of the same 
finite rank and thus it is invertible over $F:=\Fract \Sh$.
View $\fh_{\ol{0}}$ as an affine space $\mathbb{C}^n$.
Let us check that  $Y:=CS_{\nu}^{-1}$ is regular.

Retain terminology of~\ref{gamm}. 
If $\lambda\in\fh_{\ol{0}}^*$ 
is a regular point, $S(\lambda)$ is bijective and so
$Y$ is regular at $\lambda$. Since
the union of regular and subregular 
points in $\fh^*_{\ol{0}}$ is a set of codimension two, 
it is enough to verify that $Y$ is regular in a neighbourhood of 
any subregular point. Here we need the following lemma.

\subsubsection{}
\begin{lem}{DE}
Let $D=(d_{ij})_{i=1,N}^{j=1,N}$ and $E=(d_{ij})_{i=1,N}^{j=1,M}$
be two matrices, where $d_{ij},e_{ij}$ are functions 
in $z_1,\ldots,z_m$ which are regular on some neighbourhood $U$ of the origin.
Put $V:=U\cap\{z_1=0\}$. Suppose that $D$ is invertible on $U\setminus V$ 
and that for any $\lambda\in V$ one has

(a) The order of zero of $\det D$ at the point $\lambda$
is equal to $\dim\Ker D(\lambda)$,

(b) $\Ker D(\lambda)\subset \Ker E(\lambda)$.

Then $ED^{-1}$ is regular on $U$.
\end{lem}

The first assumption means that all poles of $D^{-1}$ at $\lambda$
have order one.
The proof is completely similar to one given in~\cite{klapl}.

\subsubsection{}
Let $\lambda\in\fh_{\ol{0}}^*$ be a subregular point. 
In the light of~\ref{ordzero},
the first assumption of the lemma for the matrix
$B_{\nu}$ follows from~\Cor{thmcrk}.
Recall
that $B=\int\circ S$ 
where $\int$ is an invertible map (see~\Lem{print} (ii)).
Hence the first assumption holds for $S_{\nu}$.

For the second assumption, recall that $\Ker S_{\nu}(\lambda)=
\ol{M(\lambda)}_{\lambda-\nu}$. By~\Prop{length}, 
$\ol{M(\lambda)}\cong V^{\oplus}(\tilde{\lambda})$ 
for some $\tilde{\lambda}<\lambda$. The condition
$\phi\in Z$ ensures that $\phi(\lambda)=\phi(\tilde{\lambda})$.
Take $v\in\ol{M(\lambda)}_{\lambda-\nu}$. One has
$v\in V^{\oplus}(\tilde{\lambda})_{\tilde{\lambda}-\nu'}$ where 
$\lambda-\nu=\tilde{\lambda}-\nu'$.
Therefore $z_{<\nu}v=z_{\leq \nu'}v$. Since $\nu'<\nu$ 
the induction hypothesis gives
$z_{\leq \nu'}v=v\phi=\phi(\tilde{\lambda})v$.
Thus $z_{<\nu}v=\phi(\lambda) v$ and so
$v\in \Ker C(\lambda)$. This implies
the second assumption of~\Lem{DE}.

Finally, $Y$ is regular and this completes the proof of~\Thm{thmkac}.
\qed

\subsection{}
\begin{cor}{centcon}
\begin{enumerate}
\item
$\cZ(\fq(n))=\cZ(\fsq(n))$.

\item $\cZ(\fpq(n))=\cZ(\fpsq(n))$.

\end{enumerate}
\end{cor}
\begin{proof}
Observe that $\fq(n)=\fsq(n)\oplus\mathbb{C}H$ where $H$ is odd
and $[H,\fsq(n)]\subset\fsq(n)$ (in the notation of ~\ref{hi} 
$H=H_1+\ldots+H_n$).

Take $z\in\cZ(\fq(n))$ and write $z=a+bH$ where $a,b\in\cU(\fsq(n))$.
Recall that $\cZ(\fq(n))$ is pure even and so $a$ is even and $b$ is odd.
For any $x\in\fsq(n)$ one has
$$0=[x,a+bH]=[x,a]\pm b[x,H]+[x,b]H.$$
Since $[x,a]\pm b[x,H]\in \cU(\fsq(n))$ and
$[x,b]H\in \cU(\fsq(n))H$, we conclude 
$$[x,b]=0.$$
On the other hand,
$$0=[H,a+bH]=[H,a]\pm 2bH^2+[H,b]H$$
which gives 
$$[H,b]=0.$$
Therefore $b\in\cZ(\fq(n))$ and so $b=0$ because $b$ is odd.
Hence $z=a\in\cU(\fsq(n))$. This shows that 
$\cZ(\fq(n))\subset \cU(\fsq(n))$
and so $\cZ(\fq(n))\subset \cZ(\fsq(n))$. 
The even parts of Cartan subalgebras of $\fq(n)$ and $\fsq(n)$
coincide.
By~\Thm{thmkac} one has $\HC\bigl(\cZ(\fq(n))\bigr)=
\HC\bigl(\cZ(\fsq(n))\bigr)$. 
This proves (i). 
The proof of (ii) is completely similar.
\end{proof}

\appendix\section{The algebra $\Uh$}
\label{appn}
In this section we study the algebra $\Uh$ which is a Clifford algebra
over $\Sh$. In~\ref{ZAh} we describe the centre and the anticentre of $\Uh$.
In~\ref{clla} we recall some basic facts on complex Clifford algebras.
In~\ref{intdef} we
introduce a map $\int:\Uh\to\Sh$. In~\ref{Rsigma} we compare various
constructions of dual $R$-modules. In~\ref{rednorm} we
introduce a reduced norm of an endomorphisms of $R$-module.
This is an analogue of the reduced norm 
 for endomorphisms of modules over an Azumaya algebra (see~\cite{kn}).

In~\ref{ZAh}---~\ref{rednorm} we put $A=\Sh$ and $R=\Uh$. Set 
$$n:=\dim\fh_{\ol{1}}.$$

\subsection{Notation}
Let $A$ be a commutative ring, $M$ be a projective finitely generated
 $A$-module endowed with a quadratic form $q:M\to A$. 
{\em The Clifford algebra} $\Cl(M,q)$
defined by these data is the $A$-algebra generated by $M$ with the relations:
$x^2=q(x),\ x\in M$. The Clifford algebra corresponding to a 
non-degenerate form is called non-degenerate.
We view $\Cl(M,q)$ as a superalgebra by
letting the elements of $M$ to be odd.

\subsubsection{}
One has $\Cl(M\oplus M',q+q')=\Cl(M,q)\otimes_A\Cl(M',q')$
where the right-hand side is the tensor product of superalgebras.

In what follows $M$ will be a free $A$-module of finite rank. 
In this case $\Cl(M,q)$ is free as $A$-module and admits a PBW basis.
We say that a Clifford algebra $\Cl(M,q)$ has rank $r$ if $M$ has
rank $r$ over $A$. 

All $\Cl(M,q)$-modules we consider below
are assumed to be $\mathbb{Z}_2$-graded and free over $A$.

\subsubsection{}
\label{AR}
Set $A:=\cU(\fh_{\ol{0}}),\ R:=\Uh$. Denote by $\sigma$ the antiautomorphism
defined by $x\mapsto -x$ for $x\in\fh$. It provides
 an equivalence between
the categories of left and right $R$-modules.

We can view $R$ as a Clifford algebra over $A$: $M$ is
a free $A$-module spanned by $\fh_{\ol{1}}$ and $q:M\to A$
is given by $q(H)=\frac{1}{2}[H,H]$ for $H\in\fh_{\ol{1}}$
(the corresponding bilinear form is
given on $\fh_{\ol{1}}$ by the formula $B(H,H')=[H,H']$).
The antiautomorphism $\sigma$ defined above coincides
with the restriction of 
$\sigma:\Ug\to\Ug$ defined in~\ref{sigma}.

\subsection{The anticentre of $\Uh$}
\label{ZAh}
Retain notation and definitions of Sect.~\ref{sectantie}.
It is easy to see that $\Zh=\Sh$. We describe the anticentre $\Ah$
below.

\subsubsection{}
\label{TUh}
In the notation of~\ref{antie}, $T_{\fh}=(\ad' u_0)(1)$ and
$\cA(\fh)=(\ad' u_0)(\Sh)$ for some $u_0\in \Uh$.
Any $z\in\Sh$ is central in $\Uh$ and so $(\ad' u_0)(z)=z(\ad' u_0)(1)$.
Therefore
$$\Ah=\Sh T_{\fh}.$$

\subsubsection{}
\label{Thq}
Retain notation of~\ref{hi}.
The elements $H_1,\ldots,H_n$ form a basis of $\fh_{\ol{1}}$
for $\fg=\fq(n),\fpq(n)$ and 
$H_1-H_2,H_2-H_3,\ldots, H_{n-1}-H_n$ form a basis of $\fh_{\ol{1}}$
for $\fg=\fq(n),\fpsq(n)$.
We can put
\begin{equation}
\label{Tfh}
T_{\fh}:=\left\{\begin{array}{ll}
H_1\ldots H_n & \text{ for } \fg=\fq(n),\fpq(n)\\
\sum_{i=1}^n (-1)^i H_1\ldots \hat{H_i} \ldots H_n
& \text{ for } \fg=\fsq(n),\fpsq(n)
\end{array}\right.
\end{equation}
since the right-hand side has the degree $\dim\fh_{\ol{1}}$
and is invariant with respect to the twisted adjoint action
of the basis elements;  the invariance easily follows from the formulas
$$(\ad' H_i)(H_{i_1}\ldots H_{i_r})=\left\{\begin{array}{ll}
0 & \text{ if }i\in\{i_1,\ldots,i_r\},\\
2H_iH_{i_1}\ldots H_{i_r} & \text{ if }i\not\in\{i_1,\ldots,i_r\}.
\end{array}\right.$$
One has
$$
t_{\fh}:=T^2_{\fh}=\left\{\begin{array}{ll}
\pm h_1\ldots h_n & \text{ for } \fg=\fq(n),\fpq(n)\\
\pm\sum h_1\ldots \hat{h_i}\ldots h_n
& \text{ for } \fg=\fsq(n),\fpsq(n).
\end{array}\right.
$$
 
Observe that $t_{\fh}$ is equal to the determinant
of the bilinear form $B$ introduced in~\ref{AR}.

\subsubsection{The centre of $\Uh$}
\label{ZZ}
View $\Uh$ as a non-graded algebra and denote its centre by $Z$.
The definition of the anticentre immediately gives $Z=\Zh_{\ol{0}}\oplus
\Ah_{\ol{1}}$. Thus $Z=\cU(\fh_{\ol{0}})$ if $n$ is even and $Z=\cU(\fh_{\ol{0}})
\oplus \cU(\fh_{\ol{0}}) T_{\fh}$ if $n$ is odd.

\subsection{Clifford algebras over $\mathbb{C}$.}
\label{clla}
Set $n:=\dim\fh_{\ol{1}}, A:=\Sh, R:=\Uh$.
For  $\lambda\in\fh^*_{\ol{0}}$   denote by $I(\lambda)$ the maximal
ideal of $A$ corresponding to $\lambda$. Set
$$\Cl(\lambda):=R/(R I(\lambda)).$$

Denote by $\sMat_{r,s}(\mathbb{C})$  the superalgebra of the endomorphisms
of a superspace of dimension $r+s\epsilon$.
The elements of $\fq(n)$ forms a subalgebra
in the superalgebra $\sMat_{n,n}(\mathbb{C})$; in this section
we denote this (associative) algebra by $Q(n)$.

\subsubsection{}
\label{clpi'}
Let $\Cl(m)$ be the standard  Clifford algebra: it
is generated by $\xi_1,\ldots,\xi_m$ subject to the relations
$\xi_i^2=1, \xi_i\xi_j+\xi_j\xi_i=0$ for $i\not=j$.

If the symmetric form
$B(\lambda): (H,H')\mapsto\lambda([H,H'])$ is non-degenerate 
then $\Cl(\lambda)$
is isomorphic to the standard complex Clifford algebra generated by the
image of $\fh_{\ol{1}}$. The observation in~\ref{Thq} gives
$$\Cl(\lambda)\cong\Cl(\dim\fh_{\ol{1}})\ \text{ if } t_{\fh}(\lambda)\not=0.$$

\subsubsection{}
\label{clpi}
It is well known that $\Cl(m)$ is a simple superalgebra:
it has at most two simple
modules which differ by grading ($E$ and $\Pi(E)$)
and all graded $\Cl(m)$-modules are completely reducible.

If $m$ is even then $E$ is simple as a non-graded module, 
$\dim E=2^{\frac{m}{2}}, \dim E_{\ol{0}}=\dim E_{\ol{1}}$. One has $E\not\cong\Pi(E)$ and
$\Cl(m)=\End (E)=\End (\Pi(E))$. Thus, $\Cl(m)$ is
isomorphic to the superalgebra $\sMat_{r,r}(\mathbb{C})$ 
where $r:=2^{\frac{m}{2}-1}$.

If $m$ is odd then  $E$ is not simple as a non-graded module,
$\dim E=2^{\frac{m+1}{2}}$ and $E\cong \Pi(E)$. Put 
$$r:=2^{\frac{m-1}{2}}.$$
We claim that image of $\Cl(m)$ in $\End (E)\cong\sMat_{r,r}(\mathbb{C})$
coincides with  $Q(r)$ in a suitable basis. Indeed, let $\dim\fh_{\ol{1}}=m$. Take
$\lambda$ such that $t_{\fh}(\lambda)\not=0$. Then $\Cl(\lambda)\cong\Cl(m)$
and the image of $T_{\fh}$ in $\Cl(\lambda)$ is a central
odd element whose square is the non-zero scalar $t_{\fh}(\lambda)$.
Therefore the centre of $\Cl(m)$ contains an odd element $z$ such that
$z^2=1$. Clearly, we can choose a basis in $E$ in such a way that
the matrix corresponding to $z$ is $X_{0,\id}$ 
in the notation of~\ref{qtype}. It is easy to verify that 
$$\{Y\in\sMat_{r,r}(\mathbb{C})|\ YX_{0,\id}=X_{0,\id}Y\}=Q(r).$$
Since $\dim \Cl(m)=2^m=\dim Q(r)$ we conclude
$$\Cl(m)=Q(r).$$

\subsubsection{}
\label{clpi0}
Take an arbitrary $\lambda\in\fh^*_{\ol{0}}$.
One easily sees that 
$\Cl(\lambda)\cong \Cl(m)\otimes \bigwedge(\Ker B(\lambda))$
where $m:=\dim\fh_{\ol{1}}-\dim \Ker B(\lambda)$.

As a result, $\Cl(\lambda)$ has at most two simple
modules ($E$ and $\Pi(E)$) which are simple as
$\Cl(m)$-modules (where $\Cl(m)=\Cl(m)\otimes 1\subset  
\Cl(m)\otimes \bigwedge(\Ker B(\lambda))\cong \Cl(\lambda)$). 

View $\Cl(\lambda)$ as a right module over itself 
via the multiplication. One has $\Pi(\Cl(\lambda))\cong\Cl(\lambda)$
iff $\lambda\not=0$. Indeed, if $\lambda\not=0$ there exists $H\in\fh_{\ol{1}}$
satisfying $H^2=1$; the map $x\mapsto Hx$ provides an isomorphism
$\Cl(\lambda)\iso \Pi(\Cl(\lambda))$. On the other hand,
$\Cl(0)=\bigwedge(\fh_{\ol{1}})$ has a one-dimensional socle $\bigwedge^{\top}\fh_{\ol{1}}$
which forces $\Pi(\Cl(0))\not\cong\Cl(0)$.

\subsubsection{}
\label{cliqn}
Any simple $R$-module is annihilated by $I(\lambda)$ for some
$\lambda\in\fh^*_{\ol{0}}$.
We see that for each $\lambda\in\fh^*_{\ol{0}}$ there are at most
two simple $R$-modules annihilated by $I(\lambda)$; we
may (and will) denote them by $E(\lambda),\Pi(E(\lambda))$.
The total dimension $\dim E(\lambda)$ depends
on the rank $m$ of the evaluated bilinear
form $B(\lambda)$.
For the case $\fg=\fq(n)$  one has
$m=n-|\{i: \lambda(h_i)=0\}|$.

\subsection{The map $\int$}
\label{intdef}
In this subsection we construct a linear map $\int:R\to A$
which satisfies the properties (i)-(v) of~\ref{print}.
These properties ensure that $\cB(x,y):=\int\! xy$ is
a non-degenerate invariant  bilinear form $\cB:R\otimes_A R\to A$
(the invariance means that
$\cB([x,y],z)=\cB(x,[y,z])$). This form is even (resp., odd)
if $n$ is even (resp., odd). Moreover supersymmetric that is
$\cB(x,y)=(-1)^{p(x)p(y)}\cB(y,x)$ (if $\cB$ is odd
it is equivalent to the symmetricity).
For each $\lambda\in\fh^*_{\ol{0}}$
the evaluated form 
$\cB_{\lambda}:\Cl(\lambda)\otimes_{\mathbb{C}}\Cl(\lambda)\to\mathbb{C}$
is a non-degenerate invariant bilinear form.

\subsubsection{}
\label{express}
In order to express the form $\cB$ more explicitly, we calculate
its evaluation in the case when $\Cl(\lambda)$ is non-degenerate. 

Recall that if $\Cl(\lambda)$ is non-degenerate then
$\Cl(\lambda) \cong\sMat_{r,r}(\mathbb{C})$ for even $n$ and
$\Cl(\lambda) \cong Q(n)$ for odd $n$. 
The superalgebra $\sMat_{r,r}(\mathbb{C})$ has an even non-degenerate 
invariant bilinear form $\cB'(X,Y):=\str XY$. 
The superalgebra $Q(r)$ has an odd non-degenerate invariant
bilinear form $\cB''(X,Y):=\tr' XY$ (see~\ref{qtype}
for notation). We show that for $\Cl(\lambda)$ being non-degenerate,
$\cB'=c'(\lambda)\cB_{\lambda}$ where $c'(\lambda)^2=4t_{\fh}(\lambda)$
if $n$ is even, and $\cB''=c''(\lambda)\cB_{\lambda}$ 
where $c''(\lambda)^2=t_{\fh}(\lambda)$ if $n$ is odd.
Recall that $\Cl(\lambda)$ is
non-degenerate iff $t_{\fh}(\lambda)\not=0$.

\subsubsection{}
\label{integral}
The Clifford algebra $R=\Uh$ over $A=\Sh$ has the canonical filtration:
$$\cF^0(R):=A,\ \cF^1(R):=A\fh_{\ol{1}},\ \cF^i(R):=(\cF^1(R))^i.$$
The associated graded algebra $\gr R$
is a commutative superalgebra $A\otimes\bigwedge\fh_{\ol{1}}$ with the grading:
$(\gr R)_i=A\bigwedge^i\fh_{\ol{1}}$. Take $T_{\fh}$ as in~\ref{Thq} and 
observe that $T_{\fh}\not\in\cF^{n-1}(R)$. 
There exists a unique $A$-homomorphism
$\int\!:R\to A$  such that 
$$\int(T_{\fh})=1,\ \ \ \Ker\int=\cF^{n-1}(R).$$
Since $T_{\fh}$ has the same parity as $n$,
the map $\int$ is even (resp., odd) if $n$ is even (resp., odd).

We will write $\int\!u$ instead of $\int(u)$.

Notice that  
$$\int\!u=f(\gr u)$$ where 
$f:A\otimes\bigwedge\fh_{\ol{1}}\to A$
is the $A$-homomorphism given by 
$$f(\gr T_{\fh})=1,\ \ \ \Ker f=\sum_{i=0}^{n-1}(\gr R)_i =
A\otimes\sum_{i=0}^{n-1}\bigwedge\fh_{\ol{1}}.$$

\subsubsection{}
For each $\lambda\in\fh^*_{\ol{0}}$ the evaluated map
$u\mapsto(\int\! u)(\lambda)$ induces a linear map 
$\Cl(\lambda)\to\mathbb{C}$ which we also denote by $\int$. Thus
$$\int\! u(\lambda):=(\int\! u)(\lambda)$$
where $u(\lambda)$ stands for the image of $u\in R$ in $\Cl(\lambda)$.

\subsubsection{}
Define the bilinear form $\cB:R\otimes_{A} R\to A$
by setting 
$$\cB(u,u'):=\int uu'.$$
Clearly, $\cB$ is even (resp., odd) if $n$ is even (resp., odd).
\Lem{print} below shows that $\cB$ is a non-degenerate
invariant form which is supersymmetric. It also shows
that the evaluated form  $\cB_{\lambda}$ is non-degenerate
{\em for all} $\lambda$.

\subsubsection{}
Clearly, $\int(H'_1\ldots H'_n)\in\mathbb{C}^*$
if $H'_1,\ldots, H'_n$ is any basis of $\fh_{\ol{1}}$. 
Choose a basis
$\{H'_i\}_{i\in I}$ in such a way that $\int(H'_1\ldots H'_n)=1$.
Set $I:=\{1,\ldots,n\}$.
For $J\subset I$ set $H_J:=\prod_{j\in J}H'_j$ ($H_{\emptyset}=1$)
where the factors are arranged by increasing of indices. 
The elements $H_J$ form a system of free generators of $R$ over $A$.

\subsubsection{}

We claim that
\begin{equation}
\label{intHJ}
\int\! \sigma(H_J)H_{J'}=\left\{\begin{array}{ll}
\pm 1, &
\text{ if }
J'=I\setminus J\\
0 & \text{ otherwise}.\end{array}\right.
\end{equation}
Indeed, $\int\! \sigma(H_J)H_{J'}=f(\gr\sigma(H_J)H_{J'})
=f(\gr\sigma(H_J)\gr H_{J'})$. Evidently, 
$\gr\sigma(H_J)\gr H_{I\setminus J}=\pm\gr H_I$  and
$\gr\sigma(H_J)\gr H_{J'}\in\Ker f$ for
$J'\not=I\setminus J$. 
 
\subsubsection{}
\begin{lem}{print}  For all $a,b\in R$ one has
\begin{enumerate}
\item $\int [a,b]=0$.
\item If $u\in R$ is such that $\int\! ua=0$ for all $a\in R$ then $u=0$.
\item If $u\in\Cl(\lambda)$ is such that $(\int\! ua)(\lambda)=0$ 
for all $a\in\Cl(\lambda)$ then $u=0$.
\item
$\int[a,b]c=\int\! a[b,c]$.
\item $\int\! \sigma(a)=(-1)^n\int\! a$.
\end{enumerate}
\end{lem}
\begin{proof}
The assertion (i) follows from the formula $\int\! a=f(\gr a)$
and (iv) is a reformulation of (i). The formula~(\ref{intHJ}) imply
both (ii) and (iii). The last assertion is an immediate consequence
of the formula $\sigma(T)=(-1)^nT$ which can be easily verified.
\end{proof}

\subsubsection{}
\label{propor}
Retain notation of~\ref{express}.
Take $\lambda$ such that $t_{\fh}(\lambda)\not=0$
and thus $\Cl(\lambda)$ is non-degenerate.
Let $T(\lambda)\in\Cl(\lambda)$ be the image of $T_{\fh}$.
Recall that $T(\lambda)$ commutes with the even elements
of $\Cl(\lambda)$ and that
$T(\lambda)^2=t_{\fh}(\lambda)\in \mathbb{C}$. 

Let $n$ be even. One can easily sees from~\Lem{print} (i) that
the evaluation of $\int$ on $\Cl(\lambda)\cong\sMat_{r,r}(\mathbb{C})$
is proportional to the supertrace.
The element $T(\lambda)$ belongs to the centre of the algebra
$\sMat_{r,r}(\mathbb{C})_{\ol{0}}=\Mat_r(\mathbb{C})\times \Mat_r(\mathbb{C})$.
Since $T(\lambda)\not\in \mathbb{C}$ and 
$T(\lambda)^2=t_{\fh}(\lambda)\in\mathbb{C}$ one has
$T(\lambda)=c\id\times (-c\id)$ where $c^2=t_{\fh}(\lambda)$
and $\id$ is the identity matrix in
$\Mat_r(\mathbb{C})$. In particular, $\str T(\lambda)=2c$.
Since $\int T(\lambda)=1$, we conclude $2c\int=\str$ and 
so $\cB'=2c\cB$.

Let $n$ be odd. It is easy to deduce from~\Lem{print} (i) that
the evaluation of $\int$ on
 $\Cl(\lambda)\cong Q(r)$ is proportional to the map $\tr'$.
The element $T(\lambda)\in Q(r)_{\ol{1}}$ belongs to
the centralizer of $Q(r)_{\ol{0}}$. 
Consequently, $T(\lambda)=c\id'$ where $c^2=t_{\fh}(\lambda)$
and $\id'\in Q(r)_{\ol{1}}$ corresponds to 
the identity matrix in $Q(r)_{\ol{1}}\cong\Mat_r(\mathbb{C})$.
Since $\tr' T(\lambda)=c$ we obtain $c\int=\tr'$ and thus
$\cB''=c\cB$.

\subsection{Realization of $N^*$}
\label{Rsigma}

Define on $R$ a new bimodule structure $R^{\sigma}$
via 
$$v.r:=(-1)^{p(r)p(v)}\sigma(r)v,\quad r.v:=(-1)^{p(r)p(v)}v\sigma(r)$$
where the dot stands for the new action, 
$\sigma$ is the antiautomorphism introduced in~\ref{sigma},
$r$ is an element of the algebra $R$ and
$v\in R^{\sigma}$.  

Let $U$ be a superalgebra containing $R$ such that 
the antiautomorphism $\sigma$ can be extended to $U$.
Let $N$ be a bimodule over $U$-$R$. Let
$N':=\Hom_{R_r}(N,R^{\sigma})$ be the set
of homomorphisms of {\em right} $R$-modules. 
Notice that $N'$ has the natural structure of $R$-$U$-bimodule
which we convert to $U$-$R$-bimodule structure via the antiautomorphism
$\sigma$.

\subsubsection{}
\begin{lem}{}
\begin{enumerate}
\item
The map $R^{\sigma}\to\Hom_{R_r}(R,R)$ given by $t\mapsto l_t$ where 
$l_t(r):=tr$ is  an even $R$-bimodule isomorphism.
\item
The map $R\to\Hom_{R_r}(R,R^{\sigma})$ given by
$t\mapsto l'_t$ where $l'_t(r):=(-1)^{p(r)p(t)}\sigma(r)t$
is an even $R$-bimodule isomorphism.
\end{enumerate}
\end{lem}
The proof is straightforward.

\subsubsection{}
View $N^*:=\Hom_{A_r}(N,A)$ as a $U-R$ bimodule via $\sigma$
where by $\Hom_{A_r}$ we mean the set of homomorphisms of {\em right}
$A$-modules.
Convert $U-R$ bimodules to right $(U\otimes R)$-modules
using $\sigma$. 
Define similarly $R$-bimodule structure on $R^*:=\Hom_A(R,A)$.

\subsubsection{}
\begin{lem}{}
\begin{enumerate}
\item
The map $t\mapsto \int\!l_t$ where $\int\!l_t(r):=\int(tr)$
is an isomorphism of $R$-bimodules
$R^{\sigma}\iso R^*$ if $\int$ is even and
$R^{\sigma}\iso\Pi(R^*)$ if $\int$ is odd. 

\item
Let $N$ be a $U$-$R$ bimodule.
The map 
$\psi\mapsto\int\psi$ where $(\int\psi)(x):=\int(\psi(x))$
is an isomorphism of $U-R$-bimodules
$\Hom_{R_r}(N,R^{\sigma})\to N^*$ if $\int$ is even and
$\Hom_{R_r}(N,R^{\sigma})\to \Pi(N^*)$ if $\int$ is odd.
\end{enumerate}
\end{lem}
\begin{proof}
(i) is straightforward. Now (ii) follows from the Frobenius
reciprocity. \end{proof}

\subsection{Reduced norm}\label{rednorm}
Let $A$ be a polynomial algebra and $R:=\Cl(M,q)$ be a Clifford algebra
viewed as a non-graded algebra.
Suppose that $M$ is a free $A$-module
of an even rank $2n$ and that the kernel 
of bilinear form corresponding to $q$
is zero. 
Let $N$ be an $R$-module which is free of finite rank over $A$.

We will construct a {\em reduced norm} i.e.,
a map $\Norm:\End_{R}(N)\to A$ satisfying the properties
\begin{equation}\label{normprop}
\Norm(\id)=1,\ \ \Norm(\psi\psi')=\Norm(\psi)\Norm(\psi'),
\ \ \ \Norm(\psi)^{2^n}=\det\psi
\end{equation}
where in the last formula
$\psi$ is viewed as an element of  $\End_{A}(N)$ and
$\det: \End_{A}(N)\to A$ is the determinant map.
This map is  an analogue of the reduced norm 
for endomorphisms of module over an Azumaya algebra (see~\cite{kn}).

\subsubsection{}
Set $F':=\Fract A$. The algebra $R_{F'}:=R\otimes_A F'$ is an
 Azumaya algebra over $F'$ (see~\cite{kn}). Therefore, for
a suitable Galois  extension $F/F'$, 
the algebra $R_F:=R\otimes_A F$ is isomorphic to
the matrix algebra $\Mat_{r}(F)$ where $r:=2^n$.
The module $N_{F}:=R_{F}\otimes_A N$ is an $R_F$-module
and so there is an isomorphism
$$N_{F}\iso E\otimes_F V$$
where $E$ is the simple module over the matrix algebra
$\Mat_{r}(F)$ and $V$ is a finite dimensional vector space
over $F$. Since $E$ is simple over $R_F$ any $R_F$-endomorphism
of $N_F=E\otimes_F V$ takes form $\id_E\otimes\phi'$
where $\phi\in\End(V)$. In this way, we obtain
the algebra isomorphism $\gamma:\End_{R_F}(N_F)\iso\End(V)$. 
Viewing $\End_{R}(N)$ as an $A$-subalgebra of $\End_{R_F}(N_F)$
 and set
$$\Norm\phi:=\det\gamma(\phi).$$
Clearly, $\Norm$ satisfies the first two properties of~(\ref{normprop}).
The last property follows from the formula $\dim E=r$.
It remains to verify the

\subsubsection{}
\begin{lem}{}
$\Norm\phi\in A$.
\end{lem}
\begin{proof}
The Galois group $G$ acts on $F,R_F,N_F$ leaving elements
of $F',R_{F'}$ and $N_{F'}$ stable. 
This induces the action of $G$ on $\End_{R_F}(N_F)$.
On the other hand, $G$ acts naturally on $\End(V)$.
The map $\gamma$ does not commute with these $G$-action; however,
for each $g\in G$ the commutator  $\gamma g^{-1}\gamma^{-1}g$
is an automorphism of the matrix algebra $\End(V)$ which leaves
elements of $F$ stable. Therefore $\gamma g^{-1}\gamma^{-1}g$
is an inner automorphism of $\End(V)$ and thus
for any $\psi\in\End(V)$ one has 
$\det(\gamma g^{-1}\gamma^{-1}g (\psi))=\det\psi$.
In other words, for any $\phi\in \End_{R_F}(N_F)$ one has
$$\det \gamma(g\phi)=\det g\gamma(\phi)=g(\det\gamma(\phi)).$$
Take $\phi\in\End_R(N)$. One has $g(\phi)=\phi$ for all $g\in G$.
This means that $\Norm\phi=\det\gamma(\phi)$ is $G$-stable that is
$\Norm\phi\in F'$. Observe that $\det\phi\in A$
and so, by~(\ref{normprop}) $(\Norm\phi)^r\in A$.
Since $A$ is integrally closed in $F'$, $\Norm\phi\in A$ as required.
\end{proof}

%%%%%%%%%%%%%%%%  biblio.tex


\begin{thebibliography}{MMMMM}

\bibitem[BB]{bb} A.~Beilinson, J.~Bernstein, {\em A
proof of Jantzen conjectures}, I.~M.~Gelfand Seminar, 1--50,
Adv. Soviet Math. 16, Part 1, AMS, Providence, RI, 1993.

\bibitem[DK]{dk} C.~De Concini, V.~Kac, {\em Representations
of quantum groups at roots of $1$}, in Progress in Mathematics {\bf 92},
Birkh\"auser, Boston 1990, p.471--506.


\bibitem[GK]{GK} I.~M.~Gelfand, A.~A.~Kirillov {\em The structure of
the Lie field connected with a split semisimple Lie algebra},
Funct. Analysis Appl. No. 1 (1969) p. 6--21.

\bibitem[G1]{ghost} M.~Gorelik, 
{\em On the ghost centre of Lie superalgebras}, 
Les Annales de l'institut Fourier, t.~50 (2000), p.1745--1764.


\bibitem[G2]{gcent} M.~Gorelik, {\em Kac construction of the centre $\Zg$},
Journal of non-linear Math. Phys., {\bf 11} (2004) p. 325--349.





\bibitem[Ja]{ja} J.~C.~Jantzen, {\em Moduln mit einem h\"ochsten Gewicht},
Lecture Notes in Math., {\bf 750} (1979).

\bibitem[J]{jbook} A.~Joseph, {\em Quantum groups and their primitive ideals},
Springer (1995).

\bibitem[K1]{kadv} V.~Kac, {\em Lie superalgebras}, Adv. in Math. {\bf 26}
(1977), p.8--96.



\bibitem[K2]{KFF}
V.~Kac, {\em Contravariant form on Lie algebras and superalgebras}, 
Lect. Notes in Physics, Springer, Berlin e.a, (1979), 441--445.


\bibitem[K3]{klapl} V.~Kac, {\em Laplace operators of infinite-dimensional
Lie algebras and theta functions}, Proc. Natl. Acad. Sci. USA {\bf 81}
(1984), p. 645--647.


\bibitem[K4]{k13} 
V.~Kac, {\em Highest weight representations of conformal current
algebras}, in Topological and geometrical methods in field theory,
Editors Hietarinta J. et. al., Espoo, (1986)
World Scientific Publishing, Teaneck-NJ, p. 3--15.


\bibitem[KK]{kk} V.~Kac, D.~Kazhdan, {\em 
Structure of representations with highest weight of infinite dimensional
Lie algebras},  Adv. in Math. {\bf 34} (1979), p. 97--108. 


\bibitem[Kn]{kn} M.-A.~Knus, {\em Quadratic and hermitian forms over rings},
Springer (1991).

\bibitem[LM]{lm} E.~S.~Letzter, I.~M.~Musson, {\em Complete sets of 
representations of classical Lie superalgebras}, Lett. Math. Phys.
{\bf 31} (1994), p. 247--253.

\bibitem[ML]{mc} S.~Mac Lane {\em Categories for the working mathematician},
Springer, (1971).

\bibitem[NS]{ns} M.~Nazarov, A.~Sergeev, {\em Centralizer construction
of the Yangian of the queer Lie superalgebra}, arXiv:math.RT/0404086,
to appear in Studies in Lie Theory, A. Joseph Festschrift, Progress in Math.
{\bf 243}.

\bibitem[P]{pe} I.~Penkov, {\em 
Characters of typical irreducible finite-dimensional
$\fq(n)$-modules}, Functional Anal. Appl. {\bf 20} (1986) no.1, p. 30--37.


\bibitem[R]{ro} L.~Rowen, {\em Ring theory}, vol. I, Pure and Applied Math., 
Academic Press, Boston, (1988).

\bibitem[S1]{s1} A.~Sergeev,  {\em Invariant polynomial functions on
Lie superalgebras},  C.R. Acad. Bulgare Sci.
{\bf 35} (1982),  no. 5, p. 573--576.

\bibitem[S2]{sq} A.~Sergeev, {\em The centre of enveloping algebra
for Lie superalgebra $Q(n,\mathbb{C})$}, Lett. Math. Phys. {\bf 7} (1983)
no.3, p. 177--179.


\bibitem[S3]{s3} A.~Sergeev,  {\em Invariant polynomials on simple
Lie superalgebras}, Representation Theory, {\bf 3}, (1999), p. 250--280.

\bibitem[Sh]{sh} N.~Shapovalov, {\em On a bilinear form on the
universal enveloping algebra of a complex semisimple Lie algebra},
 Functional Anal. Appl. {\bf 6} (1972), p. 307--312 (in Russian).
\end{thebibliography}
\end{document}